\documentclass[12pt,a4paper,leqno]{article}
\usepackage{amssymb,amsmath,amsthm}
\usepackage{latexsym}
\usepackage[mathscr]{eucal}
\usepackage{color}
\usepackage{comment}

\newtheorem{theorem}{Theorem}
\newtheorem{lemma}{Lemma}
\newtheorem{proposition}{Proposition}
\newtheorem{corollary}{Corollary}
\newtheorem{remark}{Remark}
\numberwithin{equation}{section}
\numberwithin{theorem}{section}
\numberwithin{lemma}{section}
\numberwithin{proposition}{section}
\numberwithin{corollary}{section}
\numberwithin{remark}{section}

\begin{document}
\title{Large time decay of the fluid-structure semigroup subject to the
Navier slip-with-friction boundary condition}
\author{Toshiaki Hishida and Bum Ja Jin}
\date{}
\maketitle
\begin{abstract}
Consider the motion of a rigid body moving into a viscous incompressible fluid that fills the whole space in 3D. 
The fluid motion obeys the Navier-Stokes system, while the motion of the body is governed by the Newton law.
Both motions are affected by each other at the fluid-solid interface through the Navier slip-with-friction 
boundary condition, which is known as a substitution for the no-slip one to resolve the no-collision paradox. 
We adopt a monolithic approach and carry out the spectral analysis to deduce the large time behavior of the
fluid-structure semigroup, where the shape of the body is allowed to be arbitrary. What is remarkable is that 
the decay structure is somewhat different under the full slip condition when the body is axially symmetric.
If, in particular, the body is a ball, we then apply the linear theory developed above to the nonlinear initial value
problem to show the stability of the rest state. It is, furthermore, asymptotically stable except for the full slip case.

\noindent
MSC: 35Q35, 76D05

\noindent
Keywords: Fluid-structure interaction, Navier boundary condition, Decay estimate. 
\end{abstract}

\section{Introduction}
\label{intro}

This paper studies
the motion, especially the large time behavior, of a single rigid body immersed in a viscous incompressible fluid 
that fills the whole space $\mathbb R^3$.
The fluid motion obeys the Navier-Stokes system in the exterior of the moving body, 
while the motion of the body is governed by the balance for linear and angular momentum, involving
the force and torque exerted by the fluid.
Both motions are affected by each other at the boundary through the Navier slip-with-friction boundary condition.
By taking a frame attached to the body with its center of mass being fixed at the origin,
the new unknowns $u(x,t)\in \mathbb R^3$, $p(x,t)\in\mathbb R$,
$\eta(t)\in\mathbb R^3$ and $\omega(t)\in\mathbb R^3$ should obey the following set of equations, see Galdi \cite{Ga02},
where $u$ and $p$ denote the velocity and pressure of the fluid, while $\eta$ and $\omega$ stand for
the translational and angular velocities of the body:
\begin{equation}
\partial_tu+(u-\eta-\omega\times x)\cdot\nabla u+\omega\times u=\Delta u-\nabla p, \quad \mbox{div $u$}=0 \quad
\mbox{in $\Omega\times (0,\infty)$},
\label{ns}
\end{equation}
\begin{equation}
\nu\cdot (u-\eta-\omega\times x)|_{\partial\Omega}=0, \quad 
[(2\mathbb Du)\nu]_\tau+\alpha (u-\eta-\omega\times x)_\tau\big|_{\partial\Omega}=0, \quad
\lim_{|x|\to\infty}u=0,
\label{bc}
\end{equation}
\begin{equation}
m\frac{d\eta}{dt}+\omega\times (m\eta)+\int_{\partial\Omega}{\mathbb S}(u,p)\nu\,d\sigma=0, \quad
J\frac{d\omega}{dt}+\omega\times (J\omega)+\int_{\partial\Omega}x\times {\mathbb S}(u,p)\nu\,d\sigma=0,
\label{newton}
\end{equation}
being completed with the initial condition at $t=0$.
The aforementioned change of frame involves a proper orthogonal matrix relating to rotation of the body, whereas one can fix the frame
attached to the center of mass of the rigid body solely by translation if, in particular, the body is a ball.
For this latter case, the equations \eqref{ns} and \eqref{newton} are, respectively, simplified as
\begin{equation}
\partial_tu+(u-\eta)\cdot\nabla u=\Delta u-\nabla p, \quad \mbox{div $u$}=0 \quad
\mbox{in $\Omega\times (0,\infty)$},
\label{ns2}
\end{equation}
\begin{equation}
m\frac{d\eta}{dt}+\int_{\partial\Omega}{\mathbb S}(u,p)\nu\,d\sigma=0, \quad
J\frac{d\omega}{dt}+\int_{\partial\Omega}x\times {\mathbb S}(u,p)\nu\,d\sigma=0.
\label{newton2}
\end{equation}
Regardless of the shape of the body, the corresponding
linearized system of equations consists of 
\begin{equation}
\partial_tu=\Delta u-\nabla p, \quad \mbox{div $u$}=0 \quad
\mbox{in $\Omega\times (0,\infty)$}
\label{stokes}
\end{equation}
and \eqref{newton2} subject to the boundary condition \eqref{bc}. 

Here, $\Omega$ is an exterior domain in $\mathbb R^3$ with smooth boundary $\partial\Omega\in C^{2,1}$ such that
${\mathscr B}=\mathbb R^3\setminus \overline{\Omega}$ is connected as well as bounded and 
$\int_{\mathscr B}x\,dx=0$.
In \eqref{bc}, 
we have denoted by $\nu$ the unit normal to $\partial\Omega$
directed toward the body ${\mathscr B}$, and by $\alpha=\alpha(x)\geq 0$ the friction coefficient at $x\in\partial\Omega$.
Moreover,
\[
\mathbb Du=\frac{\nabla u+(\nabla u)^\top}{2} 
\]
is the symmetric gradient (deformation tensor), 
where $(\cdot)^\top$ denotes
the transpose, and 
\[
\mathbb S(u,p)=2\mathbb Du-p\mathbb I
\]
stands for the Cauchy stress tensor with $\mathbb I$ being the $3\times 3$ unity matrix, so that 
$\mbox{div $\mathbb S(u,p)$}=\Delta u-\nabla p$ on account of the solenoidal condition.
Given vector field $W$ at the boundary $\partial\Omega$, we denote by $W_\tau$ the tangential component, which is defined by
\[
W_\tau=W-(\nu\cdot W)\nu=(\nu\times W)\times\nu.
\]
The condition $\mbox{\eqref{bc}}_2$ arises from
\[
[\mathbb S(u,p)\nu]_\tau+\alpha (u-\eta-\omega\times x)_\tau\big|_{\partial\Omega}=0,
\]
nevertheless, it does not involve the pressure actually.
When $\alpha>0$ is a constant, this tells us that the tangential component of $(u-\eta-\omega\times x)$ is proportional to the one of
the normal stress exerted by the fluid with the proportional coefficient $1/\alpha$, which is called the slip length.
The particular case $\alpha\equiv 0$ is the full slip condition, while the no-slip one $u|_{\partial\Omega}=\eta+\omega\times x$,
that is often adopted in many papers, is recovered formally when $\alpha\to\infty$.
Finally, in \eqref{newton} and \eqref{newton2}, $m$ and $J$ are given by
\begin{equation}
m=\int_{\mathscr B}\rho\,dx, \qquad J=\int_{\mathscr B}(|x|^2\mathbb I-x\otimes x)\rho\,dx,
\label{inertia}
\end{equation}
denoting respectively the mass of the body and tensor of inertia that is positive as well as symmetric, 
where the constant $\rho>0$ is the density of the body, while the one of the fluid may be assumed to be $1$ without loss of
generality.

Within the context of fluid-structure interaction,
there are some regimes in which the Navier slip-with-friction boundary condition
(that we call the Navier condition, for short) $\mbox{\eqref{bc}}_{1,2}$ would be suitable rather than the no-slip one.
In fact, it was proved in \cite{Hill07} (for 2D case) and \cite{HT09} (for 3D case) by Hillairet and Takahashi that no collision
can occur between a rigid ball and a wall in finite time under the no-slip condition.
Later on, G\'erard-Varet, Hillairet and Wang \cite{GHW15} made it clear that
boundary conditions affect significantly the possibility of contact by showing that,
for any weak solution (in a reasonable sense) within an appropriate setting, a rigid ball inside a bounded container
filled with a fluid touches the boundary of the container in finite time provided that the Navier condition is imposed
on both boundaries, 
namely, the fluid-body interface and the outer fixed boundary.
Although the collision is not the issue in \eqref{ns}--\eqref{newton}, the Navier condition in our specific setting
would be worth analyzing to understand another feature of dynamics induced by this boundary condition.
So far, to our knowledge, the only literature dealing with the same setting in the whole of $\mathbb R^3$
from such a viewpoint seems to be
\cite{J25} by the second author of the present paper on the linearized
system \eqref{newton2}--\eqref{stokes} subject to \eqref{bc} as well as the initial condition.
It has been proved in \cite{J25} that the linearized system enjoys a maximal $L^q$-regularity and generates a 
bounded analytic semigroup, that is called the fluid-structure semigroup,
in the $L^q$-setting together with $L^q$-$L^r$ estimate (exept for $r=\infty$) globally in time when $\alpha>0$;
however, as for the decay estimate
of gradient, the rate is less sharp and does
not allow us to proceed to the nonlinear problem.

The objective of
the present paper 
is to find a different structure  
(for the full slip case $\alpha\equiv 0$)
from the no-slip one and
to improve the $L^q$-$L^r$ decay estimate of the fluid-structure semigroup with/without the friction
in a satisfactory form, see Theorem \ref{main}; that is, so far as the fluid velocity is concerned, it is of the
same form as the result on
the Stokes semigroup in exterior domains, which was studied very well in
the pioneering work \cite{I89} by Iwashita and, later on, in \cite{MS97} by Maremonti and Solonnikov
under the no-slip condition, see also Shimada and Yamaguchi \cite{SY08} as well as \cite{H26} by the first author
under the Navier one.
Once we establish such a desired linear theory, we can apply it to the nonlinear initial value problem
for \eqref{ns2}--\eqref{newton2} subject to \eqref{bc} provided that the rigid body is a ball,
as done for the no-slip case in \cite{EMT23} by Ervedoza, Maity and Tucsnak 
and in \cite{H24} by the first author. 
As a consequence for the rigid ball, it turns out that the rest state is indeed stable and, furthermore, asymptotically stable 
with definite rate of decay of disturbance except for the full slip case $\alpha\equiv 0$, see Theorem \ref{main-nonlinear}.
As long as we adopt the semigroup approach in $L^q$, however, the problem
\eqref{ns}--\eqref{newton} 
is still out of reach, like the no-slip case \cite{EMT23,H24},
because of presence of the nonlinear term $(\omega\times x)\cdot\nabla u$ with unbounded coefficient at spatial infinity.
In this general case under the no-slip condition, by means of the energy method, Galdi \cite{Ga23} succeeded in showing that
the rest state is asymptotically stable.
If the fluid region is replaced by a bounded domain, in which the Navier condition is imposed at the interface of
fluid-solid, while the no-slip one is adopted at the outer fixed boundary, we refer to
Al Baba, Ghosh, Muha and Necasova \cite{AGMN21} on the exponential 
stability of the rest state.

For the proof of $L^q$-$L^r$ decay estimates of solutions to the initial value problem for the linearized system
\eqref{newton2}--\eqref{stokes} subject to 
\eqref{bc}, the essential step is to deduce the local energy decay properties of the semigroup, with which at hand
one can then proceed to the stage of cut-off procedures to arrive at the desired estimates as in many
papers including \cite{ES04,ES05,EMT23,H04,H16,H20,H24,H26,HS09,I89,KS98,SY08}.
To this end, we adapt the spectral analysis to the monolithic velocity, see \eqref{mono-velo} below, as in \cite{EMT23}
by Ervedoza, Maity and Tucsnak
for the no-slip case, however, there are a couple of remarkable differences.

Firstly, especially for the full slip case $\alpha\equiv 0$, the decay structure depends on the shape of the rigid body, to be precise,
as observed in subsection \ref{axisym}, a part of the rigid motion cannot decay to zero as $t\to\infty$ when the body is
axially symmetric.
Thus, in this axisymmetric case, the corresponding part  
has to be subtracted from the whole resolvent, see \eqref{resol-rep2};
then, our task is to investigate the asymptotic behavior of the remaining part near $\lambda=0$ in the complex plane,
leading to the local energy decay estimate through the formula \eqref{sg-rep} of the semigroup
in terms of the resolvent.

From the technical aspect, secondly,
the asymptotic expansion near $\lambda=0$ of the Stokes resolvent in exterior domains 
due to Iwashita \cite{I89} is used in \cite{EMT23}, while we make full use of the regularity near $\lambda=0$ along the
imaginary axis of the Stokes resolvent in exterior domains 
obtained recently by the first author \cite{H26}
(on the Oseen semigroup subject to the Navier condition, that covers the Stokes one as a particular case).
The latter approach adopted here
was employed first by Kobayashi and Shibata \cite{KS98} for the Oseen semigroup in exterior domains under the no-slip condition
and, later on, 
developed by \cite{H04} even for the Stokes semigroup on a different unbounded domain.
Uniqueness of solutions to the resolvent system, especially for $\lambda=0$, plays a crucial role 
and the statement for the full slip case should be 
given carefully in a right form that depends on whether or not the body is axially symmetric, see Proposition \ref{uni}.

It would be worth while trying the other route
developed by the first author \cite{H18,H20,H24} for the non-autonomous system to obtain the local energy decay properties 
although the rate of decay is less sharp (nevertheless, the final results are the same).
The strategy is based on the energy relation, however, unlike the no-slip case \cite{H24},
the control of the rigid motion by $\mathbb Du$ is a nontrivial issue.
This will be overcome in section \ref{alter},
which is independent of the other parts of the present paper.

Let us also mention that there is the other way of transformation to reduce the original problem to the one
in a time-independent domain.
This is a local transformation that keeps the situation far from the body as it is
in order to avoid the nonlinear term $(\omega\times x)\cdot\nabla u$ in \eqref{ns} even if the shape 
of the body is arbitrary, see Cumsille and Takahashi \cite{CT08}
on the global existence of a strong solution for small data within the $L^2$-theory (however, without any large time behavior),
and Maity and Tucsnak \cite{MT23}
even on  the large time decay in the framework of maximal $L^q$-regularity.
Both works are concerned with the no-slip case.
Since the latter paper made use of the $L^q$-$L^r$ estimates of the fluid-structure semigroup as well,
the result of the present paper might be a basis toward analysis of 
the large time behavior of the body of arbitrary shape
under the Navier condition.

The paper is organized as follows.
In the next section, after introducing function spaces as well as the fluid-structure operator, we give an observation
about the full slip case and then our main theorems (Theorem \ref{main}, Theorem \ref{main-nonlinear}). 
In section \ref{pre}, for later use, we recall some results on the Stokes resolvent in exterior domains from \cite{H26}.
Section \ref{from-to} is devoted to the spectral analysis of the fluid-structure resolvent.
Local energy decay properties are deduced in section \ref{led}.
Section \ref{pr-main}
completes the proof of the $L^q$-$L^r$ estimates. 
The final section discusses the other route of the proof mentioned above.

\section{Results}
\label{result}

\subsection{Fluid-structure operator}
\label{semi}

Let us begin with introducing fundamental notation.
Throughout this paper, we assume that both ${\mathscr B}$ and $\Omega$
are domains 
in $\mathbb R^3$ such that
\begin{equation}
\overline{\mathscr B}\subset B_1, \qquad 
\Omega=\mathbb R^3\setminus\overline{\mathscr B}, \qquad
\partial\Omega=\partial{\mathscr B}\in C^{2,1}, \qquad
\int_{\mathscr B}x\,dx=0,
\label{dom}
\end{equation}
where $B_R=B_R(0)$ stands for the open ball centered at the origin with radius $R>0$.
For $R>1$, we set
\[
\Omega_R=\Omega\cap B_R.
\]
Given a domain $D\subset\mathbb R^3$, $q\in [1,\infty]$ and integer $k\geq 0$,
the standard Lebesgue and Sobolev spaces are denoted by
$L^q(D)$ and $W^{k,q}(D)$.
We abbreviate the norm $\|\cdot\|_{q,D}=\|\cdot\|_{L^q(D)}$.
By $\langle\cdot,\cdot\rangle_D$ 
(resp. $\langle \cdot,\cdot\rangle_{\partial D}$)
we denote standard duality pairings
over the domain $D$ (resp. $\partial D$) 
in each context.
The class $C_0^\infty(D)$ consists of all $C^\infty$-functions with
compact support in $D$, then $W^{k,q}_0(D)$ denotes the completion
of $C_0^\infty(D)$ in $W^{k,q}(D)$, where $k>0$ is an integer.
Given $q\in [1,\infty]$, let $q^\prime\in [1,\infty]$ be the conjugate number
defined by $1/q^\prime+1/q=1$.
For $q\in (1,\infty)$, we define the Sobolev space of order $(-1)$ by
$W^{-1,q}(D)=W^{1,q^\prime}_0(D)^*$.
In what follows we adopt the same symbols for denoting scalar and vector (even tensor)
functions as long as there is no confusion. 
By $C^\infty_{0,\sigma}(D)$ we denote the class of all 
vector fields $u$
which are in $C_0^\infty(D)$ and satisfy $\mbox{div $u$}=0$ in $D$. 
Let $X$ and $Y$ be Banach spaces.
Then ${\mathcal L}(X,Y)$ stands for the Banach space consisting of
all bounded linear operators from $X$ to $Y$.
Finally, we denote several positive constants by $C$,
which may change from line to line.

We now introduce several function spaces especially for fluid-structure interaction problems and define
the fluid-structure operator to formulate our problem within the monolithic framework.
Let $1<q<\infty$.
The underlying space we adopt is
\begin{equation}
X_q(\mathbb R^3)=\{U\in L^q_{\rm R}(\mathbb R^3);\, \mbox{div $U=0$ in $\mathbb R^3$}\},
\label{underlying}
\end{equation}
where
\begin{equation}
L^q_{\rm R}(\mathbb R^3)=\{U\in L^q(\mathbb R^3)^3;\, U|_{\mathscr B}\in {\rm RM}\}, \qquad
{\rm RM}=\{\eta+\omega\times x;\, (\eta,\omega)\in\mathbb R^3\times\mathbb R^3\}.
\label{fra-sp}
\end{equation}
For the resolvent problem, we have to consider the complex rigid motions $\eta+\omega\times x$ with 
$(\eta,\omega)\in \mathbb C^3\times \mathbb C^3$.
In such a context, readers are asked to understand statements with obvious changes.
It is convenient to define the map
\[
i: L^q_{\rm R}(\mathbb R^3)\to L^q(\Omega)^3\times \mathbb R^3\times \mathbb R^3
\]
by
\begin{equation}
\begin{split}
&i(U):=(u,\eta,\omega)\;\; \mbox{with}  \\
&u=U|_\Omega, \qquad \eta=\frac{1}{m}\int_{\mathscr B}U(x)\rho\,dx, \qquad
\omega=J^{-1}\int_{\mathscr B}x\times U(x)\rho\,dx,
\end{split}
\label{mono}
\end{equation}
where $\rho,\, m$ and $J$ are given by \eqref{inertia}, 
then we see that
\begin{equation}
U=u\chi_\Omega+(\eta+\omega\times x)\chi_{\mathscr B}.
\label{mono-velo}
\end{equation}
Let $U\in L^q_{\rm R}(\mathbb R^3)$, $(u,\eta,\omega)=i(U)$ and 
$\mbox{div $u$}=0$ in $\Omega$.
Then $U\in X_q(\mathbb R^3)$ if and only if 
\begin{equation}
\nu\cdot (u-\eta-\omega\times x)|_{\partial\Omega}=0,  
\label{normal-va}
\end{equation}
with $\nu$ being the unit normal to $\partial\Omega$ directed toward ${\mathscr B}$.
Elements of $X_q(\mathbb R^3)$ are understood through \eqref{mono}--\eqref{normal-va} and may be called 
monolithic velocities.
The space $X_q(\mathbb R^3)$ is a Banach space endowed with norm
\[
\|U\|_{X_q(\mathbb R^3)}=\Big(\|u\|_{q,\Omega}^q+\|(\eta+\omega\times x)\rho^{1/q}\|_{q,{\mathscr B}}^q\Big)^{1/q},
\]
where $(u,\eta,\omega)=i(U)$.
Here, the constant weight $\rho$ is involved in the integral over the body ${\mathscr B}$.
When $q=2$, this is particularly relevant to describe the energy
\begin{equation}
\|U\|_{X_2(\mathbb R^3)}^2
=\int_\Omega|u(x)|^2\,dx+\int_{\mathscr B}|\eta+\omega\times x|^2\rho\,dx
=\|u\|_{2,\Omega}^2+m|\eta|^2+(J\omega)\cdot\omega
\label{ener}
\end{equation}
although it is obvious that the following three quantities are equivalent:
$\|U\|_{X_q(\mathbb R^3)}\sim \|U\|_{q,\mathbb R^3}\sim \|u\|_{q,\Omega}+|\eta|+|\omega|$.
The latter equality of \eqref{ener} is verified from the relation
\[
x\times (a\times x)=\big(|x|^2\mathbb I-x\otimes x\big)a\qquad\forall\,a\in\mathbb R^3
\]
as well as $\int_{\mathscr B}x\,dx=0$ and \eqref{inertia}.

As in the classical setting to deal with the Stokes and Navier-Stokes system, it is important to find out a right
space in which the pressure associated with the monolithic velocity lies.
Regarding this issue, Silvestre \cite{Sil02-1,Sil02-2} first proved the decomposition \eqref{decompo} below 
of the space $L^2_{\rm R}(\mathbb R^3)$ when $q=2$.
Later on, in \cite{H24} by the first author of the present paper, her result has been extended to the general $L^q$-setting.
Let us summarize the results for later use.
Note that the denseness of the class \eqref{nice-cl} below in $X_q(\mathbb R^3)$
was shown first by \cite[Appendix A.2]{EMT23}, while an alternative proof
was given by \cite{H24} within the context of the decomposition \eqref{decompo}.
\begin{lemma}
[{\cite[Proposition 3.1]{H24}}, {\cite[Theorem 3.2]{Sil02-1}}]
Let $q\in (1,\infty)$, and let
$X_q(\mathbb R^3)$ be the space given by \eqref{underlying}.

\begin{enumerate}

\item
The class
\begin{equation}
{\mathcal E}(\mathbb R^3)=\{U\in C_{0,\sigma}^\infty(\mathbb R^3);\, \mathbb DU=O\;\mbox{in ${\mathscr B}$}\}
\label{nice-cl}
\end{equation}
is dense in $X_q(\mathbb R^3)$.

\item
The space $L^q_{\rm R}(\mathbb R^3)$ given by \eqref{fra-sp}
admits the decomposition
\begin{equation}
L^q_{\rm R}(\mathbb R^3)=X_q(\mathbb R^3)\oplus Z_q(\mathbb R^3)
\label{decompo}
\end{equation}
with
\begin{equation}
\begin{split}
Z_q(\mathbb R^3)=\Big\{V\in L^q_{\rm R}(\mathbb R^3);\,
&V|_\Omega=\nabla p,\, p\in L^q_{\rm loc}(\overline{\Omega}),\, 
V|_{\mathscr B}=\eta+\omega\times x,  \\
&\eta=\frac{-1}{m}\int_{\partial\Omega}p\nu\,d\sigma,\,
\omega=-J^{-1}\int_{\partial\Omega}x\times (p\nu)\,d\sigma\Big\}.
\end{split}
\label{anni}
\end{equation}

\item
We have the relations
\[
X_q(\mathbb R^3)^*=Z_q(\mathbb R^3)^\perp=X_{q^\prime}(\mathbb R^3), \qquad
Z_q(\mathbb R^3)^*=X_q(\mathbb R^3)^\perp=Z_{q^\prime}(\mathbb R^3),
\]
where $Z_q(\mathbb R^3)^\perp$ and $X_q(\mathbb R^3)^\perp$ denote the annihilators
with respect to $\langle \cdot,\cdot\rangle_{\mathbb R^3,\rho}$. 
Here, this pairing is given by
\begin{equation}
\begin{split}
\langle U,V\rangle_{\mathbb R^3,\rho}
&=\int_\Omega u\cdot v\,dx+\int_{\mathscr B}(\eta_u+\omega_u\times x)\cdot (\eta_v+\omega_v\times x)\rho\,dx  \\
&=\langle u,v\rangle_\Omega+m\eta_u\cdot\eta_v+(J\omega_u)\cdot \omega_v
\end{split}
\label{pair}
\end{equation}
for $U\in L^q_{\rm R}(\mathbb R^3)$, $V\in L^{q^\prime}_{\rm R}(\mathbb R^3)$, $i(U)=(u,\eta_u,\omega_u)$ and
$i(V)=(v,\eta_v,\omega_v)$.

\item
Let $\mathbb P=\mathbb P_q$ be the projection from $L^q_{\rm R}(\mathbb R^3)$ onto $X_q(\mathbb R^3)$
associated with \eqref{decompo}, then 
\begin{equation}
\|\mathbb PU\|_{q,\mathbb R^3}\leq C\|U\|_{q,\mathbb R^3}
\label{proj-bdd}
\end{equation}
\begin{equation}
\langle\mathbb P_qU,V\rangle_{\mathbb R^3,\rho}=\langle U,\mathbb P_{q^\prime}V\rangle_{\mathbb R^3,\rho}
\label{sym-P}
\end{equation}
for all $U\in L^q_{\rm R}(\mathbb R^3)$ and $V\in L^{q^\prime}_{\rm R}(\mathbb R^3)$
with some constant $C>0$.
\end{enumerate}
\label{lem-decom}
\end{lemma}

There is the other way due to Wang and Xin \cite[Theorem 2.2]{WX11}, see also Bae and Jin \cite[Lemma 2.1]{BJ23},
to eliminate the pressure.
The space they decompose is not $L^q_{\rm R}(\mathbb R^3)$ but $L^q(\mathbb R^3)^3$.
By $\mathbb Q_q: L^q(\mathbb R^3)^3\to X_q(\mathbb R^3)$ we denote the projection associated with their decomposition.
Then we know $\mathbb Q_q^*=\mathbb Q_{q^\prime}$ with respect to the standard pairing $\langle\cdot,\cdot\rangle_{\mathbb R^3}$,
whereas $\mathbb P_q^*=\mathbb P_{q^\prime}$ with respect to $\langle\cdot,\cdot\rangle_{\mathbb R^3,\rho}$ as seen from
\eqref{sym-P}.
Thus, the decomposition \eqref{decompo} is 
more convenient if using 
a duality argument and it is  
preferred in the present paper.
Note that, if the density of the fluid is assumed to be $1$, then $\rho$ is regarded as the ratio of both densities of the fluid-structure
and, therefore, one cannot assume $\rho=1$ without loss, see \cite[subsection 2.1]{H24}.

Let $\alpha\in C(\partial\Omega)$.
The fluid-structure operator subject to the Navier slip-with-friction boundary condition is 
defined by
\begin{equation}
\begin{split}
D_q(A)
&=\{U\in X_q(\mathbb R^3);\, 
u\in W^{2,q}(\Omega), 
[(2\mathbb Du)\nu]_\tau+\alpha(u-\eta-\omega\times x)_\tau|_{\partial\Omega}=0\}  \\
&\quad\mbox{with $(u,\eta,\omega)=i(U)$}, \\
AU
&=\mathbb P{\mathscr A}U, \\
{\mathscr A}U
&=
(-\Delta u)\chi_\Omega
+\left[\frac{1}{m}\int_{\partial\Omega}(2\mathbb Du)\nu\,d\sigma+\left(J^{-1}\int_{\partial\Omega}y\times (2\mathbb Du)\nu\,d\sigma_y\right)\times x\right]\chi_{\mathscr B},
\end{split}
\label{ig}
\end{equation}
where $\mathbb P$ 
is the projection introduced in Lemma \ref{lem-decom}. 
When computing $\langle AU,V\rangle_{\mathbb R^3,\rho}$ with $V\in X_{q^\prime}(\mathbb R^3)$ having further conditions,
after using \eqref{sym-P}, we write $\Delta u=\mbox{div $(2\mathbb Du)$}$, perform integration by parts and take into account
the boundary condition, see \eqref{sym-A} below.
Note that $D_q(A)$ is not included in $W^{1,q}(\mathbb R^3)$ unlike the no-slip case 
$u|_{\partial\Omega}=\eta+\omega\times x$.

As in \cite{J25} by the second author, it is also possible to introduce the fluid-structure operator, $\mathbb A$, by \eqref{ig}
with the projection $\mathbb P$
replaced by the other one $\mathbb Q$ mentioned above. 
As we would expect, those operators $A$ and $\mathbb A$ coincide with each other on account of uniqueness
of solutions to the resolvent system \eqref{FS-1}--\eqref{FS-2} below with $\lambda =1$, 
which is equivalent to 
$(1+A)U=F$ as well as to $(1+\mathbb A)U=F$ in $X_q(\mathbb R^3)$.
The argument was found in \cite[Proposition 3.3]{H24} by the first author
for the no-slip case $u|_{\partial\Omega}=\eta+\omega\times x$,
but it works for the Navier case as well since we have Proposition \ref{uni} on the aforementioned uniqueness of solutions.
Hence, several results for $\mathbb A$ established by \cite{J25} 
continue to hold for $A$;
in particular, we have the a priori estimates 
(\cite[Lemma 4.7]{J25})
\begin{equation}
\begin{split}
&\|u\|_{W^{2,q}(\Omega)}\leq C\big(\|AU\|_{q,\mathbb R^3}+\|U\|_{q,\mathbb R^3}\big),  \\
&\|u\|_{W^{1,q}(\Omega)}\leq C\big(\|AU\|_{q,\mathbb R^3}^{1/2}\|U\|_{q,\mathbb R^3}^{1/2}+\|U\|_{q,\mathbb R^3}\big),
\end{split}
\label{ellip-0}
\end{equation}
for all $U\in D_q(A)$, 
where $u=U|_\Omega$,
and know that
the operator $-A$ generates a bounded analytic semigroup $\{e^{-tA}\}_{t\geq 0}$ on the space $X_q(\mathbb R^3)$
for every $q\in (1,\infty)$:
\begin{equation}
\|\partial_t^ke^{-tA}F\|_{q,\mathbb R^3}=\|A^ke^{-tA}F\|_{q,\mathbb R^3}
\leq Ct^{-k}\|F\|_{q,\mathbb R^3}
\label{bdd-anal}
\end{equation}
for all $t>0$, $F\in X_q(\mathbb R^3)$ and every integer $k\geq 0$.
The proof of the latter given by \cite{J25} was based on the maximal $L^q$-regularity estimate
as in Maremonti and Solonnikov \cite[Section 4]{MS97}.
It is also possible to give an alternative proof
of generation of an analytic semigroup as long as $\alpha\in C^1(\partial\Omega)$.
In fact, we have only to follow
Maity and Tucsnak \cite{MT18}, see also \cite[Section 3.3]{H24}, by showing the desired resolvent estimate
for large $|\lambda|$ with the aid of 
results on the exterior problem
due to Shibata and Shimada \cite{SS07},
see section \ref{pre} of the present paper.
In \cite{J25},
a contradiction argument as in Borchers and Sohr \cite{BS87}, see also
Ervedoza, Maity and Tucsnak \cite[Section 6]{EMT23},
is adapted for deduction of
the resolvent estimate near $\lambda=0$, yielding 
the uniformly boundedness of the analytic semigroup
when $\alpha>0$.
The argument still works even for the case $\alpha\equiv 0$.
Regarding this issue, a few comments
are in order.
Firstly, uniqueness of solutions to the resolvent system plays a crucial
role and this will be discussed
in Proposition \ref{uni} in a different context.
Secondly, one needs a priori estimates near $\lambda=0$ of the
resolvent system in a bounded container subject to the Navier boundary condition
(at the interface of fluid-body, while the no-slip one is adopted at the outer fixed boundary), 
which is due to Al Baba, Ghosh, Muha and Necasova \cite[Section 4]{AGMN21}.
One has to extend their result to the case when non-solenoidal condition is imposed.
As in the existing literature, being traced back to Takahashi and Tucsnak \cite{TT04},
we refer to $e^{-tA}$ obtained above as the fluid-structure semigroup
(although it is called the Stokes-structure semigroup in \cite{H24} since the Oseen-structure evolution operator
is studied there).

If we denote by $A_q$ the operator \eqref{ig} with $D_q(A)$ acting on $X_q(\mathbb R^3)$, we have
$A_q^*=A_{q^\prime}$ with respect to $\langle\cdot,\cdot\rangle_{\mathbb R^3,\rho}$
for every $q\in (1,\infty)$, which follows from the relation (\cite[Lemma 4.1]{J25})
\begin{equation}
\langle AU,V\rangle_{\mathbb R^3,\rho}=\langle U,AV\rangle_{\mathbb R^3,\rho}
=2\langle\mathbb Du,\mathbb Dv\rangle_\Omega
+\int_{\partial\Omega}\alpha(y)(\eta_u+\omega_u\times y)\cdot (\eta_v+\omega_v\times y)\,d\sigma_y
\label{sym-A}
\end{equation}
for all $U\in D_q(A)$ and $V\in D_{q^\prime}(A)$, where $(u,\eta_u,\omega_u)=i(U), (v,\eta_v,\omega_v)=i(V)$
together with the surjectivity of $\lambda +A$ with $\lambda>0$.
Hence, the semigroup $e^{-tA}$ possesses the same symmetry for each $t>0$:
\begin{equation}
\langle e^{-tA}F,G\rangle_{\mathbb R^3,\rho}=\langle F,e^{-tA}G\rangle_{\mathbb R^3,\rho}
\label{sym-semi}
\end{equation}
for all $F\in X_q(\mathbb R^3)$ and $G\in X_{q^\prime}(\mathbb R^3)$.

\subsection{Axisymmetric body case subject to the full slip condition}
\label{axisym}

We are going to show that a part of unknowns of the linearized system \eqref{newton2}--\eqref{stokes}
subject to the full slip condition \eqref{bc} with $\alpha\equiv 0$ is decoupled with 
the others, so that the number of those is essentially reduced, when the rigid body is axially symmetric.
Note that, in this case, the axis of symmetry should pass through
the origin 
since $\int_{\mathscr B}x\,dx=0$.

In this subsection we are always concerned with the full slip boundary condition
\begin{equation}
\nu\cdot (u-\eta-\omega\times x)|_{\partial\Omega}=0, \qquad
[(2\mathbb Du)\nu]_\tau\big|_{\partial\Omega}=0,
\label{full-slip}
\end{equation}
from which it follows that
\begin{equation}
(2\mathbb Du)\nu
=\big[\nu\cdot\{(2\mathbb Du)\nu\}\big]\nu,
\label{sym-defo}
\end{equation}
and, thereby,
\begin{equation}
\mathbb S(u,p)\nu=\big[\nu\cdot\{(2\mathbb Du)\nu\}-p\big]\nu
\label{sym-0}
\end{equation}
is purely normal
at $\partial\Omega$.
Let $\gamma\in\mathbb R^3$ be a unit vector: $|\gamma|=1$.
In what follows we often use the decomposition
\[
a=(\gamma\cdot a)\gamma+(\gamma\times a)\times \gamma\qquad\forall\,a\in\mathbb R^3.
\]
Suppose that the rigid body is axisymmetric about the axis $\gamma\mathbb R$.
Under this condition, the outward unit normal $\nu=\nu(x)$ to the boundary $\partial\Omega$ fulfills
\begin{equation}
\gamma\cdot (x\times \nu)=\nu\cdot (\gamma\times x)=0 \qquad\forall\, x\in\partial\Omega,
\label{sym-1}
\end{equation}
regardless of the boundary condition.
From \eqref{sym-0} together with \eqref{sym-1} we find
\begin{equation}
\gamma\cdot
\big[x\times\mathbb S(u,p)\nu\big]\Big|_{\partial\Omega}=0.
\label{sym-2}
\end{equation}
As to the matrix $J$ given by \eqref{inertia} with ${\mathscr B}$
being axisymmetric about $\gamma\mathbb R$, there is a constant $c>0$ such that
\begin{equation}
\gamma\cdot (Ja)=c(\gamma\cdot a)\qquad\forall\,a\in\mathbb R^3.
\label{sym-inertia}
\end{equation}
In fact, with use of the orthogonal matrix $Q$ satisfying $Qe_3=\gamma$, we observe
\[
\gamma\cdot (Ja)=\int_{{\mathscr B}^\prime}
\big[|y|^2(\gamma\cdot a)-\{\gamma\cdot (Qy)\}\{a\cdot (Qy)\}\big]\rho\,dy
\]
where $Q{\mathscr B}^\prime={\mathscr B}$.
Since
\[
\gamma\cdot (Qy)=e_3\cdot y=y_3
\]
and since
\[
a\cdot (Qy)=
a\cdot\big[\{\gamma\cdot (Qy)\}\gamma+\{\gamma\times (Qy)\}\times \gamma
\big]
=y_3(a\cdot\gamma)+a\cdot \big[Q\{(e_3\times y)\times e_3\}\big],
\]
we furnish
\[
\gamma\cdot (Ja)
=\int_{{\mathscr B}^\prime}(|y|^2-y_3^2)\rho\,dy\; (\gamma\cdot a)
-\int_{{\mathscr B}^\prime}y_3\{(e_3\times y)\times e_3\}\rho\,dy\cdot (Q^\top a),
\]
in which the second integral is identically zero.
We thus obtain \eqref{sym-inertia}.
Similarly, it is also verified that $J\gamma=c\gamma$ with the same constant $c>0$ as in \eqref{sym-inertia}.

We now dot-multiply the equation $\mbox{\eqref{newton2}}_2$
by $\gamma$
and take into account \eqref{sym-2}--\eqref{sym-inertia} to get
\[
\frac{d}{dt}(\gamma\cdot\omega)=0,
\]
yielding 
\begin{equation}
\gamma\cdot \omega(t)=\gamma\cdot\omega(0).
\label{const-ii}
\end{equation}
In this way, our problem is reduced to \eqref{newton2}--\eqref{stokes} subject to \eqref{full-slip},
in which the equation $\mbox{\eqref{newton2}}_2$ for $\omega$
and $\mbox{\eqref{full-slip}}_1$ are replaced by
\begin{equation}
J\frac{d}{dt}(\gamma\times \omega)\times\gamma+
\left(\gamma\times \int_{\partial\Omega}x\times \mathbb S(u,p)\nu\,d\sigma\right)\times\gamma=0, 
\label{reduced-axisym}
\end{equation}
and by
\[
\nu\cdot \big[u-\eta-\{(\gamma\times\omega)\times\gamma\}
\times x\big]\Big|_{\partial\Omega}
=0,
\]
respectively, where the constant axial part
$(\gamma\cdot \omega(0))\gamma$ does not appear in this latter on account of 
\eqref{sym-1}.

If the body possesses two distinct axes of axial symmetry,
then it should be a ball, 
for which
$\nu$ is parallel to $(-x)$ at $\partial\Omega$,
see Falocchi and Gazzola \cite[Lemma 1]{FG22}.
Then one immediately see from \eqref{sym-0} that
\begin{equation}
x\times \mathbb S(u,p)\nu|_{\partial\Omega}=0.
\label{ball}
\end{equation}
Thus the equation of the balance for angular momentum leads us to
\begin{equation}
\omega(t)=\omega(0),
\label{const-i}
\end{equation}
so that our problem is reduced to \eqref{newton2}--\eqref{stokes} subject to \eqref{full-slip},
in which the equation $\mbox{\eqref{newton2}}_2$ for $\omega$ is
absent and 
$\mbox{\eqref{full-slip}}_1$ is replaced by
$\nu\cdot(u-\eta)|_{\partial\Omega}=0$.

\subsection{Main theorems}
\label{main-thm}

We are in a position to give our main theorem on the large time decay of the fluid-structure semigroup.
For the cases (i)--(iii) below in 
Theorem \ref{main}, 
given $F\in L^q_{\rm R}(\mathbb R^3)$ and $(f,\kappa,\mu)=i(F)$, see \eqref{mono}, it is convenient to introduce
\begin{equation}
\begin{split}
&\mbox{case (i)}\qquad
F_{\rm i}=f\chi_\Omega+\kappa\chi_{\mathscr B}, \\
&\mbox{case (ii)}\qquad
F_{\rm ii}=f\chi_\Omega+\big[\kappa+\{(\gamma\times\mu)\times \gamma\}\times x\big]\chi_{\mathscr B},
\end{split}
\label{IC-modi}
\end{equation}
and
\begin{equation}
\widetilde F=
\left\{
\begin{array}{ll}
F_{\rm i}&\mbox{case (i)}, \\
F_{\rm ii}\qquad&\mbox{case (ii)}, \\
F&\mbox{case (iii)},
\end{array}
\right.
\label{data-modi}
\end{equation} 
by taking into account the observation in the previous subsection.
Note that $\|\widetilde F\|_{q,\mathbb R^3}$ is equivalent to
\[
\left\{
\begin{array}{ll}
\|f\|_{q,\Omega}+|\kappa| &\mbox{case (i)}, \\
\|f\|_{q,\Omega}+|\kappa|+|\gamma\times \mu|\qquad &\mbox{case (ii)}, \\
\|f\|_{q,\Omega}+|\kappa|+|\mu|&\mbox{case (iii)},
\end{array}
\right.
\]
where $1<q<\infty$ and that $\|\widetilde F\|_{q,\mathbb R^3}\leq C\|F\|_{q,\mathbb R^3}$.
\begin{theorem}
Suppose that $\partial\Omega\in C^{2,1}$ and that $\alpha\in C^1(\partial\Omega)$ is nonnegative.
Let $q\in (1,\infty)$.
Consider
\begin{equation*}
\begin{split}
\mbox{\rm case (i)}\quad
&\mbox{$\alpha\equiv 0$, and ${\mathscr B}$ is a ball}; \\
\mbox{\rm case (ii)}\quad
&\mbox{$\alpha\equiv 0$, and ${\mathscr B}$ is axisymmetric about a single axis 
$\gamma\mathbb R$}; \\ 
\mbox{\rm case (iii)}\quad
&\mbox{otherwise},
\end{split}
\end{equation*}
where $\gamma\in \mathbb R^3$ with $|\gamma|=1$.

\begin{enumerate}
\item
For $F\in X_q(\mathbb R^3)$, 
we set
\[
U(t)=e^{-tA}F, \qquad 
(u,\eta,\omega)=i(U(t)),
\]
see \eqref{mono}.
Let $r\in [q,\infty]$.
There are constants $C=C(\alpha,q,r,\Omega)>0$ and $C^\prime=C^\prime(\alpha,q,\Omega)>0$ such that
\begin{equation}
\|u(t)\|_{r,\Omega}
\leq Ct^{-\frac{3}{2}(\frac{1}{q}-\frac{1}{r})}\|\widetilde F\|_{q,\mathbb R^3}
\label{decay-1}
\end{equation}
\begin{equation}
\left.
\begin{array}{rl}
\mbox{\rm case (i)}\quad
&|\eta(t)|  \\ 
\mbox{\rm case (ii)}\quad
&|\eta(t)|+|\gamma\times \omega(t)|  \\ 
\mbox{\rm case (iii)}\quad
&|\eta(t)|+|\omega(t)|
\end{array}
\right\}
\leq
C^\prime (1+t)^{-\frac{3}{2q}}\|\widetilde F\|_{q,\mathbb R^3}
\label{decay-body}
\end{equation}
for all $t>0$ and $F\in X_q(\mathbb R^3)$,
where $\widetilde F$ is given by \eqref{data-modi}.

Moreover, for every integer $k\geq 1$ and $r\in [q,\infty)$, there is a constant $C=C(\alpha,q,r,k,\Omega)>0$ such that
\begin{equation}
\|\partial_t^k e^{-tA}F\|_{r,\mathbb R^3}\leq Ct^{-k-\frac{3}{2}(\frac{1}{q}-\frac{1}{r})}\|\widetilde F\|_{q,\mathbb R^3}
\label{decay-tempo}
\end{equation}
for all $t>0$ and $F\in X_q(\mathbb R^3)$.

\item
Let $u(t)$ be as in the previous item.
Let $j\in \{1,2\}$ and
\begin{equation}
\left\{
\begin{array}{ll}
r\in [q,\infty] \qquad &\mbox{for $j=1$}, \\
r\in [q,\infty) &\mbox{for $j=2$}.
\end{array}
\right.
\label{exp-grad}
\end{equation}
Then there is a constant $C=C(\alpha,q,r,j,\Omega)>0$ such that
\begin{equation}
\|\nabla^j u(t)\|_{r,\Omega}
\leq Ct^{-\frac{j}{2}-\frac{3}{2}(\frac{1}{q}-\frac{1}{r})}
(1+t)^{\max\{\frac{j}{2}-\frac{3}{2r},0\}}\|\widetilde F\|_{q,\mathbb R^3}
\label{decay-2}
\end{equation}
for all $t>0$ and $F\in X_q(\mathbb R^3)$.

\item
Let $r\in [q,\infty]$.
Then there is a constant $C=C(\alpha,q,r,\Omega)>0$ such that
\begin{equation}
\|e^{-tA}\mathbb P\mbox{\rm div $G$}\|_{r,\mathbb R^3}\leq Ct^{-\frac{1}{2}-\frac{3}{2}(\frac{1}{q}-\frac{1}{r})}
(1+t)^{\max\{\frac{3}{2q}-1,0\}}
\|G\|_{q,\mathbb R^3}
\label{decay-3}
\end{equation}
for all $t>0$ and $G\in L^q(\mathbb R^3)^{3\times 3}$
with
$(G\nu)|_{\partial\Omega}=0$
as well as $\mbox{\rm div $G$}\in L^p_{\rm R}(\mathbb R^3)$ with some $p\in (1,\infty)$
satisfying $\widetilde{\mbox{\rm div $G$}}=\mbox{\rm div $G$}$, where $\widetilde{(\cdot)}$ is given by \eqref{data-modi}.

\item
Let $j\in \{1,2\}$, and 
let $r$ be as in \eqref{exp-grad}.
Then there is a constant $C=C(\alpha,q,r,j,\Omega)>0$ such that
\begin{equation}
\|\nabla^je^{-tA}\mathbb P\mbox{\rm div $G$}\|_{r,\Omega}
\leq Ct^{-\frac{1+j}{2}-\frac{3}{2}(\frac{1}{q}-\frac{1}{r})}
(1+t)^{\max\{\frac{j}{2}-\frac{3}{2r},0\}+\max\{\frac{3}{2q}-1,0\}}
\|G\|_{q,\mathbb R^3}
\label{decay-4}
\end{equation}
for all $t>0$ and $G$ with the same properties as in the previous item.
\end{enumerate}
\label{main}
\end{theorem}
\begin{remark}
The rate of decay \eqref{decay-2} for 
spatial derivatives is rephrased as
\begin{equation*}
\|\nabla u(t)\|_{r,\Omega}
\leq C\|\widetilde F\|_{q,\mathbb R^3}\left\{
\begin{array}{ll}
t^{-\frac{1}{2}-\frac{3}{2}(\frac{1}{q}-\frac{1}{r})}\qquad &(r\leq 3), \\
t^{-\frac{3}{2q}}&(r>3),
\end{array}
\right.
\end{equation*}
\begin{equation*}
\|\nabla^2 u(t)\|_{r,\Omega}
\leq C\|\widetilde F\|_{q,\mathbb R^3}\left\{
\begin{array}{ll}
t^{-1-\frac{3}{2}(\frac{1}{q}-\frac{1}{r})}\qquad &(r\leq \frac{3}{2}), \\
t^{-\frac{3}{2q}}&(r>\frac{3}{2}),
\end{array}
\right.
\end{equation*}
for all $t\geq 1$ and $F\in X_q(\mathbb R^3)$.
This rate is the same as the one for the Stokes semigroup in exterior domains 
subject to the no-slip or Navier condition
as shown in several papers, in particular, by
Iwashita \cite{I89} for the first time, and later on, by Maremonti and Solonnikov \cite{MS97}. 
For the gradient estimate of the Stokes semigroup, the rate is sharp 
as observed by \cite{MS97}, see also \cite{H11},
however, the optimality is not clear for the fluid-structure semigroup.
\label{rem-opti}
\end{remark}
\begin{remark}
Set $U(t)=e^{-tA}\mathbb P\mbox{\rm div $G$}$.
The assumptions on $G$ for 
\eqref{decay-3} imply that
$\widetilde{U(t)}=U(t)$, see 
\eqref{IC-modi}--\eqref{data-modi},
because of \eqref{obje-0} and \eqref{proj-ang2} below.
This is why the whole part decays to zero in \eqref{decay-3} even for the cases (i) and (ii).
Estimates \eqref{decay-3}--\eqref{decay-4}
are useful to study the nonlinear problem, in particular, to deal with the nonlinear term
$\eta\cdot\nabla u$, see \cite[(5.15)]{H24}.
\label{rem-div}
\end{remark}
\begin{remark}
As described in \cite[Remark 2.1, Remark 4.1]{H26} for the Oseen semigroup in exterior domains subject to the Navier boundary condition,
unfortunately, the dependence of the constants on the friction
coefficient $\alpha$ is not clear in \eqref{decay-1}--\eqref{decay-tempo} and \eqref{decay-2}--\eqref{decay-4}.
This is because one does not see such dependence of the constants  
in \eqref{ext-conti}--\eqref{ext-3}, see Proposition \ref{prop-ext} below.
Such dependence of the constant in \eqref{bdd-anal} is also unclear.
Suppose that $\alpha\geq 0$ is a constant.
Since we have the same result as in Theorem \ref{main} for the case (iii) when the no-slip condition is imposed, see \cite{EMT23,H24},
we might expect that those constants mentioned above would be taken uniformly for $\alpha\to\infty$.
\label{rem-dep-fric}
\end{remark}

Let us consider the nonlinear initial value problem for
\eqref{ns2}--\eqref{newton2} subject to \eqref{bc},
where $\mbox{\eqref{bc}}_1$ is rewritten as
$\nu\cdot (u-\eta)|_{\partial\Omega}=0$, as well as
\[
u(\cdot,0)=u_0, \qquad \eta(0)=\eta_0, \qquad \omega(0)=\omega_0,
\]
when the rigid body is a ball.
For the full slip case $\alpha\equiv 0$,
as in the linearized system, see the end of subsection \ref{axisym}, we see from \eqref{ball} that $\omega(t)=\omega_0$.
Set
\[
H(U)=\mathbb P\big[\{(\eta-u)\cdot\nabla u\}\chi_\Omega\big]
=\mathbb P\mbox{div $\{u\otimes (\eta-u)\chi_\Omega\}$},
\]
where the latter equality follows from $\nu\cdot (u-\eta)|_{\partial\Omega}=0$. 
Then, by use of the monolithic velocity \eqref{mono-velo}, the problem is formulated as
\begin{equation}
\frac{dU}{dt}+AU=H(U), \quad t>0; \qquad U(0)=U_0,
\label{evo-eq}
\end{equation}
where $U_0=u_0\chi_\Omega+(\eta_0+\omega_0\times x)\chi_{\mathscr B}$.

The following result on the nonlinear problem tells us that
the rest state is stable 
as in \eqref{stable} and, except for the full slip case, it is even asymptotically stable, 
see \eqref{NS-decay}--\eqref{ang-decay} below.
\begin{theorem}
Suppose that ${\mathscr B}$ is a ball and that $\alpha\in C^1(\partial\Omega)$ is nonnegative.
There is a constant $\delta=\delta(\alpha,\Omega)>0$ such that if $U_0\in X_3(\mathbb R^3)$ satisfies 
$\|\widetilde U_{0}\|_{3,\mathbb R^3}<\delta$, see \eqref{IC-modi}--\eqref{data-modi},
then problem \eqref{evo-eq} admits a unique strong solution
\[
U\in C([0,\infty);\,X_3(\mathbb R^3))\cap C((0,\infty);\, D_3(A))\cap C^1((0,\infty);\, X_3(\mathbb R^3))
\]
which enjoys
\begin{equation}
\left\{
\begin{array}{ll}
\|u(t)\|_{3,\Omega}+|\eta(t)|\leq C\|U_{0,{\rm i}}\|_{3,\mathbb R^3}, \quad 
\omega(t)=\omega_0\qquad&\mbox{if $\alpha\equiv 0$},  \\
\|u(t)\|_{3,\Omega}+|\eta(t)|+|\omega(t)|\leq C\|U_0\|_{3,\mathbb R^3}&\mbox{otherwise},
\end{array}
\right.
\label{stable}
\end{equation}
for all $t>0$ with some constant $C=C(\alpha,\Omega)>0$ and
\begin{equation}
\|u(t)\|_{q,\Omega}=o(t^{-1/2+3/2q}), \qquad
\|\nabla u(t)\|_{r,\Omega}+|\eta(t)|=o(t^{-1/2})
\label{NS-decay}
\end{equation}
as $t\to\infty$ for every $q,\, r\in [3,\infty]$, 
where $(u,\eta,\omega)=i(U(t))$, 
see \eqref{mono}.
Moreover, except for the full slip case $\alpha\equiv 0$, we have
\begin{equation}
|\omega(t)|=o(t^{-1/2})
\label{ang-decay}
\end{equation}
as well for $t\to\infty$.
\label{main-nonlinear}
\end{theorem}

With Theorem \ref{main} at hand, we just follow the same way as in the proof of \cite[Theorem 2.2]{H24}
to conclude Theorem \ref{main-nonlinear}.
Although we may omit the details, a few comments are in order. 

Firstly,
unlike the no-slip case studied in \cite{H24},
$\nabla e^{-tA}F$ belongs merely to $L^r(\Omega)$ as in \eqref{decay-2}, however, 
no longer belongs to any $L^r(\mathbb R^3)$; nonetheless, this is not a problem at all.

Secondly,
the reason why the decay of $\|\nabla u(t)\|_{\infty,\Omega}$ is also deduced in \eqref{NS-decay}
is that we have \eqref{decay-2} and \eqref{decay-4} even for $(j,r)=(1,\infty)$
(which was not available in \cite{H24}).
Note that one needs a slight modification of the proof of \cite[(5.8)]{H24} 
(in which $\|\nabla (\cdot)\|_{r,\mathbb R^3}$ has to be replaced by $\|\nabla (\cdot)\|_{r,\Omega}$ in our regime)
when $r=\infty$;
indeed, look at \cite[(5.15)]{H24} in which 
$\|\nabla (\cdot)\|_{3,\mathbb R^3}$ is still replaced by $\|\nabla(\cdot)\|_{\infty,\Omega}$ and
further splitting 
$\int_{t/2}^{t-1}+\int_{t-1}^t$ (case $s=0$) is needed, then one observes that the norm $\|\nabla (\cdot)\|_{\infty,\Omega}$
of the latter integral decays even better than $t^{-1/2}$ by use of estimate 
$\|H(U)\|_{p,\mathbb R^3}\leq\|\nabla u\|_{p,\Omega}(|\eta|+\|u\|_{\infty,\Omega})$ with some $p\in (3,\infty)$.

The last comment concerns the full slip case $\alpha\equiv 0$,
in which the system is degenerate as mentioned before.
What should be done is to deduce suitable estimates of
\begin{equation}
\begin{split}
U_{\rm i}(t)
&=u(t)\chi_\Omega+\eta(t)\chi_{\mathscr B} 
=U(t)-(\omega_0\times x)\chi_{\mathscr B}  \\
&=e^{-tA}U_0-(\omega_0\times x)\chi_{\mathscr B}
+\int_0^te^{-(t-\tau)A}H(U_{\rm i}(\tau))\,d\tau,
\end{split}
\label{int-eq-i}
\end{equation}
where we notice that $H(U)$ is determined by $U_{\rm i}$:
$H(U)=H(U_{\rm i})$.
By \eqref{decay-1} along with \eqref{decay-body} for the case (i), the difference
$e^{-tA}U_0-(\omega_0\times x)\chi_{\mathscr B}$ is fine, see also \eqref{obje-0} and \eqref{S-decay2-noch} below.
When ${\mathscr B}$ is a ball, we observe
from \eqref{anni} that 
$i(V)=(\nabla p,\eta,0)$ for all $V\in Z_q(\mathbb R^3)$, see \eqref{Z-ang} below. 
Hence, in view of the decomposition \eqref{decompo}, 
$H(U_{\rm i})|_{\mathscr B}$ is also translatory, 
so that $e^{-(t-\tau)A}H(U_{\rm i})$ enjoys desired decay
properties obtained in Theorem \ref{main}.
By use of ${\mathcal S}(t)$ defined by \eqref{obje-0} below,
the integral equation \eqref{int-eq-i} is rewritten as
\[
U_{\rm i}(t)={\mathcal S}(t)U_0+\int_0^t {\mathcal S}(t-\tau)H(U_{\rm i}(\tau))\,d\tau,
\]
see also \eqref{S-rela}.

\section{Preparatory results on the exterior problem}
\label{pre}

In this section we give several results on the exterior problem, which are useful to show Theorem \ref{main}. 
Let $\Omega$ be an exterior domain in $\mathbb R^3$ with smooth boundary $\partial\Omega\in C^{2,1}$, and let $q\in (1,\infty)$.
Given $f\in L^q(\Omega)$, let us consider the resolvent problem for the Stokes system subject to the Navier boundary condition
\begin{equation}
\begin{array}{ll}
\lambda u-\Delta u+\nabla p=f, \qquad \mbox{div $u$}=0\quad&\mbox{in $\Omega$}, \\
\nu\cdot u|_{\partial\Omega}=0, \quad [(2\mathbb Du)\nu]_\tau+\alpha u_\tau|_{\partial\Omega}=0, &
\displaystyle{\lim_{|x|\to\infty}u=0},
\end{array}
\label{St}
\end{equation}
that is rewritten as
\[
(\lambda+A_\Omega)u=P_\Omega f
\]
in the underlying space
$L^q_\sigma(\Omega)=\{u\in L^q(\Omega);\, \mbox{div $u$}=0\;\mbox{in $\Omega$}, \nu\cdot u|_{\partial\Omega}=0\}$,
in which the class $C_{0,\sigma}^\infty(\Omega)$ is dense, 
where $P_\Omega: L^q(\Omega)\to L^q_\sigma(\Omega)$ 
stands for the Fujita-Kato projection associated with the Helmholtz decomposition
of $L^q$-vector fields (\cite{FM77, Mi82, SiS92})
\[
L^q(\Omega)=L^q_\sigma(\Omega)\oplus \{\nabla p\in L^q(\Omega);\, p\in L^q_{\rm loc}(\overline{\Omega})\}
\]
and the Stokes operator $A_\Omega$ in exterior domains is defined by
\[
D_q(A_\Omega)=\{u\in L^q_\sigma(\Omega)\cap W^{2,q}(\Omega);\, [(2\mathbb Du)\nu]_\tau+\alpha u_\tau|_{\partial\Omega}=0\},
\quad A_\Omega=-P_\Omega\Delta.
\]

Shibata and Shimada \cite{SS07} deduced the resolvent estimate, which implies that the Stokes operator $-A_\Omega$
generates an analytic semigroup of class $(C_0)$ on the space $L^q_\sigma(\Omega)$ 
for every $q\in (1,\infty)$ when $\alpha\geq 0$ is a constant.
Even for the general case $\alpha\in C^1(\partial\Omega)$, 
the same result holds true
as a corollary to their results; in fact,
thanks to the generalized resolvent estimate due to 
\cite[Theorem 1.3]{SS07} for the system \eqref{St} 
subject to the full slip condition ($\alpha\equiv 0$)
in which the homogeneous one is 
replaced by the inhomogeneous condition $[(2\mathbb Du)\nu]_\tau=h$, 
we apply their estimate to \eqref{St} with $\alpha\in C^1(\partial\Omega)$ by regarding $h=-\alpha u_\tau$, 
so that the desired a priori estimate can be available for large $|\lambda|$.

Subsequently, Shimada and Yamaguchi \cite{SY08} succeeded in obtaining the $L^q$-$L^r$ decay estimates
of the Stokes semigroup.
Indeed they discussed the full slip case $\alpha\equiv 0$, but their argument still works for constant $\alpha\geq 0$. 
Recently, the first author \cite{H26} has studied the $L^q$-$L^r$ decay estimates of the Oseen semigroup,
which recovers the result of \cite{SY08} on the Stokes semigroup as a particular case,
by means of the same strategy as in Kobayashi and Shibata \cite{KS98}.

Following the notation adopted in \cite{H26},
for $q\in (1,\infty)$ and $d\geq 3$, 
we denote by $L^q_{[d]}(\Omega)$
the space of $L^q$-vector fields vanishing outside the ball $B_d$, that is, 
\begin{equation}
L^q_{[d]}(\Omega):=\{f\in L^q(\Omega);\, f(x)=0\; \mbox{a.e. $\Omega\setminus B_d$}\}.
\label{data-sp}
\end{equation}
For every $\lambda\in\mathbb C\setminus (-\infty,0)$ and $f\in L^q_{[d]}(\Omega)$, Proposition 5.2 of
\cite{H26} provides us with a solution $(u,p)$ to \eqref{St}, that is unique up to constants for the pressure, within the class
\[
(u,p)\in W^{2,q}_{\rm loc}(\overline\Omega)\times W^{1,q}_{\rm loc}(\overline\Omega), \qquad
u\in L^r(\Omega), \qquad \nabla p\in L^s(\Omega)
\]
for some $r,\, s\in (1,\infty)$.
Given $R\geq 3$, let us single out the pressure $p$ 
in such a way that
\begin{equation}
\int_{\Omega_R}p\,dx=0.
\label{single-pre}
\end{equation}
With this solution $(u,p)$ to \eqref{St} at hand,
we define the solution operator by
\begin{equation}
(S_u(\lambda), S_p(\lambda)): L^q_{[d]}(\Omega)\ni
f\mapsto (u,p), \qquad \lambda\in \mathbb C\setminus (-\infty,0).
\label{ext-sol}
\end{equation}
Since the issue here is the asymptotic behavior of the resolvent near $\lambda=0$,
it is important to introduce \eqref{ext-sol} for $\lambda=0$ as well.
If, in particular, $\lambda\in \mathbb C\setminus (-\infty,0]$, then we know that $\lambda\in \rho(-A_\Omega)$
(\cite[Proposition 5.3]{H26}) together with
\[
S_u(\lambda)f=(\lambda+A_\Omega)^{-1}P_\Omega f.
\]
Set
\begin{equation}
V_u(\tau)f=\partial_\tau S_u(i\tau)f, \qquad
V_p(\tau)f=\partial_\tau S_p(i\tau)f, \qquad \tau\in\mathbb R\setminus\{0\}
\label{ext-sol2}
\end{equation}
along the imaginary axis for $f\in L^q_{[d]}(\Omega)$.
For later use, we summarize usuful regularity properties of those solution operators 
in the following proposition, which are 
obtained in \cite{H26}. 
Since the analysis here is independent of generation of the analytic semigroup $e^{-tA_\Omega}$, 
that is relating to the asymptotic behavior of the resolvent for large $|\lambda|$, the $C^1$-regularity of the friction $\alpha$ 
plays no role (although it is needed for the well-posendness of the initial value problem as mentioned above). 
\begin{proposition}
[\cite{H26}]
Suppose that $\alpha\in C(\partial\Omega)$ is nonnegative.
Let $q\in (1,\infty)$, $d\geq 3$ and $R\geq 3$.
We set
\[
Y_j={\mathcal L}(L^q_{[d]}(\Omega),W^{j,q}(\Omega_R)), \qquad j\in \{1,2\}
\]
for simplicity of notation, where $L^q_{[d]}(\Omega)$ is given by \eqref{data-sp}.
Then the operators introduced in \eqref{ext-sol}--\eqref{ext-sol2}
enjoy the following propoerties:
For every compact set $K\subset\overline{\mathbb C_+}$ and $\theta\in (0,\frac{1}{2})$,
there is a constant $C=C(\alpha,K,\theta,q,d,R,\Omega)>0$ such that
\begin{equation}
\|S_u(\lambda)-S_u(\lambda^\prime)\|_{Y_2}
+\|S_p(\lambda)-S_p(\lambda^\prime)\|_{Y_1}
\leq C|\lambda-\lambda^\prime|^\theta
\label{ext-conti}
\end{equation}
for all $\lambda,\, \lambda^\prime\in K$. 
Moreover, there is a constant $C=C(\alpha,q,d,R,\Omega)>0$ such that 
\begin{equation}
\sup_{\lambda\in\overline{\mathbb C_+}, |\lambda|\leq 4}
\Big(\|S_u(\lambda)\|_{Y_2}+\|S_p(\lambda)\|_{Y_1}\big)
\leq C
\label{ext-0}
\end{equation}
\begin{equation}
\int_{-4}^4\Big(\|V_u(\tau)\|_{Y_2}+\|V_p(\tau)\|_{Y_1}\Big)\,d\tau\leq C
\label{ext-1}
\end{equation}
\begin{equation}
\sup_{|\tau|\leq 4}\Big(\|S_u(i(\tau+h))-S_u(i\tau)\|_{Y_2}+\|S_p(i(\tau+h)-S_p(i\tau)\|_{Y_1}\Big)\leq C|h|^{1/2}
\label{ext-2}
\end{equation}
\begin{equation}
\int_{-2}^2\Big(\|V_u(\tau+h)-V_u(\tau)\|_{Y_2}+\|V_p(\tau+h)-V_p(\tau)\|_{Y_1}\Big)\leq C|h|^{1/2}
\label{ext-3}
\end{equation}
for all 
$h\in \mathbb R$ with $|h|\leq 1$. 
\label{prop-ext}
\end{proposition}

\begin{proof}
Estimates \eqref{ext-1} and \eqref{ext-3} for $V_u(\tau)$ are found in
\cite[Proposition 5.4]{H26}.
In \cite{H26}, the regularity of the resolvent with values in $W^{2,q}(\Omega_3)$ is discussed, however,
it is obvious that all the estimates there hold true even though $\Omega_3$ is replaced by $\Omega_R$.
Regarding the others, one has only to look into details of \cite{H26}, in which the solution operator \eqref{ext-sol} 
is constructed by using solutions to the resolvent system in the whole space and in a bounded domain near the body. 
In fact, in view of the structure of the parametrix, see \cite[(5.22)]{H26},
the continuity 
\eqref{ext-conti} 
and \eqref{ext-2} as well as \eqref{ext-0} for $S_u(\lambda)$
follow from \cite[(3.10), (3.14), (4.12), (4.13), (5.26), (5.32), (5.34), the corresponding 
equality to (5.31) for $\lambda\in K$, the corresponding
estimate to (5.32) for $\lambda\in K$]{H26}.

Finally, by 
\eqref{single-pre} 
we have
\[
\|p\|_{W^{1,q}(\Omega_R)}\leq C\|\nabla p\|_{q,\Omega_R}, \qquad
\|\partial_\tau p\|_{W^{1,q}(\Omega_R)}\leq C\|\nabla\partial_\tau p\|_{q,\Omega_R}
\]
which combined with the equations
$\lambda u-\Delta u+\nabla p=f$ and 
\[
i\tau\partial_\tau u+iu-\Delta\partial_\tau u+\nabla\partial_\tau p=0
\]
imply all the estimates \eqref{ext-conti}--\eqref{ext-3} for the pressure $S_p(\lambda)$ and $V_p(\tau)$
on account of those for $S_u(\lambda)$ and $V_u(\tau)$.
The proof is complete.
\end{proof}

\section{From the exterior resolvent to the fluid-structure resolvent}
\label{from-to}

In this section, following the formulation developed by Ervedoza, Maity and Tucsnak \cite[Section 4]{EMT23},
we investigate the regularity of the fluid-structure resolvent near $\lambda=0$.
The asymptotic expansion of the Stokes resolvent in exterior domains subject to the no-slip condition due to Iwashita \cite{I89}
is used in \cite{EMT23}, whereas we intend to make full use of 
Proposition \ref{prop-ext}.

For $q\in (1,\infty)$ and $d\geq 3$, we set
\begin{equation}
L^q_{{\rm R},[d]}(\mathbb R^3):=\{F\in L^q_{\rm R}(\mathbb R^3);\, F(x)=0\;\mbox{a.e. $\mathbb R^3\setminus B_d$}\}.
\label{data-sp2}
\end{equation}
Given $F\in L^q_{{\rm R},[d]}(\mathbb R^3)$
and $(f,\kappa,\mu)=i(F)$, yielding $f\in L^q_{[d]}(\Omega)$, see \eqref{mono} and \eqref{data-sp},
let us consider
\[
(\lambda+A)U=\mathbb PF
\]
that is equivalent to the set of the equations 
\begin{equation}
\begin{split}
&\lambda u-\Delta u+\nabla p=f, \quad \mbox{div $u$}=0 \quad\mbox{in $\Omega$}, \\
&\nu\cdot (u-\eta-\omega\times x)|_{\partial\Omega}=0, \quad 
[(2\mathbb Du)\nu]_\tau+\alpha (u-\eta-\omega\times x)_\tau|_{\partial\Omega}=0, \quad
\lim_{|x|\to\infty}u=0,  \\
\end{split}
\label{FS-1}
\end{equation}
\begin{equation}
\lambda\eta+\frac{1}{m}\int_{\partial\Omega}{\mathbb S}(u,p)\nu\,d\sigma=\kappa, \qquad
\lambda \omega+J^{-1}\int_{\partial\Omega}x\times {\mathbb S}(u,p)\nu\,d\sigma=\mu. 
\label{FS-2}
\end{equation}

As observed in subsection \ref{axisym}, 
the system \eqref{FS-1}--\eqref{FS-2} for the particular cases (i) and (ii) in Theorem \ref{main}
is, respectively, rewritten as follows.
Consider the case (ii), then by \eqref{sym-2} together with
$\gamma\cdot (J^{-1}a)=c^{-1}(\gamma\cdot a)$ for all $a\in\mathbb R^3$, see \eqref{sym-inertia},
we are led to
\begin{equation}
\begin{split}
&\lambda u-\Delta u+\nabla p=f, \quad \mbox{div $u$}=0 \quad\mbox{in $\Omega$}, \\
&\nu\cdot \big[u-\eta-\{(\gamma\times \omega)\times\gamma\}\times x\big]\Big|_{\partial\Omega}=0, \quad 
[(2\mathbb Du)\nu]_\tau\big|_{\partial\Omega}=0, \quad
\lim_{|x|\to\infty}u=0,  \\
\end{split}
\label{FS-ii-1}
\end{equation}
\begin{equation}
\begin{split}
&\lambda\eta+\frac{1}{m}\int_{\partial\Omega}{\mathbb S}(u,p)\nu\,d\sigma=\kappa, \\
&\lambda (\gamma\times \omega)\times \gamma+
J^{-1}\int_{\partial\Omega}x\times \mathbb S(u,p)\nu\,d\sigma
=(\gamma\times \mu)\times \gamma, 
\end{split}
\label{FS-ii-2}
\end{equation}
as well as
\begin{equation}
\gamma\cdot \omega=\frac{\gamma\cdot\mu}{\lambda}
\label{FS-ii-3}
\end{equation}
unless $\lambda=0$.
When $\lambda=0$, we have a solution only if $\gamma\cdot\mu=0$.
Note that
\begin{equation}
J^{-1}\int_{\partial\Omega}x\times \mathbb S(u,p)\nu\,d\sigma
=\left(\gamma\times J^{-1}\int_{\partial\Omega}x\times {\mathbb S}(u,p)\nu\,d\sigma\right)\times\gamma
\label{torque-ii}
\end{equation}
in \eqref{FS-ii-2}.

As to the case (i),
it follows from \eqref{ball} that
\begin{equation}
\begin{split}
&\lambda u-\Delta u+\nabla p=f, \quad \mbox{div $u$}=0 \quad\mbox{in $\Omega$}, \\
&\nu\cdot (u-\eta)|_{\partial\Omega}=0, \quad 
[(2\mathbb Du)\nu]_\tau|_{\partial\Omega}=0, \quad
\lim_{|x|\to\infty}u=0,  \\
\end{split}
\label{FS-i-1}
\end{equation}
\begin{equation}
\lambda\eta+\frac{1}{m}\int_{\partial\Omega}{\mathbb S}(u,p)\nu\,d\sigma=\kappa,
\label{FS-i-2}
\end{equation}
as well as
\begin{equation}
\omega=\frac{\mu}{\lambda}
\label{FS-i-3}
\end{equation}
unless $\lambda=0$.
When $\lambda=0$, there is a solution only if $\mu=0$.

Let us take a standard lifting function 
\begin{equation}
\ell(\eta,\omega)=\frac{1}{2}\mbox{rot $\Big(\phi(x)(\eta\times x-|x|^2\omega)\Big)$}
\label{lift}
\end{equation}
of the rigid motion $\eta+\omega\times x$ at $\partial\Omega$ in \eqref{FS-1},
where a cut-off function $\phi$ is fixed in such a way that
\begin{equation}
\phi\in C_0^\infty(B_3), \qquad 0\leq \phi\leq 1,\qquad \phi =1 \;\mbox{on $B_2$}.
\label{cut}
\end{equation}
Then we have
\begin{equation}
\begin{split}
&\mbox{div $\ell(\eta,\omega)$}=0, \qquad
\ell(\eta,\omega)|_{B_2}=\eta+\omega\times x, \\
&\ell(\eta,\omega)\in C_0^\infty(B_3),\qquad 
\|\ell(\eta,\omega)\|_{W^{2,q}(\Omega)}\leq C_q(|\eta|+|\omega|)
\end{split}
\label{lift-est}
\end{equation}
for every $q\in (1,\infty)$.
It is convenient to introduce the 
lifting operator
\[
\ell: (\eta,\omega)\mapsto \ell(\eta,\omega).
\]
Then, for the case (ii) (resp. case (i)), one may take
the lifting function
\begin{equation}
\ell\big(\eta,(\gamma\times\omega)\times\gamma\big)\qquad
\Big(\mbox{resp. $\ell(\eta,0)$}\Big)
\label{lift2}
\end{equation}
which fulfills \eqref{lift-est} with obvious changes, respectively.

Let $\lambda\in \mathbb C\setminus (-\infty,0)$.
By use of the pair of the solution operators \eqref{ext-sol} to \eqref{St} defined on $L^q_{[d]}(\Omega)$, a unique solution to \eqref{FS-1}
is described as
\begin{equation}
u=u_1+u_2, \qquad p=p_1+p_2
\label{decouple}
\end{equation}
with
\begin{equation}
\begin{split}
&u_1=\ell(\eta,\omega)+S_u(\lambda)\Big(\Delta\ell(\eta,\omega)-\lambda\ell(\eta,\omega)\Big), \qquad
p_1=S_p(\lambda)\Big(\Delta \ell(\eta,\omega)-\lambda\ell(\eta,\omega)\Big), \\
&u_2=S_u(\lambda)f, \qquad p_2=S_p(\lambda)f, \\
\end{split}
\label{decouple2}
\end{equation}
which solve
\begin{equation}
\begin{split}
&\lambda u_1-\Delta u_1+\nabla p_1=0, \quad \mbox{div $u_1$}=0 \quad\mbox{in $\Omega$}, \\
&\nu\cdot (u_1-\eta-\omega\times x)|_{\partial\Omega}=0, \quad 
[(2\mathbb Du_1)\nu]_\tau+\alpha (u_1-\eta-\omega\times x)_\tau|_{\partial\Omega}=0, \quad
\lim_{|x|\to\infty}u_1=0,  \\
\end{split}
\label{FS-auxi-1}
\end{equation}
and
\begin{equation}
\begin{split}
&\lambda u_2-\Delta u_2+\nabla p_2=f, \quad \mbox{div $u_2$}=0 \quad\mbox{in $\Omega$}, \\
&\nu\cdot u_2|_{\partial\Omega}=0, \quad 
[(2\mathbb Du_2)\nu]_\tau+\alpha u_{2,\tau}|_{\partial\Omega}=0, \quad
\lim_{|x|\to\infty}u_2=0. \\
\end{split}
\label{FS-auxi-2}
\end{equation}

Let us define the matrix $T_\lambda\in \mathbb C^{6\times 6}$ by
\begin{equation}
T_\lambda\left(
\begin{array}{c}
\eta \\
\omega
\end{array}
\right)=
\left(
\begin{array}{c}
m^{-1}\int_{\partial\Omega}\mathbb S(u_1,p_1)\nu\,d\sigma \\
J^{-1}\int_{\partial\Omega}x\times \mathbb S(u_1,p_1)\nu\,d\sigma
\end{array}
\right),
\label{auxi-mat}
\end{equation}
then 
our problem 
is reduced to finding $(\eta,\omega)$ to the equation
\begin{equation}
(\lambda+T_\lambda)\left(
\begin{array}{c}
\eta \\
\omega
\end{array}
\right)=
\left(
\begin{array}{c}
\kappa-m^{-1}\int_{\partial\Omega}\mathbb S(u_2,p_2)\nu\,d\sigma \\
\mu-J^{-1}\int_{\partial\Omega}x\times\mathbb S(u_2,p_2)\nu\,d\sigma
\end{array}
\right)
\label{solid}
\end{equation}
through \eqref{FS-2} with \eqref{decouple}--\eqref{FS-auxi-2}.
Once we verify the invertibility of $\lambda+T_\lambda$ (see Corollary \ref{cor-reg} below)
to obtain $(\eta,\omega)$,
we are led to a representation of the resolvent
\begin{equation}
(\lambda+A)^{-1}\mathbb PF=u \chi_\Omega+(\eta+\omega\times x)\chi_{\mathscr B}
\label{resol-rep}
\end{equation}
in terms of $(S_u(\lambda),S_p(\lambda))$ via $u$ given by \eqref{decouple}--\eqref{decouple2}.
We will see that this argument works well for all the cases (i)--(iii) in Theorem \ref{main} as long as 
$\lambda\in\mathbb C\setminus (-\infty,0]$,
however, in order to show the asymptotic behavior of \eqref{resol-rep} near $\lambda =0$ through
the one of $(\lambda+T_\lambda)^{-1}$,
the aforementioned procedure should cover $\lambda=0$ as well.
We will show in Corollary \ref{cor-reg} below that this is indeed possible only for the case (iii) in Theorem \ref{main}.
The other cases (i) and (ii) are treated 
by use of \eqref{FS-i-1}--\eqref{FS-i-2} and \eqref{FS-ii-1}--\eqref{FS-ii-2}, respectively, in the following way.

Still with the matrix $T_\lambda$ given by \eqref{auxi-mat} but with 
the other lifting function \eqref{lift2}
for $(u_1,p_1)$ in \eqref{decouple2},
the equation \eqref{solid} is replaced by 
\begin{equation}
(\lambda+T_\lambda)\left(
\begin{array}{c}
\eta \\
(\gamma\times \omega)\times\gamma
\end{array}
\right)=
\left(
\begin{array}{c}
\kappa-m^{-1}\int_{\partial\Omega}\mathbb S(u_2,p_2)\nu\,d\sigma \\
(\gamma\times \mu)\times\gamma-J^{-1}\int_{\partial\Omega}x\times\mathbb S(u_2,p_2)\nu\,d\sigma
\end{array}
\right)
\label{solid-ii}
\end{equation}
for the case (ii), which arises from \eqref{FS-ii-2}, and by
\begin{equation}
(\lambda+T_\lambda)\left(
\begin{array}{c}
\eta \\
0
\end{array}
\right)=
\left(
\begin{array}{c}
\kappa-m^{-1}\int_{\partial\Omega}\mathbb S(u_2,p_2)\nu\,d\sigma \\
0   
\end{array}
\right)
\label{solid-i}
\end{equation}
for the case (i), where the latter three rows of 
both sides of
\eqref{solid-i} are identically zero 
since
\[
\int_{\partial\Omega}x\times \mathbb S(u_k,p_k)\nu\,d\sigma=0 \qquad (k=1,2)
\]
follows from $[(\mathbb Du_k)\nu]_\tau|_{\partial\Omega}=0$,
see \eqref{ball} through \eqref{sym-0}.
Hence, \eqref{solid-i} is essentially the three-dimensional equation and corresponds exactly to \eqref{FS-i-2}.

By taking into account \eqref{FS-ii-3} and \eqref{FS-i-3}, 
we have the other representations of the resolvent
\begin{equation}
\left\{
\begin{array}{rl}
\mbox{case (i)}\quad&
(\lambda+A)^{-1}\mathbb PF={\mathcal Q}_{\rm i}(\lambda)F_{\rm i}+\left(\frac{\mu}{\lambda}\times x\right)\chi_{\mathscr B}, \\ 
\quad&{\mathcal Q}_{\rm i}(\lambda)F_{\rm i}=u\chi_\Omega+\eta\chi_{\mathscr B}  \\
\quad&\quad \\
\mbox{case (ii)}\quad&
(\lambda+A)^{-1}\mathbb PF={\mathcal Q}_{\rm ii}(\lambda)F_{\rm ii}+\left(\frac{\gamma\cdot\mu}{\lambda}\gamma\times x\right)
\chi_{\mathscr B},    \\
\quad& {\mathcal Q}_{\rm ii}(\lambda)F_{\rm ii}=u\chi_\Omega+
\big[\eta+\{(\gamma\times \omega)\times\gamma\}\times x\big]\chi_{\mathscr B}
\end{array}
\right.
\label{resol-rep2}
\end{equation}
in terms of $(S_u(\lambda),S_p(\lambda))$ via \eqref{decouple}--\eqref{decouple2}, in which the lifting function
$\ell(\eta,\omega)$ should be respectively replaced by \eqref{lift2}, as well as \eqref{solid-ii}--\eqref{solid-i}.
Here, $F_{\rm i}$ and $F_{\rm ii}$ are given by \eqref{IC-modi}.
Our strategy for those cases is to derive the asymptotic behavior of ${\mathcal Q}_{\rm i}(\lambda)$ and 
${\mathcal Q}_{\rm ii}(\lambda)$ in \eqref{resol-rep2}
near $\lambda=0$.
Note that
\begin{equation}
\begin{split}
&\mbox{case (i)}\qquad
{\mathcal Q}_{\rm i}(\lambda)F_{\rm i}=(\lambda+A)^{-1}\mathbb PF_{\rm i}, \\
&\mbox{case (ii)}\qquad
{\mathcal Q}_{\rm ii}(\lambda)F_{\rm ii}=(\lambda+A)^{-1}\mathbb PF_{\rm ii}.
\end{split}
\label{resol-rep3}
\end{equation}
In fact, for the case (ii), we see from \eqref{sym-1} that $(\gamma\times x)\chi_{\mathscr B}$ fulfills \eqref{normal-va},
thereby,
$(\gamma\times x)\chi_{\mathscr B}\in X_q(\mathbb R^3)$
and 
that
\[ 
(\lambda+A)^{-1}\mathbb P\big[\{(\gamma\cdot\mu)\gamma\times x\}\chi_{\mathscr B}\big]
=\left(\frac{\gamma\cdot\mu}{\lambda}\gamma\times x\right)\chi_{\mathscr B},
\]
by taking into account \eqref{FS-ii-3} and Proposition \ref{uni} below on the uniqueness.
This observation implies \eqref{resol-rep3}.
The case (i) is verified similarly.

To begin with, we give several useful estimates of $T_\lambda$ near $\lambda=0$.
\begin{lemma}
Suppose that $\alpha\in C(\partial\Omega)$ is nonnegative.
For every compact set $K\subset \overline{\mathbb C_+}$ and $\theta\in (0,\frac{1}{2})$, there is a constant
$C=C(\alpha,K,\theta,\Omega)$ such that
the matrix $T_\lambda\in \mathbb C^{6\times 6}$ defined by \eqref{auxi-mat} satisfies
\begin{equation}
|T_\lambda-T_{\lambda^\prime}|_{\mathbb C^{6\times 6}}\leq C|\lambda-\lambda^\prime|^\theta
\label{T-conti}
\end{equation}
for all $\lambda,\, \lambda^\prime\in K$.

We set
\[
L(\tau)=\frac{d}{d\tau}T_{i\tau}, \qquad \tau\in\mathbb R\setminus\{0\}.
\]
Then there is a constant $C=C(\alpha,\Omega)>0$ such that $T_\lambda$ and $L(\tau)$ enjoy
\begin{equation}
\sup_{\lambda\in\overline{\mathbb C_+}, |\lambda|\leq 4}|T_\lambda|_{\mathbb C^{6\times 6}}\leq C
\label{mat-est}
\end{equation}
\begin{equation}
\int_{-4}^4 \left|L(\tau)\right|_{\mathbb C^{6\times 6}}\,d\tau\leq C
\label{auxi-5}
\end{equation}
\begin{equation}
\sup_{|\tau|\leq 4}|T_{i(\tau+h)}-T_{i\tau}|_{\mathbb C^{6\times 6}}\leq C|h|^{1/2}
\label{auxi-conti}
\end{equation}
\begin{equation}
\int_{-2}^2 |L(\tau+h)-L(\tau)|_{\mathbb C^{6\times 6}}\,d\tau\leq C|h|^{1/2}
\label{auxi-6}
\end{equation}
for all $h\in\mathbb R$ with $|h|\leq 1$.
\label{lem-mat}
\end{lemma}

\begin{proof}
We fix $q\in (1,\infty)$ and $R\geq 3$ arbitrarily, then we have
$\ell(\eta,\omega)\in L^q_{[3]}(\Omega)$, see
\eqref{lift} and \eqref{lift-est}.
In view of \eqref{decouple2}--\eqref{FS-auxi-1} and \eqref{auxi-mat},
we use the trace estimate to deduce
\begin{equation*}
\begin{split}
&\quad\left|
\big(T_{\lambda}-T_{\lambda^\prime}\big)
\left(
\begin{array}{c}
\eta \\
\omega
\end{array}
\right)
\right|_{\mathbb C^6}  \\
&\leq C\Big(\|(S_u(\lambda)-S_u(\lambda^\prime))\Delta\ell(\eta,\omega)\|_{W^{2,q}(\Omega_R)}
+|\lambda-\lambda^\prime|\|S_u(\lambda)\ell(\eta,\omega)\|_{W^{2,q}(\Omega_R)}  \\
&\quad +|\lambda^\prime|\|(S_u(\lambda)-S_u(\lambda^\prime))\ell(\eta,\omega)\|_{W^{2,q}(\Omega_R)} 
+\|(S_p(\lambda)-S_p(\lambda^\prime))\Delta\ell(\eta,\omega)\|_{W^{1,q}(\Omega_R)}  \\
&\quad +|\lambda-\lambda^\prime|\|S_p(\lambda)\ell(\eta,\omega)\|_{W^{1,q}(\Omega_R)}  
+|\lambda^\prime|\|(S_p(\lambda)-S_p(\lambda^\prime))\ell(\eta,\omega)\|_{W^{1,q}(\Omega_R)}\Big),
\end{split}
\end{equation*}
which combined with \eqref{ext-conti} 
as well as \eqref{lift-est} leads us to \eqref{T-conti} and, therefore, 
\eqref{mat-est} as well.
By the similar computation
along the imaginary axis, it is readily seen that \eqref{ext-2} implies \eqref{auxi-conti} for $|h|\leq 1$.

Likewise, we find from \eqref{auxi-mat} that
\begin{equation*}
\begin{split}
\left| L(\tau)
\left(\begin{array}{c}\eta \\ \omega\end{array}\right)\right|_{\mathbb C^6}
&\leq C\Big(\|V_u(\tau)\big(\Delta \ell(\eta,\omega)-i\tau\ell(\eta,\omega)\big)\|_{W^{2,q}(\Omega_R)}
+\|S_u(i\tau)\ell(\eta,\omega)\|_{W^{2,q}(\Omega_R)}  \\
&\quad +\|V_p(\tau)\big(\Delta(\eta,\omega)-i\tau\ell(\eta,\omega)\big)\|_{W^{1,q}(\Omega_R)}
+\|S_p(i\tau)\ell(\eta,\omega)\|_{W^{1,q}(\Omega_R)}\Big),
\end{split}
\end{equation*}
where $V_u(\tau)$ and $V_p(\tau)$ are as in \eqref{ext-sol2}. 
We thus conclude \eqref{auxi-5} from \eqref{ext-0}--\eqref{ext-1} and \eqref{lift-est}. 

Since
\begin{equation*}
\begin{split}
&\quad\left|\big(L(\tau+h)-L(\tau)\big)
\left(
\begin{array}{c}
\eta \\
\omega
\end{array}
\right)\right|_{\mathbb C^6}  \\
&\leq C
\Big(
\|(V_u(i(\tau+h))-V_u(i\tau))\big(\Delta\ell(\eta,\omega)-i(\tau+h)\ell(\eta,\omega)\big)\|_{W^{2,q}(\Omega_R)}  \\
&\quad +|h|\|V_u(i\tau)\ell(\eta,\omega)\|_{W^{2,q}(\Omega_R)}
+\|(S_u(i(\tau+h))-S_u(\tau))\ell(\eta,\omega)\|_{W^{2,q}(\Omega_R)}  \\
&\quad +\|(V_p(i(\tau+h))-V_p(i\tau))\big(\Delta\ell(\eta,\omega)-i(\tau+h)\ell(\eta,\omega)\big)\|_{W^{1,q}(\Omega_R)}  \\
&\quad +|h|\|V_p(i\tau)\ell(\eta,\omega)\|_{W^{1,q}(\Omega_R)}
+\|(S_p(i(\tau+h))-S_p(\tau))\ell(\eta,\omega)\|_{W^{1,q}(\Omega_R)}\Big),
\end{split}
\end{equation*}
we utilize \eqref{ext-1}--\eqref{ext-3} to find \eqref{auxi-6}.
The proof is complete.
\end{proof}

Let us show Proposition \ref{uni} below on
the uniqueness of solutions to the fluid-structure resolvent systems \eqref{FS-1}--\eqref{FS-2},
\eqref{FS-ii-1}--\eqref{FS-ii-2} and \eqref{FS-i-1}--\eqref{FS-i-2}
with the aid of the following lemma,
which is found in \cite{AR14,FG22,KLT25,W03}.
\begin{lemma}
[{\cite[Theorem 1]{FG22}}, {\cite[Corollary 7.3]{KLT25}}]
Let ${\mathscr B}$ be a bounded
domain as in \eqref{dom}.
Define the subspace $M$ of $\mathbb C^6$ by
\begin{equation}
M=\{(\eta,\omega)\in \mathbb C^3\times \mathbb C^3;\; \nu\cdot(\eta+\omega\times x)|_{\partial{\mathscr B}}=0\}.
\label{rigid-tan}
\end{equation}
Then, $\dim M=0$ if and only if ${\mathscr B}$ is not axially symmetric. 

Let $\gamma\in\mathbb R^3$ satisfy $|\gamma|=1$.
If ${\mathscr B}$ is axisymmetric about a single axis $\gamma\mathbb R$ (which should pass through the origin since
$\int_{\mathscr B}x\,dx=0$), 
then $\dim M=1$, to be precise, 
$M=\{(0,c\gamma)\in\mathbb C^3\times \mathbb C^3;\; c\in\mathbb C\}$.

If ${\mathscr B}$ is axisymmetric about two distinct axes, then it must be a ball 
and $\dim M=3$, to be precise, 
$M=\{(0,\omega)\in\mathbb C^3\times\mathbb C^3\}$. 
\label{lem-rigid-tan}
\end{lemma}
\begin{proposition}
Suppose that $\alpha\in C(\partial\Omega)$ is nonnegative and that $\lambda\in\mathbb C\setminus (-\infty,0)$.

\begin{enumerate}
\item
Let $q,\,r,\,s\in (1,\infty)$. 
Then the only solution
$(u,p,\eta,\omega)\in W^{2,q}_{\rm loc}(\overline\Omega)\times W^{1,q}_{\rm loc}(\overline\Omega)\times\mathbb C^3\times \mathbb C^3$
to \eqref{FS-1}--\eqref{FS-2} with $(f,\kappa,\mu)=(0,0,0)$ 
satisfying
\begin{equation}
u\in L^r(\Omega), \qquad 
\nabla p\in L^s(\Omega)
\label{weak-sum}
\end{equation}
is $(u,p,\eta,\omega)=(0,p_\infty,0,0)$
with some constant $p_\infty\in\mathbb C$ except 
when $\lambda=0$ for the cases \emph{(i)} and \emph{(ii)} in Theorem \ref{main}.

\item
Consider the case \emph{(ii)}.
Then the only solution $(u,p,\eta,(\gamma\times\omega)\times\gamma)$
of the same class as in the previous item
to \eqref{FS-ii-1}--\eqref{FS-ii-2}
with $(f,\kappa,(\gamma\times\mu)\times\gamma)=(0,0,0)$ is
$(u,p,\eta,(\gamma\times \omega)\times\gamma)=(0,p_\infty,0,0)$ with some constant $p_\infty\in\mathbb C$.

\item 
Consider the case \emph{(i)}.
Then the only solution $(u,p,\eta)$ 
of the same class as in the item 1
to \eqref{FS-i-1}--\eqref{FS-i-2}
with $(f,\kappa)=(0,0)$
is $(u,p,\eta)=(0,p_\infty,0)$ with some constant $p_\infty\in\mathbb C$.
\end{enumerate}
\label{uni}
\end{proposition}

\begin{proof}
First of all, we observe
$(u,p)\in W^{2,2}_{\rm loc}(\overline{\Omega})\times W^{1,2}_{\rm loc}(\overline{\Omega})$ (even though $q$ is close to $1$)
by bootstrap argument after subtracting the lifting function \eqref{lift} with the aid of the regularity theory for the Stokes system
subject to the Navier boundary condition (\cite{A21,AR14,SS07}) as well as a cut-off procedure.
On the other hand, weak summability \eqref{weak-sum} at infinity
implies that
$(u,p-p_\infty)$ with some constant $p_\infty\in\mathbb C$ behaves like the fundamental solution
to the Stokes resolvent system (or to the Stokes system when $\lambda =0$) in the whole space $\mathbb R^3$.
This is verified by reducing the exterior problem to the one in the whole space by a cut-off technique.

Consider the problem \eqref{FS-1}--\eqref{FS-2}
with $(f,\kappa,\mu)=(0,0,0)$.
By taking into account the behavior of $(u,p-p_\infty)$ at infinity mentioned above,
we multiply the equation
\[
\lambda u-\mbox{div $\mathbb S(u,p-p_\infty)$}=0
\]
by $\phi_\rho\overline{u}$ and let $\rho\to\infty$ to find that the third term of \eqref{go-to-ene} below vanishes
and, therefore, that the energy relation is available,
where $\phi_\rho(x)=\phi(x/\rho)$ with $\phi$ being given by \eqref{cut}.
In fact, by integration by parts and by using the Navier boundary condition, where we note
$[(2\mathbb Du)\nu]_\tau=[\mathbb S(u,p-p_\infty)\nu]_\tau$
at $\partial\Omega$, we obtain
\begin{equation}
\begin{split}
&\lambda\int_\Omega\phi_\rho |u|^2\,dx
+2\int_\Omega\phi_\rho |\mathbb Du|^2\,dx
+\int_{B_{3\rho}\setminus \overline{B_{2\rho}}}\mathbb S(u,p-p_\infty):(\overline{u}\otimes \nabla\phi_\rho)\,dx  \\
&\qquad +\int_{\partial\Omega}\alpha(x) |u-\eta-\omega\times x|^2\,d\sigma
=\int_{\partial\Omega}(\overline{\eta}+\overline{\omega}\times x)\cdot [\mathbb S(u,p-p_\infty)\nu]\,d\sigma.
\end{split}
\label{go-to-ene}
\end{equation}
Since $\int_{\partial\Omega}\nu\,d\sigma=\int_{\partial\Omega}x\times \nu\,d\sigma=0$, the system \eqref{FS-2}
with $(\kappa,\mu)=(0,0)$ may be rewritten as
\[
\lambda m\eta+\int_{\partial\Omega}\mathbb S(u,p-p_\infty)\nu\,d\sigma=0, \qquad
\lambda J\omega+\int_{\partial\Omega}x\times \mathbb S(u,p-p_\infty)\nu\,d\sigma=0.
\]
Multiplying them by $\overline{\eta}$ and $\overline{\omega}$, respectively,
adding the resulting equalities to \eqref{go-to-ene}, and then letting $\rho\to\infty$ to furnish
%
\begin{equation}
\lambda\Big[\|u\|_{2,\Omega}^2+m|\eta|^2+(J\omega)\cdot\overline{\omega}\Big]
+2\|\mathbb Du\|_{2,\Omega}^2+\int_{\partial\Omega}\alpha(x)|u-\eta-\omega\times x|^2\,d\sigma=0.
\label{energy}
\end{equation}
Note that $(J\omega)\cdot\overline{\omega}>0$ 
since $J\in\mathbb R^{3\times 3}$ is symmetric and positive. 
Then it is obvious that $(u,\eta,\omega)=(0,0,0)$, with $p$ being a constant,
if $\lambda\in\mathbb C\setminus (-\infty,0]$.

Let $\lambda=0$ in \eqref{FS-1}--\eqref{FS-2}.
From \eqref{energy} with $\alpha\geq 0$ we find $\mathbb Du=O$, 
which implies $u=0$ on account of the condition at infinity.
Unless $\alpha\in C(\partial\Omega)$ is identically zero,
there is a nonempty open subset $\Gamma$ of $\partial\Omega$ 
(with respect to the induced topology)
such that $\eta+\omega\times x=0$ for every $x\in\Gamma$.
If $\omega$ were nonzero vector, then $\omega$ and $x-x^\prime$ would be linearly dependent 
for every $x,\, x^\prime\in\Gamma$ with $x\neq x^\prime$,
yielding $\Gamma$ being a line segment, a contradiction.
We are thus led to $(\eta,\omega)=(0,0)$.

For the full slip case $\alpha\equiv 0$ as well as $\lambda=0$ still in \eqref{FS-1}--\eqref{FS-2},
the only condition which the rigid motion fulfills is that
$\nu\cdot(\eta+\omega\times x)|_{\partial{\mathscr B}}=0$.
Then, by Lemma \ref{lem-rigid-tan}, we find $(\eta,\omega)=(0,0)$ 
provided the body ${\mathscr B}$ is not axially symmetric, that is covered by
the case (iii) in Theorem \ref{main}.

Let us turn to the cases (i) and (ii). 
Then the same computation as in deduction of \eqref{energy} implies that \eqref{FS-i-1}--\eqref{FS-i-2} lead to
\[
\lambda \Big[\|u\|_{2,\Omega}^2+m|\eta|^2\Big]+2\|\mathbb Du\|_{2,\Omega}^2=0
\]
for the case (i) and that \eqref{FS-ii-1}--\eqref{FS-ii-2} yield
\[
\lambda \Big[\|u\|_{2,\Omega}^2+m|\eta|^2
+\big(J\{(\gamma\times \omega)\times \gamma\}\big)\cdot\{(\gamma\times\overline{\omega})\times\gamma\}\Big]
+2\|\mathbb Du\|_{2,\Omega}^2=0
\]
for the case (ii).
It thus suffices to consider the case $\lambda=0$, for which we have $\mathbb Du=0$ and, thereby, $u=0$, yielding
\begin{equation*}
\begin{array}{ll}
\nu\cdot\eta|_{\partial{\mathscr B}}=0 
&\mbox{case (i)}, \\
\nu\cdot\big[\eta+\{(\gamma\times \omega)\times\gamma\}\times x\big]\Big|_{\partial{\mathscr B}}=0 \qquad
&\mbox{case (ii)}.
\end{array}
\end{equation*}
From the latter half of Lemma \ref{lem-rigid-tan} it follows that
nontrivial solutions are not allowed in those cases.
The proof is complete.
\end{proof}
\begin{corollary}
Suppose that $\alpha\in C(\partial\Omega)$ is nonnegative.
We assume that
$\lambda\in \mathbb C\setminus (-\infty,0)$.
Let $T_\lambda\in\mathbb C^{6\times 6}$ be the matrix defined by \eqref{auxi-mat}.

\begin{enumerate}
\item
The matrix $\lambda+T_\lambda$ is invertible except when $\lambda=0$ 
for the cases \emph{(i)} and \emph{(ii)} in Theorem \ref{main}.

\item
Let $\lambda=0$ for those cases above. 
Set
\begin{equation}
\begin{array}{rll}
E&=\{(\eta,0)\in \mathbb C^3\times \mathbb C^3\} 
&\mbox{\rm case (i)},  \\
E_\gamma&
=\{(\eta,\omega)\in\mathbb C^3\times \mathbb C^3;\; \gamma\cdot\omega=0\}\qquad
&\mbox{\rm case (ii)}.
\end{array}
\label{auxi-sp}
\end{equation}
Then we have
\begin{equation}
\begin{array}{ll}
(\ker T_0)\cap E=\{(0,0)\}&\mbox{\rm case (i)}, \\
(\ker T_0)\cap E_\gamma=\{(0,0)\}\qquad &\mbox{\rm case (ii)}.
\end{array}
\label{cond-uni}
\end{equation}

\item
For the cases {\rm (i)} and {\rm (ii)} in Theorem \ref{main}, 
the range of $\lambda+T_\lambda$, that is restricted to $E$ and $E_\gamma$ respectively, is given by
\begin{equation}
\begin{array}{ll}
{\rm R}\big((\lambda+T_\lambda)|_{E}\big)=E  &\mbox{\rm case (i)}, \\
{\rm R}\big((\lambda+T_\lambda)|_{E_\gamma}\big)=E_\gamma \qquad&\mbox{\rm case (ii)},
\end{array}
\label{range}
\end{equation}
for all $\lambda\in\mathbb C\setminus (-\infty,0)$ including $\lambda=0$.
\end{enumerate}
\label{cor-reg}
\end{corollary}

\begin{proof}
Suppose
\[
(\lambda+T_\lambda)
\left(
\begin{array}{c}
\eta \\
\omega
\end{array}
\right)=
\left(
\begin{array}{c}
0 \\
0
\end{array}
\right).
\]
In view of the definition \eqref{auxi-mat} of $T_\lambda$, we see that
$(u_1,p_1,\eta,\omega)$ obeys
\eqref{FS-auxi-1} together with \eqref{FS-2} with $(\kappa,\mu)=(0,0)$.
Since $(u_1,p_1)$ given by \eqref{decouple2} is of class specified in the item {\it 1} of Proposition \ref{uni},
we know from this proposition that $(\eta,\omega)=(0,0)$ except when $\lambda=0$ for the cases (i) and (ii).

Let $\alpha\equiv 0$ and consider those cases, then we see that $\ker T_0=M$
being given by \eqref{rigid-tan}.
By Lemma \ref{lem-rigid-tan} we obtain
\begin{equation}
\ker T_0=\left\{
\begin{array}{ll}
\{(0,\omega)\in \mathbb C^3\times\mathbb C^3\} 
&\mbox{case {\rm (i)}},  \\
\{(0,c\gamma)\in\mathbb C^3\times\mathbb C^3;\; c\in\mathbb C\}\qquad &\mbox{case {\rm (ii)}}, \\
\end{array}
\right.
\label{ker-T0}
\end{equation}
which immediately concludes \eqref{cond-uni}.

Finally, by the definition \eqref{auxi-mat} of $T_\lambda$, it is given by
\[
T_\lambda\left(
\begin{array}{c}
\eta \\
(\gamma\times \omega)\times\gamma
\end{array}
\right)=\left(
\begin{array}{c}
m^{-1}\int_{\partial\Omega}\mathbb S(u_1,p_1)\nu\,d\sigma \\
J^{-1}\int_{\partial\Omega}x\times \mathbb S(u_1,p_1)\nu\,d\sigma
\end{array}
\right)
\]
for the case (ii), where $(u_1,p_1)$ is as in \eqref{decouple2}, however, with the lifting function $\ell(\eta,\omega)$ being 
replaced by \eqref{lift2}.
Taking into account \eqref{torque-ii},
we observe that
$(\lambda+T_\lambda)E_\gamma\subset E_\gamma$.
We know
$\dim {\rm R}(T_0)=5$ from \eqref{ker-T0}, while $\dim {\rm R}\big((\lambda+T_\lambda)|_{E_\gamma}\big)=5$
for $\lambda\in\mathbb C\setminus (-\infty,0]$ is obvious since $\lambda+T_\lambda$ is invertible.
We are thus led to \eqref{range} for the case (ii).
The case (i) is discussed similarly.
The proof is complete.
\end{proof}

Let $\lambda\in \mathbb C\setminus (-\infty,0)$ and recall the equation \eqref{solid}.
We find 
\begin{equation}
\left(
\begin{array}{c}
\eta \\
\omega
\end{array}
\right)=(\lambda+T_\lambda)^{-1}\left(
\begin{array}{c}
\kappa-m^{-1}\int_{\partial\Omega}\mathbb S(u_2,p_2)\nu\,d\sigma  \\
\mu-J^{-1}\int_{\partial\Omega}x\times\mathbb S(u_2,p_2)\nu\,d\sigma
\end{array}
\right)
\label{rigid-inv}
\end{equation}
from Corollary \ref{cor-reg}
for all the cases in Theorem \ref{main}; moreover, $\lambda=0$ is also included here except for the cases (i) and (ii).
For those cases, 
since the right-hand sides of \eqref{solid-ii}--\eqref{solid-i} belong to the range \eqref{range}
on account of \eqref{torque-ii} for the case (ii),
the conditional uniqueness \eqref{cond-uni} 
allows us to include $\lambda=0$, respectively, in
\begin{equation}
\left(
\begin{array}{c}
\eta \\
0
\end{array}
\right)=\big[(\lambda+T_\lambda)|_{E}\big]^{-1}\left(
\begin{array}{c}
\kappa-m^{-1}\int_{\partial\Omega}\mathbb S(u_2,p_2)\nu\,d\sigma  \\
0 
\end{array}
\right)
\label{rigid-inv-i}
\end{equation}
for the case (i), and in
\begin{equation}
\left(
\begin{array}{c}
\eta \\
(\gamma\times\omega)\times \gamma
\end{array}
\right)=\big[(\lambda+T_\lambda)|_{E_\gamma}\big]^{-1}\left(
\begin{array}{c}
\kappa-m^{-1}\int_{\partial\Omega}\mathbb S(u_2,p_2)\nu\,d\sigma  \\
(\gamma\times\mu)\times\gamma-J^{-1}\int_{\partial\Omega}x\times\mathbb S(u_2,p_2)\nu\,d\sigma
\end{array}
\right)
\label{rigid-inv-ii}
\end{equation}
for the case (ii),
where $E$ and $E_\gamma$ are given by \eqref{auxi-sp}.
We see from \eqref{range} that
$\big[(\lambda+T_\lambda)|_{E}\big]^{-1}=(\lambda+T_\lambda)^{-1}|_E$
for the case (i) and
$\big[(\lambda+T_\lambda)|_{E_\gamma}\big]^{-1}=(\lambda+T_\lambda)^{-1}|_{E_\gamma}$
for the case (ii) unless $\lambda=0$.

With \eqref{rigid-inv}--\eqref{rigid-inv-ii} at hand,
we introduce the operator
${\mathcal Q}_{\rm b}(\lambda)$ by 
\begin{equation}
{\mathcal Q}_{\rm b}(\lambda): 
\left\{
\begin{array}{ll}
F_{\rm i}
\mapsto (\eta,0) &\mbox{case {\rm (i)}},   \\
F_{\rm ii}
\mapsto
(\eta,(\gamma\times\omega)\times\gamma)\qquad &\mbox{case {\rm (ii)}}, \\
F\mapsto (\eta,\omega) &\mbox{case {\rm (iii)}},
\end{array}
\right.
\label{rigid}
\end{equation}
which provides a unique solution to \eqref{solid}, \eqref{solid-ii} and \eqref{solid-i}, respectively, 
for $\lambda\in \mathbb C\setminus (-\infty,0)$ and $F\in L^q_{{\rm R},[d]}(\mathbb R^3)$, where 
$F_{\rm i}$ and $F_{\rm ii}$ are given by \eqref{data-modi}.

The following proposition provides us with the regularity of the inverse matrices 
in \eqref{rigid-inv}--\eqref{rigid-inv-ii}
near $\lambda=0$ along the imaginary axis.
\begin{proposition}
Suppose that $\alpha\in C(\partial\Omega)$ is nonnegative. 
Set
\[
K(\tau)=\frac{d}{d\tau}(i\tau+T_{i\tau})^{-1}, \qquad \tau\in\mathbb R\setminus\{0\}.
\]
For the case \emph{(iii)} in Theorem \ref{main},
there is a constant $C=C(\alpha,\Omega)>0$ such that
\begin{equation}
\sup_{\lambda\in\overline{\mathbb C_+}, |\lambda|\leq 4}
|(\lambda+T_\lambda)^{-1}|_{\mathbb C^{6\times 6}}
\leq C
\label{auxi-0}
\end{equation}
\begin{equation}
\int_{-4}^4 |K(\tau)|_{\mathbb C^{6\times 6}}\,d\tau\leq C
\label{auxi-2}
\end{equation}
\begin{equation}
\sup_{|\tau|\leq 3}\big|\big(i(\tau+h)+T_{i(\tau+h)}\big)^{-1}-(i\tau+T_{i\tau})^{-1}\big|_{\mathbb C^{6\times 6}}\leq C|h|^{1/2}
\label{auxi-3}
\end{equation}
\begin{equation}
\int_{-2}^2 |K(\tau+h)-K(\tau)|_{\mathbb C^{6\times 6}}\,d\tau
\leq C|h|^{1/2}
\label{auxi-4}
\end{equation}
for all $h\in\mathbb R$ with $|h|\leq 1$.
\label{prop-auxi}

Let $E$ and $E_\gamma$ be the subspaces of $\mathbb C^3\times\mathbb C^3$ given by \eqref{auxi-sp}, and set
\[
K_{\rm i}(\tau)=\frac{d}{d\tau}\big[(i\tau+T_{i\tau})|_E\big]^{-1}=K(\tau)|_E, \qquad
K_{\rm ii}(\tau)=\frac{d}{d\tau}\big[(i\tau+T_{i\tau})|_{E_\gamma}\big]^{-1}=K(\tau)|_{E_\gamma} 
\]
for $\tau\in\mathbb R\setminus\{0\}$.
For the case \emph{(i)}, the operators
$\big[(\lambda+T_\lambda)|_E\big]^{-1}$ and $K_{\rm i}(\tau)$ enjoy the same estimates on 
${\mathcal L}(E,E)$
as in \eqref{auxi-0}--\eqref{auxi-4}.
Likewise, for the case \emph{(ii)}, the operators
$\big[(\lambda+T_\lambda)|_{E_\gamma}\big]^{-1}$ and $K_{\rm ii}(\tau)$ fulfill those on 
${\mathcal L}(E_\gamma,E_\gamma)$.
\end{proposition}

\begin{proof}
Consider the case (iii).
Owing to \eqref{T-conti},
by the standard argument using the Neumann series, 
the map
\[
\{\lambda\in\overline{\mathbb C_+};\; |\lambda|\leq 4\}\ni \lambda\mapsto
(\lambda+T_\lambda)^{-1}
\]
is continuous; 
actually, this is true even at $\lambda=0$ because the cases (i) and (ii) are ruled out.
In fact, since
\[
(\lambda+T_\lambda)^{-1}=T_0^{-1}\big[1+(\lambda+T_\lambda-T_0)T_0^{-1}\big]^{-1}
=T_0^{-1}\sum_{k=0}^\infty\big\{-(\lambda+T_\lambda-T_0)T_0^{-1}\big\}^k
\]
for $\lambda\in \overline{\mathbb C_+}$ with $|\lambda|$ being small enough, we have
\[
(\lambda+T_\lambda)^{-1}-T_0^{-1}
=T_0^{-1}\big\{-(\lambda+T_\lambda-T_0)\big\}T_0^{-1}\sum_{k=0}^\infty\big\{-(\lambda+T_\lambda-T_0)T_0^{-1}\big\}^k,
\]
which implies
\[
|(\lambda+T_\lambda)^{-1}-T_0^{-1}|_{\mathbb C^{6\times 6}}\leq C\big(|\lambda|+|T_\lambda-T_0|_{\mathbb C^{6\times 6}}\big),
\]
yielding the continuity of $(\lambda+T_\lambda)^{-1}$ at $\lambda=0$. 
Obviously, the same thing holds at every point on $\overline{\mathbb C_+}$.
We thus obtain \eqref{auxi-0}.

Since 
\begin{equation}
K(\tau)=-(i\tau+T_{i\tau})^{-1}\big(i+L(\tau)
\big)(i\tau+T_{i\tau})^{-1},
\label{K-rep}
\end{equation}
where $L(\tau)=\frac{d}{d\tau}T_{i\tau}$,
we infer \eqref{auxi-2} 
from \eqref{auxi-5} along with 
\eqref{auxi-0}.
By virtue of \eqref{auxi-conti} together with \eqref{auxi-0},
we at once find \eqref{auxi-3} from
\begin{equation}
\begin{split}
&\quad\big(i(\tau+h)+T_{i(\tau+h)}\big)^{-1}-(i\tau+T_{i\tau})^{-1}  \\
&=-\big(i(\tau+h)+T_{i(\tau+h)}\big)^{-1}\big(ih+T_{i(\tau+h)}-T_{i\tau}\big)(i\tau+T_{i\tau})^{-1}.
\end{split}
\label{inv-diff}
\end{equation}
In view of \eqref{K-rep}, we collect 
\eqref{auxi-5}, \eqref{auxi-6}
\eqref{auxi-0} and \eqref{auxi-3}
to infer \eqref{auxi-4}.

We turn to the case (ii).
Then we have the inverse $\big[T_0|_{E_\gamma}\big]^{-1}: E_\gamma\to E_\gamma$ 
by \eqref{cond-uni}--\eqref{range}, which implies
\begin{equation}
\sup_{\lambda\in \overline{\mathbb C_+}, |\lambda|\leq 4}
\left|\big[(\lambda+T_\lambda)|_{E_\gamma}\big]^{-1}\right|_{{\mathcal L}(E_\gamma,E_\gamma)}   
\leq C
\label{auxi-0-ii}
\end{equation}
by the same reasoning as in the case (iii).
Similarly to \eqref{K-rep} and \eqref{inv-diff}, we have
\begin{equation*}
K_{\rm ii}(\tau)=-\big[(i\tau+T_{i\tau})|_{E_\gamma}\big]^{-1}\big(i+L(\tau)
\big)\big[(i\tau+T_{i\tau})|_{E_\gamma}\big]^{-1}
\end{equation*}
and
\begin{equation*}
\begin{split}
&\quad\left[\big(i(\tau+h)+T_{i(\tau+h)}\big)\big|_{E_\gamma}\right]^{-1}-\left[(i\tau+T_{i\tau})\big|_{E_\gamma}\right]^{-1}   \\
&=-\left[\big(i(\tau+h)+T_{i(\tau+h)}\big)\big|_{E_\gamma}\right]^{-1}\big(ih+T_{i(\tau+h)}-T_{i\tau}\big)
\left[(i\tau+T_{i\tau})\big|_{E_\gamma}\right]^{-1}.
\end{split}
\end{equation*}
With these at hand, we gather \eqref{auxi-0-ii} and \eqref{auxi-5}--\eqref{auxi-6}
to obtain the corresponding estimates on ${\mathcal L}(E_\gamma,E_\gamma)$
to \eqref{auxi-2}--\eqref{auxi-4}.
The similar argument is applied to the case (i) as well.
The proof is complete.
\end{proof}
\begin{corollary}
Suppose that $\alpha\in C(\partial\Omega)$ is nonnegative.
Let $q\in (1,\infty)$ and $d\geq 3$.
For
$\lambda\in \mathbb C\setminus (-\infty,0)$,
let ${\mathcal Q}_{\rm b}(\lambda)$ 
be the operator defined by \eqref{rigid}.
We set
\[
{\mathscr R}_{\rm b}(\tau)=
\partial_\tau{\mathcal Q}_{\rm b}(i\tau), \qquad \tau\in\mathbb R\setminus\{0\},
\]
and 
\[
{\mathcal Y}=
\left\{
\begin{array}{ll}
{\mathcal L}(L^q_{[d]}(\Omega)\times {\rm RM}_{(\kappa,0)\in E}, E) &\mbox{\rm case (i)}, \\
{\mathcal L}(L^q_{[d]}(\Omega) \times {\rm RM}_{(\kappa,\mu)\in E_\gamma}, E_\gamma) \qquad&\mbox{\rm case (ii)}, \\
{\mathcal L}(L^q_{{\rm R},[d]}(\mathbb R^3), \mathbb C^3\times \mathbb C^3) &\mbox{\rm case (iii)},
\end{array}
\right.
\]
for three cases in Theorem \ref{main},
where
\begin{equation}
{\rm RM}_{(\kappa,\mu)\in G}:=\{\kappa+\mu\times x;\; (\kappa,\mu)\in G\}
\label{spe-rigid}
\end{equation}
for $G\subset \mathbb C^3\times\mathbb C^3$ with $G=E, \,E_\gamma$ given by \eqref{auxi-sp}.
Then there is a constant $C=C(\alpha,q,d,\Omega)>0$ such that
\begin{equation}
\sup_{\lambda\in\overline{\mathbb C_+}, |\lambda|\leq 4}
\|{\mathcal Q}_{\rm b}(\lambda)\|_{\mathcal Y}
\leq C
\label{rigid-0}
\end{equation}
\begin{equation}
\int_{-4}^4\|{\mathscr R}_{\rm b}(\tau)\|_{\mathcal Y}\,d\tau\leq C
\label{rigid-1}
\end{equation}
\begin{equation}
\sup_{|\tau|\leq 3}\|{\mathcal Q}_{\rm b}(i(\tau+h))-{\mathcal Q}_{\rm b}(i\tau)\|_{\mathcal Y}
+\int_{-2}^2\|{\mathscr R}_{\rm b}(\tau+h)-{\mathscr R}_{\rm b}(\tau)\|_{\mathcal Y}\,d\tau\leq C|h|^{1/2}
\label{rigid-2}
\end{equation}
for all $h\in\mathbb R$ with $|h|\leq 1$.
\label{cor-body}
\end{corollary}

\begin{proof}
By the trace estimate and by collecting  \eqref{ext-0}--\eqref{ext-3} for $(u_2,p_2)$ as well as \eqref{auxi-0}--\eqref{auxi-4},
we are led to \eqref{rigid-0}--\eqref{rigid-2}.
\end{proof}

Corollary \ref{cor-body} gives us 
the desired regularity for the rigid motion of the body.
It remains to show the one for the velocity $u$ of the fluid to the 
fluid-structure resolvent systems \eqref{FS-1}--\eqref{FS-2}, 
\eqref{FS-ii-1}--\eqref{FS-ii-2} and \eqref{FS-i-1}--\eqref{FS-i-2}.
To this end, 
for $\lambda\in \mathbb C\setminus (-\infty,0)$,
we introduce the operator ${\mathcal Q}_{\rm f}(\lambda)$ 
by
\begin{equation}
{\mathcal Q}_{\rm f}(\lambda): 
\widetilde F
\mapsto u,
\label{flow}
\end{equation}
where $u$ is determined by \eqref{decouple}--\eqref{decouple2} through the lifting function 
\eqref{lift}/\eqref{lift2} and $\widetilde F$ is given by \eqref{IC-modi}--\eqref{data-modi}.
Then the regularity of ${\mathcal Q}_{\rm f}(\lambda)$ 
near $\lambda=0$ reads
\begin{proposition}
Suppose that $\alpha\in C(\partial\Omega)$ is nonnegative.
Let $q\in (1,\infty)$, $d\geq 3$ and $R\geq 3$.
For $\lambda\in \mathbb C\setminus (-\infty,0)$,
let ${\mathcal Q}_{\rm f}(\lambda)$ be the operator \eqref{flow}.
We set
\begin{equation}
{\mathscr R}_{\rm f}(\tau)=\partial_\tau{\mathcal Q}_{\rm f}(i\tau), \qquad \tau\in \mathbb R\setminus\{0\}
\label{diff-flow}
\end{equation}
and 
\begin{equation*}
{\mathcal X}=
\left\{
\begin{array}{ll}
{\mathcal L}(L^q_{[d]}(\Omega)\times {\rm RM}_{(\kappa,0)\in E}, W^{2,q}(\Omega_R)) &\mbox{\rm case (i)}, \\
{\mathcal L}(L^q_{[d]}(\Omega) \times {\rm RM}_{(\kappa,\mu)\in E_\gamma}, W^{2,q}(\Omega_R)) \qquad&\mbox{\rm case (ii)}, \\
{\mathcal L}(L^q_{{\rm R},[d]}(\mathbb R^3), W^{2,q}(\Omega_R)) &\mbox{\rm case (iii)},
\end{array}
\right.
\]
for three cases in Theorem \ref{main},
where $E$ and $E_\gamma$ are given by \eqref{auxi-sp}, see also \eqref{spe-rigid}.
Then there is a constant $C=C(\alpha,q,d,R,\Omega)>0$ such that
%
\begin{equation}
\sup_{\lambda\in\overline{\mathbb C_+}, |\lambda|\leq 4}
\|{\mathcal Q}_{\rm f}(\lambda)\|_{\mathcal X}\leq C
\label{flow-0}
\end{equation}
\begin{equation}
\int_{-4}^4\|{\mathscr R}_{\rm f}(\tau)\|_{\mathcal X}\,d\tau\leq C
\label{flow-1}
\end{equation}
\begin{equation}
\sup_{|\tau|\leq 3}\|{\mathcal Q}_{\rm f}(i(\tau+h))-{\mathcal Q}_{\rm f}(i\tau)\|_{\mathcal X}
+\int_{-2}^2\|{\mathscr R}_{\rm f}(\tau+h)-{\mathscr R}_{\rm f}(\tau)\|_{\mathcal X}\,d\tau\leq C|h|^{1/2}
\label{flow-2}
\end{equation}
for all $h\in\mathbb R$ with $|h|\leq 1$.
\label{prop-flow}
\end{proposition}

\begin{proof}
In view of \eqref{decouple}--\eqref{decouple2} and 
\[
\partial_\tau u=\partial_\tau\ell+V_u(\tau)
(\Delta\ell-i\tau\ell +f)+S_u(i\tau)
(\Delta\partial_\tau\ell-i\tau\partial_\tau\ell-i\ell)
\]
along the imaginary axis $\lambda=i\tau$ with \eqref{ext-sol2} and the abbreviation
\[
\ell=\left\{
\begin{array}{ll}
\ell(\eta,0)&\mbox{case (i)}, \\
\ell(\eta,(\gamma\times\omega)\times\gamma)\qquad&\mbox{case (ii)},  \\
\ell(\eta,\omega)&\mbox{case (iii)},
\end{array}
\right.
\]
where $\ell\in L^q_{[3]}(\Omega)$,
the proof of \eqref{flow-0}--\eqref{flow-2} is accomplished by 
Proposition \ref{prop-ext}, \eqref{lift-est} and
Corollary \ref{cor-body}.
\end{proof}

We combine the result on the rigid motion (Corollary \ref{cor-body}) with the one on the flow (Proposition \ref{prop-flow})
to furnish Corollary \ref{cor-reg-resol} below.
To describe the consequence, let us
introduce the operator ${\mathcal Q}(\lambda)$ by
\begin{equation}
{\mathcal Q}(\lambda):
\widetilde F\mapsto (\lambda+A)^{-1}\mathbb P\widetilde F=
\left\{
\begin{array}{ll}
{\mathcal Q}_{\rm i}(\lambda)F_{\rm i} & \mbox{\rm case (i)}, \\
{\mathcal Q}_{\rm ii}(\lambda)F_{\rm ii} & \mbox{\rm case (ii)}, \\
(\lambda+A)^{-1}\mathbb P F \qquad& \mbox{\rm case (iii)},
\end{array}
\right.
\label{mix}
\end{equation}
for $\lambda\in \mathbb C\setminus (-\infty,0)$ and $F\in L^q_{{\rm R},[d]}(\mathbb R^3)$,
where $\widetilde F$,
${\mathcal Q}_{\rm i}(\lambda)F_{\rm i}$ and ${\mathcal Q}_{\rm ii}(\lambda)F_{\rm ii}$ are given by 
\eqref{IC-modi}--\eqref{data-modi} and \eqref{resol-rep2},
see also \eqref{resol-rep} and \eqref{resol-rep3}.
\begin{corollary}
Suppose that $\alpha\in C(\partial\Omega)$ is nonnegative.
Let $q\in (1,\infty)$, $d\geq 3$ and $R\geq 3$.
For $\lambda\in\mathbb C\setminus (-\infty,0)$, 
let ${\mathcal Q}(\lambda)$ be the operator defined by \eqref{mix}.
We set
\begin{equation}
{\mathscr R}(\tau)=\partial_\tau {\mathcal Q}(i\tau), \qquad \tau\in\mathbb R\setminus \{0\}
\label{diff-mix}
\end{equation}
and
\[
{\mathcal Z}=
\left\{
\begin{array}{ll}
{\mathcal L}(L^q_{[d]}(\Omega)\times {\rm RM}_{(\kappa,0)\in E}, L^q(B_R))&\mbox{\rm case (i)}, \\
{\mathcal L}(L^q_{[d]}(\Omega)\times {\rm RM}_{(\kappa,\mu)\in E_\gamma}, L^q(B_R))\qquad&\mbox{\rm case (ii)}, \\
{\mathcal L}(L^q_{{\rm R},[d]}(\mathbb R^3), L^q(B_R))&\mbox{\rm case (iii)},
\end{array}
\right.
\]
for three cases in Theorem \ref{main}, where $E$ and $E_\gamma$ are given by \eqref{auxi-sp}, see also \eqref{spe-rigid}.
Then there is a constant $C=C(\alpha,q,d,R,\Omega)>0$ such that
\begin{equation}
\sup_{\lambda\in\overline{\mathbb C_+}, |\lambda|\leq 4}
\|{\mathcal Q}(\lambda)\|_{\mathcal Z}\leq C
\label{mix-0}
\end{equation}
\begin{equation}
\int_{-4}^4\|{\mathscr R}(\tau)\|_{\mathcal Z}\,d\tau\leq C
\label{mix-1}
\end{equation}
\begin{equation}
\sup_{|\tau|\leq 3}\|{\mathcal Q}(i(\tau+h))-{\mathcal Q}(i\tau)\|_{\mathcal Z}
+\int_{-2}^2\|{\mathscr R}(\tau+h)-{\mathscr R}(\tau)\|_{\mathcal Z}\,d\tau\leq C|h|^{1/2}
\label{mix-2}
\end{equation}
for all $h\in\mathbb R$ with $|h|\leq 1$.
\label{cor-reg-resol}
\end{corollary}
\begin{remark}
Unlike the no-slip case (\cite{EMT23, H24}),
the fluid-structure resolvent $(\lambda+A)^{-1}\mathbb P\widetilde F$ 
can no longer, in general, bolong to $W^{1,q}(B_R)$ across the boundary $\partial\Omega$ under the Navier condition,
so that we have no chance to replace $L^q(B_R)$ by $W^{1,q}(B_R)$ in \eqref{mix-0}--\eqref{mix-2}.  
\label{rem-reg-resol}
\end{remark}

\section{Local energy decay}
\label{led}

This section is devoted to a crucial step on local energy decay estimates of the fluid-structure semigroup.
Based on the spectral analysis developed in the previous section, such estimates for compactly supported data are deduced
in subsection \ref{sub-led} after some remarks on the semigroup in subsection \ref{sub-rem-sg}.
Successively, 
those for general data are established in subsection \ref{sub-led2}. 
There is an alternative route toward the former with less rate of decay, that is discussed in subection \ref{alter};
then, the latter procedure is the same.

\subsection{Remarks on the fluid-structure semigroup}
\label{sub-rem-sg}

As we have observed in subsection \ref{axisym}, a part of the fluid-structure semigroup does not decay
under the full slip condition (case $\alpha\equiv0$) when the rigid body ${\mathscr B}$ is axially symmetric. 
By taking into account \eqref{const-ii} for the case (ii) and \eqref{const-i} for the case (i) in Theorem \ref{main},
the object that we should discuss is
\begin{equation}
\begin{split}
{\mathcal S}(t)F
&:=\widetilde{e^{-tA}F}=e^{-tA}F-
\left\{
\begin{array}{ll}
(\omega(t)\times x)\chi_{\mathscr B}& \mbox{case (i)} \\
\{(\gamma\cdot\omega(t))\gamma\times x\}\chi_{\mathscr B}\qquad &\mbox{case (ii)}  \\
0&\mbox{case (iii)}
\end{array}
\right.  \\
&=e^{-tA}F-\left\{
\begin{array}{ll}
(\mu\times x)\chi_{\mathscr B}&\mbox{case (i)} \\
\{(\gamma\cdot\mu)\gamma\times x\}\chi_{\mathscr B}\qquad&\mbox{case (ii)} \\
0&\mbox{case (iii)}
\end{array}
\right\}
=e^{-tA}\widetilde F
\end{split}
\label{obje-0}
\end{equation}
for $F\in X_q(\mathbb R^3)$
with $i(F)=(f,\kappa,\mu)$, see \eqref{mono},
where $(u,\eta,\omega)=i(e^{-tA}F)$, while
$\widetilde{(\cdot)}$ is defined by \eqref{IC-modi}--\eqref{data-modi}. 
The last equality ${\mathcal S}(t)F=e^{-tA}\widetilde F$ in \eqref{obje-0} is readily seen from
\begin{equation}
\begin{split}
&\mbox{case (i)}\qquad
(\mu\times x)\chi_{\mathscr B}\in X_q(\mathbb R^3),\quad
e^{-tA}\big[(\mu\times x)\chi_{\mathscr B}\big]=(\mu\times x)\chi_{\mathscr B},   \\
&\mbox{case (ii)}\qquad
(\gamma\times x)\chi_{\mathscr B}\in X_q(\mathbb R^3), \quad
e^{-tA}\big[
(\gamma\times x)\chi_{\mathscr B}\big]=
(\gamma\times x)\chi_{\mathscr B}, 
\end{split}
\label{S-rela2}
\end{equation}
for all $\mu\in\mathbb R^3$ and $q\in (1,\infty)$
by noting that, for the case (ii),
$(\gamma\times x)\chi_{\mathscr B}$ fulfills \eqref{normal-va} because of \eqref{sym-1}.

For later use, let us verify the semigroup property
\begin{equation}
{\mathcal S}(t){\mathcal S}(s)={\mathcal S}(t+s), \qquad t,\, s\geq 0
\label{S-semi}
\end{equation}
on $X_q(\mathbb R^3)$ even for the cases (i) and (ii).
In fact, since
we see from \eqref{obje-0} that
\begin{equation}
\begin{split}
&\mbox{case (i)}\qquad
e^{-tA}F={\mathcal S}(t)F\qquad\mbox{if $i(F)=(f,\kappa,0)$},  \\
&\mbox{case (ii)}\qquad
e^{-tA}F={\mathcal S}(t)F\qquad\mbox{if $i(F)=(f,\kappa,\mu)$ with $\gamma\cdot\mu=0$},
\end{split}
\label{S-rela}
\end{equation}
we find
\begin{equation*}
\begin{split}
{\mathcal S}(t){\mathcal S}(s)F
&={\mathcal S}(t)\big[e^{-sA}F-\{(\gamma\cdot\mu)\gamma\times x\}\chi_{\mathscr B}\big]  \\
&=e^{-tA}\big[e^{-sA}F-\{(\gamma\cdot\mu)\gamma\times x\}\chi_{\mathscr B}\big]  \\
&=e^{-(t+s)A}F-\{(\gamma\cdot\mu)\gamma\times x\}\chi_{\mathscr B}  \\
&={\mathcal S}(t+s)F
\end{split}
\end{equation*}
for all $F\in X_q(\mathbb R^3)$, yielding \eqref{S-semi} for the case (ii).
The other case (i) is verified in the same manner.

In what follows we will use the relation
\begin{equation}
\widetilde{{\mathcal S}(t)F}={\mathcal S}(t)F={\mathcal S}(t)\widetilde F
\label{S-rela3}
\end{equation}
for all $F\in X_q(\mathbb R^3)$,
which follows from \eqref{obje-0} and the definition \eqref{data-modi} of $\widetilde{(\cdot)}$.
It is also noted that 
\begin{equation}
\langle {\mathcal S}(t)F,G\rangle_{\mathbb R^3,\rho}
=\langle F,{\mathcal S}(t)G\rangle_{\mathbb R^3,\rho}
\label{S-sym}
\end{equation}
for all $F\in X_q(\mathbb R^3)$ and $G\in X_{q^\prime}(\mathbb R^3)$.
In fact, for the case (ii), it follows from \eqref{pair} and \eqref{obje-0} that
\[
\langle {\mathcal S}(t)F,G\rangle_{\mathbb R^3,\rho}
=\langle e^{-tA}F,G\rangle_{\mathbb R^3,\rho}
-\big(J(\gamma\cdot\mu_F)\gamma\big)\cdot\mu_G,
\]
where $(f,\kappa_F,\mu_F)=i(F)$ and $(g,\kappa_G,\mu_G)=i(G)$.
Since we know \eqref{sym-semi} 
and since \eqref{sym-inertia} together with symmetry of $J$ implies
\begin{equation*}
\big(J(\gamma\cdot\mu_F)\gamma\big)\cdot\mu_G
=(\gamma\cdot\mu_F)\gamma\cdot(J\mu_G)
=c(\gamma\cdot\mu_F)(\gamma\cdot\mu_G)
=\mu_F\cdot \big(J(\gamma\cdot\mu_G)\gamma\big),
\end{equation*}
we are led to \eqref{S-sym} for the case (ii).
The other case (i) is discussed similarly.

Let us deduce regularizing estimates of ${\mathcal S}(t)$ near the initial time $t=0$:
\begin{equation}
\|{\mathcal S}(t)F\|_{q,\mathbb R^3}
+t^{j/2}\|\nabla^j{\mathcal S}(t)F\|_{q,\Omega}\leq C\|\widetilde F\|_{q,\mathbb R^3}
\label{reg-S0}
\end{equation}
\begin{equation}
\left\{
\begin{array}{ll}
t^{3/2q}\|{\mathcal S}(t)F\|_{\infty,\Omega}\leq C\|\widetilde F\|_{q,\mathbb R^3} 
&(q>3/2)  \\
t^{3/2q+1/2}\|\nabla{\mathcal S}(t)F\|_{\infty,\Omega}\leq C\|\widetilde F\|_{q,\mathbb R^3} \qquad &(q>3)
\end{array}
\right.
\label{reg-S}
\end{equation}
for all $t\in (0,4)$, $F\in X_q(\mathbb R^3)$ and  $j\in \{1,2\}$, 
where $q\in (1,\infty)$ for \eqref{reg-S0}, while $q$ is specified as above for \eqref{reg-S} (this restriction can be
removed by use of the semigroup property, but it will be done in section \ref{pr-main}).
In fact, by \eqref{bdd-anal} with $k=0,\,1$, we have
\[
\|{\mathcal S}(t)F\|_{q,\mathbb R^3}+t\|A{\mathcal S}(t)F\|_{q,\mathbb R^3}
=\|e^{-tA}\widetilde F\|_{q,\mathbb R^3}+t\|Ae^{-tA}\widetilde F\|_{q,\mathbb R^3}\leq C\|\widetilde F\|_{q,\mathbb R^3} 
\]
for all $t>0$ and $F\in X_q(\mathbb R^3)$,
which implies that
\begin{equation}
\left.
\begin{array}{ll}
\mbox{case (i)}&|\eta(t)|  \\
\mbox{case (ii)}&|\eta(t)|+|\gamma\times \omega(t)| \\
\mbox{case (iii)}\qquad&|\eta(t)|+|\omega(t)|
\end{array}
\right\}
\leq C\|\widetilde F\|_{q,\mathbb R^3}
\label{reg-rigid}
\end{equation}
for all $t>0$, where $(u,\eta,\omega)=i(U(t))$ and $U(t)=e^{-tA}F$, because
\begin{equation}
\big({\mathcal S}(t)F\big)|_{\mathcal B}
=\left\{
\begin{array}{ll}
\eta(t)\chi_{\mathscr B}&\mbox{case (i)},  \\
\big[\eta(t)+\{(\gamma\times \omega(t))\times \gamma\}\times x\big]\chi_{\mathscr B}\qquad &\mbox{case (ii)},  \\
(\eta(t)+\omega(t)\times x)\chi_{\mathscr B}&\mbox{case (iii)},
\end{array}
\right.
\label{S-body}
\end{equation}
and that
\begin{equation}
t^{j/2}\|u(t)\|_{W^{j,q}(\Omega)}\leq C\|\widetilde F\|_{q,\mathbb R^3}
\label{reg-S2}
\end{equation}
for all $t\in (0,4)$
and $j\in \{1,2\}$ by virtue of \eqref{ellip-0}.
We thus obtain \eqref{reg-S0}.
It follows from the interpolation inequality
due to Crispo and Maremonti \cite[Theorem 2.1]{CM04}
(although the standard Gagliardo-Nirenberg inequality is enough for our aim)
that
\[
\|u(t)\|_{\infty,\Omega}\leq C
\|\nabla^2 u(t)\|_{q,\Omega}^{3/2q}
\|u(t)\|_{q,\Omega}^{1-3/2q}\qquad (q>3/2)
\]
and that
\[
\|\nabla u(t)\|_{\infty,\Omega}\leq C
\|\nabla^2u(t)\|_{q,\Omega}^{3/q}
\|\nabla u(t)\|_{q,\Omega}^{1-3/q}\qquad (q>3).
\]
Then \eqref{reg-S2} combined with them leads us to \eqref{reg-S}.
Further regularizing estimates will be studied along with decay properties in section \ref{pr-main}.

In the next subsection we take the initial data $F$ from $L^q_{{\rm R},[d]}(\mathbb R^3)$ 
with $(f,\kappa,\mu)=i(F)$, see \eqref{data-sp2},
and study the decay property of
\begin{equation}
{\mathcal S}(t)\mathbb PF=e^{-tA}\mathbb PF-\left\{
\begin{array}{ll}
(\mu\times x)\chi_{\mathscr B}&\mbox{case (i)} \\
\{(\gamma\cdot\mu)\gamma\times x\}\chi_{\mathscr B}\qquad&\mbox{case (ii)} \\
0&\mbox{case (iii)}
\end{array}
\right\}=
e^{-tA}\mathbb P\widetilde F,
\label{obje}
\end{equation}
where we have used
\begin{equation}
\left\{
\begin{array}{ll}
\mu^\prime=\mu&\mbox{case (i)}, \\
\gamma\cdot\mu^\prime=\gamma\cdot\mu\qquad&\mbox{case (ii)},
\end{array}
\right.
\label{proj-ang}
\end{equation}
for $i(\mathbb PF)=(f^\prime,\kappa^\prime,\mu^\prime)$, 
which together with $\mbox{\eqref{S-rela2}}_1$ yields
\begin{equation}
\widetilde{\mathbb PF}=\mathbb P\widetilde F
\label{proj-ang2}
\end{equation}
see \eqref{IC-modi}--\eqref{data-modi}.
Indeed, it follows from \eqref{sym-1} and \eqref{sym-inertia} 
that the augular velocity $\omega$ for each element of $Z_q(\mathbb R^3)$,
see \eqref{anni}, satisfies
\begin{equation}
\begin{split}
&\mbox{case (i)}\qquad
\omega=-J^{-1}\int_{\partial\Omega}x\times (p\nu)\,d\sigma=0, \\
&\mbox{case (ii)}\qquad
\gamma\cdot\omega=\gamma\cdot\left(-J^{-1}\int_{\partial\Omega}x\times (p\nu)\,d\sigma\right)=0,
\end{split}
\label{Z-ang}
\end{equation}
which imply \eqref{proj-ang}.

\subsection{Estimates for compactly supported data}
\label{sub-led}

Let us start with the following lemma,
which gives us the relation between the regularity of functions and the rate of decay of those
(inverse) Fourier transform.
The origin is traced back to Shibata \cite{S83} and, later on, an elementary proof of his lemma is given by the first author \cite{H04}.
Estimate of the form \eqref{est-four} below is found in \cite{H16}; indeed, it would be easier to apply to various situations.
\begin{lemma}
[{\cite[Lemma 7.3]{H16}}]
Let $X$ be a Banach space with norm $\|\cdot\|$ and $w\in L^1(\mathbb R;X)$.
Set
\begin{equation}
W(t)=\int_{-\infty}^\infty e^{i\tau t}w(\tau)\,d\tau.
\label{fourier}
\end{equation}
Then the function $W(t)$ enjoys
\begin{equation}
\|W(t)\|\leq C\int_{-\infty}^\infty\left\|w\left(\tau+\frac{1}{t}\right)-w(\tau)\right\|\,d\tau
\label{est-four}
\end{equation}
for all $t\in\mathbb R\setminus \{0\}$.
\label{lem-four}
\end{lemma}

\begin{proof}
We give a quite short proof for completeness.
Let $h\in\mathbb R$ be a parameter, then we have
\[
W(t)=\int_{-\infty}^\infty e^{i(\tau+h)t}w(\tau+h)\,d\tau
=e^{iht}\left[\int_{-\infty}^\infty e^{i\tau t}\big(w(\tau+h)-w(\tau)\big)\,d\tau+W(t)\right].
\]
Given $t\in\mathbb R\setminus \{0\}$, we take $h=1/t$ to obtain \eqref{est-four}.
\end{proof}

In view of \eqref{resol-rep2}--\eqref{resol-rep3}
for the cases (i) and (ii) in Theorem \ref{main}, it is convenient to use the operator 
${\mathcal Q}(\lambda)$ given by \eqref{mix} in order to describe \eqref{obje}:
\begin{equation}
{\mathcal S}(t)\mathbb PF=e^{-tA}\mathbb P\widetilde F=
\frac{-1}{2\pi it}\int_{-\infty}^\infty e^{i\tau t}\partial_\tau {\mathcal Q}(i\tau)\widetilde F 
\,d\tau
\label{sg-rep}
\end{equation}
which is valid in $L^q(B_R)$ as well as in $W^{1,q}(\Omega_R)$,
see Remark \ref{rem-reg-resol}.
In fact, starting from the standard Dunford integral of 
${\mathcal Q}(\lambda)\widetilde F=(\lambda+A)^{-1}\mathbb P\widetilde F$,
in which the path of integration may hit the origin ($\lambda =0$) on account of the boundedness \eqref{flow-0} and \eqref{mix-0},
we perform integration by parts and then move
the path 
to the imaginary axis to derive \eqref{sg-rep},
that is actually summable near $\tau=0$ due to \eqref{flow-1} and \eqref{mix-1}.
The proof of the local energy decay property \eqref{1st-led} below is based on the represenation \eqref{sg-rep}.
The idea goes back to Kobayashi and Shibata \cite{KS98}, see also
\cite{ES04,H04,H16,H26,HS09}.
\begin{proposition}
Suppose that $\alpha\in C^1(\partial\Omega)$ is nonnegative.
Let $q\in (1,\infty)$, $d\geq 3$ and $R\geq 3$. 
Then there is a constant $C=C(\alpha,q,d,R,\Omega)>0$ such that
\begin{equation}
\|{\mathcal S}(t)\mathbb PF\|_{q,B_R}+\|{\mathcal S}(t)\mathbb PF\|_{W^{2,q}(\Omega_R)}
\leq Ct^{-3/2}
\|\widetilde F\|_{q,\mathbb R^3}
\label{1st-led}
\end{equation}
for all $t\geq 1$ 
and $F\in L^q_{{\rm R},[d]}(\mathbb R^3)$, see \eqref{data-sp2}, where $\widetilde F$ is given by \eqref{IC-modi}--\eqref{data-modi}.
\label{prop-led}
\end{proposition}

\begin{proof}
With use of a cut-off function $\psi_0\in C^\infty(\mathbb R;[0,1])$
fulfilling $\psi_0(\tau)=1$ for $|\tau|\leq 1$ and $\psi_0(\tau)=0$ for $|\tau|\geq 2$,
we divide the integral \eqref{sg-rep} into two parts
\begin{equation*}
\begin{split}
&{\mathcal S}(t)\mathbb PF=\frac{-1}{2\pi it}\Big(W_0(t)+W_\infty(t)\Big), \\  
&W_k(t)=\int_{-\infty}^\infty e^{i\tau t}\psi_k(\tau)
{\mathscr R}(\tau)\widetilde F\,d\tau, \qquad k\in \{0,\infty\},
\end{split}
\end{equation*}
where $\psi_\infty=1-\psi_0$, while ${\mathscr R}(\tau)$ is given by \eqref{diff-mix}.
It suffices to derive the rate $t^{-1/2}$ of decay of $W_k(t)$. 

Let us apply Lemma \ref{lem-four} to
\[
w(\tau)=\psi_0(\tau){\mathscr R}(\tau)\widetilde F, \qquad
X=W^{2,q}(\Omega_R), 
\]
Recall that ${\mathscr R}(\tau)\widetilde F={\mathscr R}_{\rm f}(\tau)\widetilde F$ 
on the fluid region $\Omega$, see \eqref{diff-flow}.
By \eqref{est-four} we infer
\begin{equation*}
\begin{split}
\|W_0(t)\|_{W^{2,q}(\Omega_R)}
&\leq C\int_{-\infty}^\infty\left|\psi_0\left(\tau+\frac{1}{t}\right)-\psi_0(\tau)\right|
\left\|{\mathscr R}_{\rm f}\left(\tau+\frac{1}{t}\right)\widetilde F
\right\|_{W^{2,q}(\Omega_R)}\,d\tau   \\
&\qquad +C\int_{-\infty}^\infty |\psi_0(\tau)|
\left\|{\mathscr R}_{\rm f}\left(\tau+\frac{1}{t}\right)\widetilde F
-{\mathscr R}_{\rm f}(\tau)\widetilde F
\right\|_{W^{2,q}(\Omega_R)}\,d\tau
\end{split}
\end{equation*}
for $t>0$. 
If in particular $t\geq 1$,
then it follows from \eqref{flow-1}--\eqref{flow-2} that
\[
\|W_0(t)\|_{W^{2,q}(\Omega_R)}\leq C\big(t^{-1}+t^{-1/2}\big)\|\widetilde F\|_{q,\mathbb R^3}.
\]
Here, $X=L^q(B_R)$ is also possible since we have \eqref{mix-1}--\eqref{mix-2} for ${\mathscr R}(\tau)$.
The heart of the proof is that the modulus $|h|^{1/2}$ 
of continuity in \eqref{flow-2}/\eqref{mix-2} brings us the rate $t^{-1/2}$
of decay in the second term above.

Let us turn to the other integral $W_\infty(t)$ away from $\tau=0$, which is treated easily by integration by parts;
in fact, we observe
\[
W_\infty(t)=\frac{-1}{it}\int_{-\infty}^\infty e^{i\tau t}\big\{\psi_\infty(\tau)\partial_\tau{\mathscr R}(\tau)\widetilde F
+\psi_\infty^\prime(\tau){\mathscr R}(\tau)\widetilde F\big\}\,d\tau.
\]
Thanks to \eqref{ellip-0},
we find 
\begin{equation*}
\begin{split}
&\quad \|\partial_\tau{\mathscr R}(\tau)\widetilde F\|_{W^{2,q}(\Omega)}
+\|\partial_\tau{\mathscr R}(\tau)\widetilde F\|_{q,\mathbb R^3}  \\
&=
2\left(\|(i\tau+A)^{-3}\mathbb P\widetilde F\|_{W^{2,q}(\Omega)}+\|(i\tau+A)^{-3}\mathbb P\widetilde F\|_{q,\mathbb R^3}\right)  \\
&\leq C\left(\|A(i\tau+A)^{-3}\mathbb P\widetilde F\|_{q,\mathbb R^3}+\|(i\tau+A)^{-3}\mathbb P\widetilde F\|_{q,\mathbb R^3}\right)   \\
&\leq C|\tau|^{-2}\|\widetilde F\|_{q,\mathbb R^3}
\end{split}
\end{equation*}
for $|\tau|\geq 1$.
Thus, the integral above is actually summable both in $W^{2,q}(\Omega_R)$ and $L^q(B_R)$
to derive the rate $t^{-1}$ of decay of $W_\infty(t)$.
The proof is complete.
\end{proof}

\subsection{Estimates for general data}
\label{sub-led2}

For the proof of
\eqref{decay-1} and \eqref{decay-2}, as in several pepers
\cite{ES05,EMT23,H04,H16,H20,H24,HS09,I89,KS98,MS97,SY08},
it would be standard to split those
into the estimate near the body 
and the one at large distance. 
This goes back to Iwashita \cite{I89} for the Stokes semigroup in exterior domains.
The present subsection studies the former, while the latter step will be established in the next section.
We need to adapt slightly but carefully the methodology to our circumstance in order that the demonstration covers all the cases
(i)--(iii) in Theorem \ref{main}.

We employ a cut-off procedure in both steps mentioned above, 
in which the Bogovskii operator (\cite{Bog79, BS90, Ga-book})
is needed to recover the solenoidal condition.
Consider the boundary value problem 
\[
\mbox{div $w$}=g\;\;\mbox{in $D$}, \qquad w|_{\partial D}=0,
\]
where $D$ is a bounded domain with Lipschitz boundary $\partial D$.
Let $1<q<\infty$.
Given $g\in L^q(D)$ with compatibility condition $\int_Dg\,dx=0$,
we have many solutions, among which a particular one discovered by Bogovskii \cite{Bog79}
is useful because of fine properties of his solution.
The solution operator which provides his solution, called the Bogovskii operator, is well defined as follows:
there is a linear operator $\mathbb B=\mathbb B_D: C_0^\infty(D)\to C_0^\infty(D)^3$ such that
\begin{equation}
\|\nabla^{k+1}\mathbb Bg\|_{q,D}\leq C\|\nabla^kg\|_{q,D}
\label{bog-est1}
\end{equation}
for all $q\in (1,\infty)$, integer $k\geq 0$ and $g\in C_0^\infty(D)$, where $C=C(D,q,k)>0$ is invariant with respect to dilation
of the domain $D$, and that
\begin{equation*}
\mbox{div $(\mathbb Bg)$}=g\qquad\mbox{if}\;
\int_Dg(x)\,dx=0.
\end{equation*}
By continuity, $\mathbb B$ extends uniquely to a bounded operator from $W^{k,q}_0(D)$ to $W^{k+1,q}_0(D)$.
In \cite[Theorem 2.5]{GHH06} Geissert, Heck and Hieber proved
\begin{equation}
\|\mathbb Bg\|_{q,D}\leq C\|g\|_{W^{1,q^\prime}(D)^*}
\label{bog-est2}
\end{equation}
for all $q\in (1,\infty)$ and $g\in W^{1,q^\prime}(D)^*$, where $1/q^\prime+1/q=1$.
This is not true if $W^{1,q^\prime}(D)^*$ is replaced by $W^{-1,q}(D)$, see \cite[Chapter III]{Ga-book}, who nevertheless showed
\begin{equation}
\|\mathbb B[\mbox{div $G$}]\|_{q,D}\leq C\|G\|_{q,D}
\label{bog-est3}
\end{equation}
for all $q\in (1,\infty)$ and $G\in L^q(D)^3$ satisfying $\nu\cdot G|_{\partial D}=0$ as well as $\mbox{div $G$}\in L^q(D)$,
see \cite[Theorem III.3.4]{Ga-book}.
Indeed, there is a term arising from cut-off procedure, for which one needs either \eqref{bog-est2} or \eqref{bog-est3}.

We start with the following proposition on decay estimates of $\nabla^j {\mathcal S}(t)F$ for $j=0,\, 1$.
The case $j=2$ is discussed in Proposition \ref{prop-hess-led}.
\begin{proposition}
Suppose that $\alpha\in C^1(\partial\Omega)$ is nonnegative.
Let $q\in (1,\infty)$ and $R\geq 3$.
Then there is a constant $C=C(\alpha,q,R,\Omega)>0$ such that
\begin{equation}
\|{\mathcal S}(t)F\|_{q,B_R}+\|{\mathcal S}(t)F\|_{W^{1,q}(\Omega_R)}\leq Ct^{-3/2q}\|\widetilde F\|_{q,\mathbb R^3}
\label{2nd-led}
\end{equation}
for all $t\geq 1$ and $F\in X_q(\mathbb R^3)$, where $\widetilde F$ is given by \eqref{IC-modi}--\eqref{data-modi}.
\label{prop-led2}
\end{proposition}

\begin{proof}
It suffices to show that the left-hand side of \eqref{2nd-led} is bounded from above by $Ct^{-3/2q}\|F\|_{q,\mathbb R^3}$
for all $F\in X_q(\mathbb R^3)$
since ${\mathcal S}(t)F={\mathcal S}(t)\widetilde F$, see \eqref{S-rela3}.

Given $F\in X_q(\mathbb R^3)$
with $(f,\kappa,\mu)=i(F)$, 
by using the heat semigroup $e^{t\Delta}$ on the whole space $\mathbb R^3$,
we consider
\begin{equation*}
\begin{split}
&U_0(t)=(1-\phi)u_0(t)+\mathbb B[u_0(t)\cdot\nabla \phi],  \\ 
&u_0(t)=e^{t\Delta}F:=\int_{\mathbb R^3}\Gamma(\cdot-y,t)F(y)\,dy, \\
\end{split}
\end{equation*}
where 
\[
\Gamma(x,t)=(4\pi t)^{-3/2}e^{-|x|^2/4t},
\]
$\phi$ is the cut-off function as in \eqref{cut}
and $\mathbb B=\mathbb B_D$
denotes the Bogovskii operator
in the domain $D=B_3\setminus \overline{B_1}$.
Note that $\int_D
u_0\cdot\nabla\phi\,dx=0$ follows from 
$\mbox{div $u_0$}=0$
and that $\mathbb B[u_0\cdot\nabla \phi]$ may be understood as its extension
defined on $\mathbb R^3$ by setting zero outside $D$.
By the well-known estimates of the heat semigroup along with \eqref{bog-est1}, we have
\begin{equation}
\|U_0(t)\|_{W^{1,q}(B_R)}\leq C\|u_0(t)\|_{W^{1,\infty}(\mathbb R^3)}\leq Ct^{-3/2q}\|F\|_{q,\mathbb R^3}
\label{est-app}
\end{equation}
for all $t\geq 1$.
Set
\[
V(t)=e^{-tA}F-U_0(t), \qquad
(v,\eta,\omega)=i(V(t))
\]
and let $p(t)$ be the pressure associated with $e^{-tA}F$. 
The rigid motion $\eta+\omega\times x$ associated with $V(t)$ is exactly the same as the one determined by $e^{-tA}F$ 
since $U_0(t)|_{\mathscr B}=0$.
In view of \eqref{obje-0}, we may focus on
\begin{equation}
\widetilde{V(t)}
=V(t)-
\left\{
\begin{array}{ll}
(\mu\times x)\chi_{\mathscr B}&\mbox{case (i)} \\
\{(\gamma\cdot\mu)\gamma\times x\}\chi_{\mathscr B}\qquad &\mbox{case (ii)}  \\
0&\mbox{case (iii)}
\end{array}
\right\}
={\mathcal S(t)}F-U_0(t),
\label{focus}
\end{equation}
see \eqref{data-modi}, since we know \eqref{est-app}.

We see that $v,\, p,\, \eta$ and $\omega$ should obey
\begin{equation*}
\begin{split}
&\partial_tv=\Delta v-\nabla p+G(t), \quad\mbox{div $v$}=0, \quad\mbox{in $\Omega\times (0,\infty)$}, \\
&\nu\cdot (v-\eta-\omega\times x)|_{\partial\Omega}=0, \quad
[(2\mathbb Dv)\nu]_\tau+\alpha (v-\eta-\omega\times x)_\tau|_{\partial\Omega}=0, \quad
\lim_{|x|\to\infty}v=0, \\
&m\frac{d\eta}{dt}+\int_{\partial\Omega}\mathbb S(v,p)\nu\,d\sigma=0, \quad
J\frac{d\omega}{dt}+\int_{\partial\Omega}x\times \mathbb S(v,p)\nu\,d\sigma=0, \\
&v(\cdot,0)=f^\prime, 
\quad \eta(0)=\kappa, \quad \omega(0)=\mu,
\end{split}
\end{equation*}
which is converted into
\begin{equation}
V(t)=
e^{-tA}F^\prime+\int_0^te^{-(t-s)A}G(s)\,ds,
\label{focus-int1}
\end{equation}
where
$F^\prime=f^\prime\chi_\Omega+(\kappa+\mu\times x)\chi_{\mathscr B}\in X_q(\mathbb R^3)$
with $f^\prime=\phi f-\mathbb B[f\cdot\nabla \phi]$, 
for which it follows from \eqref{1st-led} that
\begin{equation}
\|{\mathcal S}(t)F^\prime\|_{q,B_R}+\|{\mathcal S}(t)F^\prime\|_{W^{1,q}(\Omega_R)}\leq Ct^{-3/2}\|F\|_{q,\mathbb R^3}
\label{1st-apply1}
\end{equation}
for all $t\geq 1$, while the remaining term $G(t)$ is given by
\[
G(t)=-2\nabla\phi\cdot\nabla u_0(t)-(\Delta\phi)u_0(t)
-\mathbb B[\partial_tu_0(t)\cdot\nabla\phi]
+\Delta \mathbb B[u_0(t)\cdot\nabla\phi]
\]
that fulfills
\[
\mbox{div $G(t)$}=\Delta p=0, \qquad
\mathbb PG(t)=G(t), \qquad
\widetilde{G(t)}=G(t),
\]
see \eqref{data-modi}, on account of
$G(t)|_{\mathscr B}=0$, as well as
\begin{equation}
\begin{split}
\|G(t)\|_{q,\mathbb R^3}
&\leq C\big(\|u_0(t)\|_{W^{1,q}(D)}+\|\Delta u_0(t)\|_{W^{-1,q}(D)}\big)  \\
&\leq Ct^{-1/2}(1+t)^{-3/2q+1/2}\|F\|_{q,\mathbb R^3}
\end{split}
\label{est-gar}
\end{equation}
for all $t>0$
by use of \eqref{bog-est1}--\eqref{bog-est2};
indeed, we observe
\[
|\langle \partial_tu_0(t)\cdot\nabla\phi, \psi\rangle_D|
\leq \|\partial_tu_0(t)\|_{W^{-1,q}(D)}\|(\nabla\phi)\psi\|_{W^{1,q^\prime}_0(D)}
\leq C\|\partial_tu_0(t)\|_{W^{-1,q}(D)}\|\psi\|_{W^{1,q^\prime}(D)}
\]
for all $\psi\in W^{1,q^\prime}(D)$, which together with \eqref{bog-est2} implies
\[
\|\mathbb B[\partial_tu_0(t)\cdot\nabla\phi]\|_{q,D}\leq C\|\partial_tu_0(t)\|_{W^{-1,q}(D)}=
\|\Delta u_0(t)\|_{W^{-1,q}(D)}.
\]

By \eqref{focus} and by $e^{-(t-s)A}G(s)={\mathcal S}(t-s)G(s)$, which follows from $G(s)|_{\mathscr B}=0$,
the equation \eqref{focus-int1} is rewritten as 
\begin{equation}
\widetilde{V(t)}={\mathcal S}(t)F^\prime
+\int_0^t{\mathcal S}(t-s)G(s)\,ds.
\label{focus-int2}
\end{equation}
Since we know \eqref{1st-apply1}, our task is to show the desired estimate of the Duhamel term.
We combine \eqref{1st-led} with \eqref{est-gar} to infer
\begin{equation*}
\begin{split}
&\quad \int_0^t\|{\mathcal S}(t-s)G(s)\|_{q,B_R}\,ds
+\int_0^t\|{\mathcal S}(t-s)G(s)\|_{W^{1,q}(\Omega_R)}\,ds  \\
&\leq C\|F\|_{q,\mathbb R^3}\int_0^t(t-s)^{-1/2}(1+t-s)^{-1}s^{-1/2}(1+s)^{-3/2q+1/2}\,ds  \\
&\leq Ct^{-3/2q}\|F\|_{q,\mathbb R^3}
\end{split}
\end{equation*}
for all $t\geq 1$ readily by splitting the integral into $\int_0^{t/2}+\int_{t/2}^t$.
The proof is complete.
\end{proof}

By \eqref{mono} and \eqref{S-body}, 
as in the following corollary,
estimate of $\|{\mathcal S}(t)F\|_{q,B_R}$ 
in \eqref{2nd-led} along with \eqref{reg-rigid}
immediately implies \eqref{decay-body}.
The second corollary below provides us with estimates of $\partial_t{\mathcal S}(t)F$ and the associated pressure
when the initial data are taken from $D_q(A)$.
\begin{corollary}
Suppose that $\alpha\in C^1(\partial\Omega)$ is nonnegative.
Let $q\in (1,\infty)$.
For $F\in X_q(\mathbb R^3)$, we set $U(t)=e^{-tA}F$ and $(u,\eta,\omega)=i(U(t))$.
Then there is a constant $C=C(\alpha,q,\Omega)>0$ such that
\begin{equation}
\begin{split}
\left.
\begin{array}{rl}
\mbox{\rm case (i)}&\qquad|\eta(t)| \\
\mbox{\rm case (ii)}&\qquad|\eta(t)|+|\gamma\times\omega(t)| \\
\mbox{\rm case (iii)}&\qquad|\eta(t)|+|\omega(t)|
\end{array}
\right\}
\leq Ct^{-3/2q}\|\widetilde F\|_{q,\mathbb R^3}
\end{split}
\label{dacay-body2}
\end{equation}
for all $t\geq 1$ and $F\in X_q(\mathbb R^3)$, where $\widetilde F$ is given by \eqref{IC-modi}--\eqref{data-modi}, and thereby
\eqref{decay-body} in Theorem \ref{main}.
\label{cor-rigid}
\end{corollary}
\begin{corollary}
Suppose that $\alpha\in C^1(\partial\Omega)$ is nonnegative.
Let $q\in (1,\infty)$ and $R\geq 3$.
Set $U(t)={\mathcal S}(t)F$ with $F\in D_q(A)$, $u(t)=U(t)|_{\Omega}$ and single out the pressure $p(t)$ associated with $U(t)$ in such a way that
$\int_{\Omega_R}p(t)\,dx=0$.
Then there is a constant $C=C(\alpha,q,R,\Omega)>0$ such that
\begin{equation}
\|\partial_tu(t)\|_{q,\Omega_R}+\|u(t)\|_{W^{1,q}(\Omega_R)}+\|p(t)\|_{q,\Omega_R}
\leq C(1+t)^{-3/2q}\big(\|A\widetilde F\|_{q,\mathbb R^3}+\|\widetilde F\|_{q,\mathbb R^3}\big)
\label{sm-led}
\end{equation}
for all $t\geq 0$ and $F\in D_q(A)$, where $\widetilde F$ is given by \eqref{IC-modi}--\eqref{data-modi}. 
\label{cor-led}
\end{corollary}

\begin{proof}
First of all, in view of 
\eqref{ig}, 
$F\in D_q(A)$ if and only if $\widetilde F\in D_q(A)$; for, $\alpha\equiv 0$ for the cases (i) and (ii).
Then we have $A\widetilde F=AF$.
Moreover, we see from \eqref{ig} again that
%
\begin{equation}
\widetilde{AF}=AF
\label{til-eq}
\end{equation}
for all $F\in D_q(A)$.
In fact,
it follows from \eqref{sym-defo}, \eqref{sym-1} and \eqref{sym-inertia} 
that the angular velocity $\omega$ for ${\mathscr A}F$ fulfills
\begin{equation*}
\begin{split}
&\mbox{case (i)}\qquad
\omega=J^{-1}\int_{\partial\Omega}y\times (2\mathbb Df)\nu\,d\sigma_y=0,   \\
&\mbox{case (ii)}\qquad
\gamma\cdot\omega=\gamma\cdot\left(
J^{-1}\int_{\partial\Omega}y\times (2\mathbb Df)\nu\,d\sigma_y\right)=0,
\end{split}
\end{equation*}
where $f=F|_{\Omega}$,
which together with \eqref{Z-ang} implies that so do the angular velocity for $AF=\mathbb P{\mathscr A}F$.
We thus observe \eqref{til-eq}.

From \eqref{obje-0} we employ \eqref{2nd-led} along with \eqref{til-eq} to deduce
\begin{equation}
\begin{split}
\|\partial_t{\mathcal S}(t)F\|_{q,B_R}
&=\|\partial_te^{-tA}\widetilde F\|_{q,B_R}
=\|e^{-tA}A\widetilde F\|_{q,B_R}
=\|{\mathcal S}(t)A\widetilde F\|_{q,B_R}  \\
&\leq C(1+t)^{-3/2q}\|A\widetilde F\|_{q,\mathbb R^3}
\end{split}
\label{t-diff-led}
\end{equation}
for all $t\geq 0$.
On account of \eqref{ellip-0} near $t=0$ and \eqref{2nd-led} for large $t$, we obtain
\begin{equation}
\|u(t)\|_{W^{1,q}(\Omega_R)}\leq C(1+t)^{-3/2q}\big(\|A\widetilde F\|_{q,\mathbb R^3}+\|\widetilde F\|_{q,\mathbb R^3}\big)
\label{sm-led-0}
\end{equation}
for all $t\geq 0$.
By the choice of the pressure and by using the equation we find
\[
\|p(t)\|_{q,\Omega_R}\leq C\|\nabla p(t)\|_{W^{-1,q}(\Omega_R)}
\leq C\|\partial_tu(t)\|_{q,\Omega_R}+C\|\nabla u(t)\|_{q,\Omega_R},
\]
which combined with \eqref{t-diff-led}--\eqref{sm-led-0} leads us to \eqref{sm-led}. 
The proof is complete.
\end{proof}

Let us proceed to decay estimates of $\nabla^2{\mathcal S}(t)F$.
To this end, 
Lemma \ref{lem-apri} below is useful.
Since the following lemma itself is independent of generation of the fluid-structure semigroup, the assumption on the friction coefficient 
may be allowed to be less.
\begin{lemma}
Suppose that $\alpha\in C(\partial\Omega)$ is nonnegative.
Let $q\in(1,\infty)$ and $R\geq 3$.
Given $U\in D_q(A)$, we set $u=U|_\Omega$.
Then there is a constant $C=C(\alpha,q,R,\Omega)>0$ such that
\begin{equation}
\|u\|_{W^{2,q}(\Omega_R)}\leq C\big(\|AU\|_{q,\Omega_{2R}}+\|\nabla u\|_{q,G_R}+\|U\|_{q,B_{2R}}\big)
\label{ap-est}
\end{equation}
for all $U\in D_q(A)$, where $G_R=B_{2R}\setminus\overline{B_R}$.
\label{lem-apri}
\end{lemma}

\begin{proof}
Given $U\in D_q(A)$, we set $F=AU$.
Let $(u,\eta,\omega)=i(U)$ and $(f,\kappa,\mu)=i(F)$, see \eqref{mono}.
Then there is a pressure $p$ such that
$(u,p,\eta,\omega)$ obeys \eqref{FS-1}--\eqref{FS-2} with $\lambda=0$, where $p$ is singled out in such a way that
\[
\int_{G_R}p\,dx=0,
\]
which together with the equation implies
\begin{equation}
\|p\|_{q,G_R}\leq C\|\nabla p\|_{W^{-1,q}(G_R)}\leq C\big(\|f\|_{q,G_R}+\|\nabla u\|_{q,G_R}\big).
\label{p-apri}
\end{equation}
Let us take a cut-off function $\phi$ satisfying
\begin{equation*}
\phi\in C_0^\infty(B_{2R}), \qquad
0\leq \phi\leq 1, \qquad
\phi=1\;\; \mbox{on $B_{R+1}$},
\end{equation*}
and consider
\[
v=\phi u-\mathbb B[u\cdot\nabla \phi], \qquad p_v=\phi p,
\]
where $\mathbb B=\mathbb B_{G_R}$ denotes the Bogovskii operator. 
We here note that
\begin{equation*}
\int_{G_R}u\cdot\nabla\phi\,dx 
=\int_{\partial B_R}\frac{-x}{R}\cdot u\,d\sigma
=\int_{\partial\Omega}\nu\cdot u\,d\sigma=\int_{\partial\Omega}\nu\cdot (\eta+\omega\times x)\,d\sigma=0 %
\end{equation*}
and that $u\cdot\nabla\phi\in W^{2,q}_0(G_R)$ implies $\mathbb B[u\cdot\nabla\phi]\in W^{3,q}_0(G_R)$ and, thereby,
$\mathbb Dv$ as well as $v$ vanishes at $\partial B_{2R}$.
Then the pair of $v$ and $p_v$ is a solution to
\begin{equation*}
\begin{split}
&-\Delta v+\nabla p_v=g, \quad\mbox{div $v$}=0 \quad \mbox{in $\Omega_{2R}$}, \\
&\nu\cdot (v-\eta-\omega\times x)|_{\partial\Omega}=0, \quad
[(2\mathbb Dv)\nu]_\tau+\alpha(v-\eta-\omega\times x)_\tau|_{\partial\Omega}=0, \\
&\nu\cdot v|_{\partial B_{2R}}=0, \quad [(2\mathbb Dv)\nu]_\tau+\alpha v_\tau|_{\partial B_{2R}}=0,
\end{split}
\end{equation*}
where 
\[
g=\phi f-2\nabla\phi\cdot\nabla u-(\Delta \phi)u
+\Delta\mathbb B[u\cdot\nabla\phi]+(\nabla\phi)p,
\]
which satisfies
\begin{equation}
\|g\|_{q,\Omega_{2R}}\leq C\big(\|f\|_{q,\Omega_{2R}}+\|\nabla u\|_{q,G_R}+\|u\|_{q,G_R}\big)
\label{g-apri}
\end{equation}
by taking into account \eqref{p-apri}.
Using the lifting fuction $\ell(\eta,\omega)$ given by \eqref{lift}, we describe the solution $v$
of the form $v=\ell(\eta,\omega)+w$, then $w$ 
should obey, with the same pressure $p_v$,
\begin{equation*}
\begin{split}
&-\Delta w+\nabla p_v=g+\Delta \ell(\eta,\omega), \quad\mbox{div $w$}=0 \quad \mbox{in $\Omega_{2R}$}, \\
&\nu\cdot w|_{\partial\Omega_{2R}}=0, \quad
[(2\mathbb Dw)\nu]_\tau+\alpha w_\tau|_{\partial\Omega_{2R}}=0.  \\
\end{split}
\end{equation*}

We now rely on the a priori estimate of this boundary value problem in the bounded domain $\Omega_{2R}$, due to
Acevedo, Amrouche, Conca and Ghosh \cite[Theorem 5.10]{A21} under even less assumption on the friction coefficient $\alpha(x)$,
to obtain
\[
\|w\|_{W^{2,q}(\Omega_{2R})}+\|p_v\|_{W^{1,q}(\Omega_{2R})}\leq C\|g+\Delta\ell(\eta,\omega)\|_{q,\Omega_{2R}}.
\]
This combined with \eqref{g-apri} as well as \eqref{mono} and \eqref{lift-est} yield
\begin{equation*}
\begin{split}
\|v\|_{W^{2,q}(\Omega_{2R})}
&\leq C\big(\|f\|_{q,\Omega_{2R}}+\|\nabla u\|_{q,G_R}+\|u\|_{q,G_R}+|\eta|+|\omega|\big)  \\
&\leq C\big(\|AU\|_{q,\Omega_{2R}}+\|\nabla u\|_{q,G_R}+\|U\|_{q,B_{2R}}\big).
\end{split}
\end{equation*}
Since $u=v$ on $\Omega_R$, we have completed the proof of \eqref{ap-est}.
\end{proof}

With \eqref{ap-est} at hand, we readily get the following.
\begin{proposition}
Under the same assumptions as in Corollary \ref{cor-led},
there is a constant $C=C(\alpha,q,R,\Omega)>0$ such that
\begin{equation}
\|u(t)\|_{W^{2,q}(\Omega_R)}+\|p(t)\|_{W^{1,q}(\Omega_R)}
\leq C(1+t)^{-3/2q}\big(\|A\widetilde F\|_{q,\mathbb R^3}+\|\widetilde F\|_{q,\mathbb R^3}\big)
\label{hess-led}
\end{equation}
for all $t\geq 0$ and $F\in D_q(A)$.

Moreover, there is a constant $C=C(\alpha,q,R,\Omega)>0$ such that
\begin{equation}
\|{\mathcal S}(t)F\|_{W^{2,q}(\Omega_R)}\leq Ct^{-3/2q}\|\widetilde F\|_{q,\mathbb R^3}
\label{hess-led2}
\end{equation}
for all $t\geq 1$ and $F\in X_q(\mathbb R^3)$.
If, in particular, the exponent $q$ is specified as below,
then we have
\begin{equation}
\|{\mathcal S}(t)F\|_{\infty, \Omega_R}\leq Ct^{-3/2q}\|\widetilde F\|_{q,\mathbb R^3}\qquad (q>3/2)
\label{sup-led}
\end{equation}
\begin{equation}
\|\nabla {\mathcal S}(t)F\|_{\infty, \Omega_R}\leq Ct^{-3/2q}\|\widetilde F\|_{q,\mathbb R^3}\qquad (q>3)
\label{sup-gr-led}
\end{equation}
for all $t\geq 1$ and $F\in X_q(\mathbb R^3)$.
\label{prop-hess-led}
\end{proposition}

\begin{proof}
By the same reasoning as in \eqref{t-diff-led} we have
\begin{equation*}
\|A{\mathcal S}(t)F\|_{q,B_{2R}}\leq C(1+t)^{-3/2q}\|A\widetilde F\|_{q,\mathbb R^3}
\end{equation*}
for all $t\geq 0$.
This together with \eqref{2nd-led} and \eqref{sm-led} imply \eqref{hess-led} for $u(t)$
by setting $U={\mathcal S}(t)F$ in \eqref{ap-est}.
Then we use the equation and \eqref{sm-led} to find \eqref{hess-led} for $p(t)$.

Since 
\begin{equation}
\widetilde{{\mathcal S}(1)F}={\mathcal S}(1)F=e^{-A}\widetilde F\in D_q(A),
\label{sm-ef0}
\end{equation}
see \eqref{S-rela3},
along with
\begin{equation}
\|A\widetilde{{\mathcal S}(1)F}\|_{q,\mathbb R^3}\leq C\|\widetilde F\|_{q,\mathbb R^3},
\label{sm-ef}
\end{equation}
we see from \eqref{S-semi} and \eqref{hess-led} that
\begin{equation*}
\begin{split}
\|{\mathcal S}(t)F\|_{W^{2,q}(\Omega_R)}
&=\|{\mathcal S}(t-1){\mathcal S}(1)F\|_{W^{2,q}(\Omega_R)}  \\
&\leq Ct^{-3/2q}\big(\|A\widetilde{{\mathcal S}(1)F}\|_{q,\mathbb R^3}+\|\widetilde{{\mathcal S}(1)F}\|_{q,\mathbb R^3}\big)  \\
&\leq Ct^{-3/2q}\|\widetilde F\|_{q,\mathbb R^3}
\end{split}
\end{equation*}
for all $t\geq 1$ and $F\in X_q(\mathbb R^3)$.
Finally, it is obvious that \eqref{hess-led2} leads to \eqref{sup-led}--\eqref{sup-gr-led} by the embedding relation.
The proof is complete.
\end{proof}

\section{Proof of Theorem \ref{main}}
\label{pr-main}

To complete the proof of Theorem \ref{main}, the final step is to show the following estimate near spatial infinity.
\begin{proposition}
Suppose that $\alpha\in C^1(\partial\Omega)$ is nonnegative.
Let $q\in (1,\infty)$ and $R\geq 3$.
Then there is a constant $C=C(\alpha,q,R,\Omega)>0$ such that
\begin{equation}
\|{\mathcal S}(t)F\|_{q,\mathbb R^3\setminus B_R}\leq C\|\widetilde F\|_{q,\mathbb R^3}
\label{est-far0}
\end{equation}
\begin{equation}
\|\nabla^j{\mathcal S}(t)F\|_{q,\mathbb R^3\setminus B_R}
\leq Ct^{-\min\{j/2,3/2q\}}\|\widetilde F\|_{q,\mathbb R^3}
\label{est-far}
\end{equation}
for all $t>3$, $F\in X_q(\mathbb R^3)$ and $j\in \{1,2\}$, 
where $\widetilde F$ is given by \eqref{IC-modi}--\eqref{data-modi}. 

If, in particular, $q\in (3/2,\infty)$, 
then we have
\begin{equation}
\|{\mathcal S}(t)F\|_{W^{1,\infty}(\mathbb R^3\setminus B_R)}\leq Ct^{-3/2q}\|\widetilde F\|_{q,\mathbb R^3} 
\label{sup-far}
\end{equation}
for all $t>3$ and $F\in X_q(\mathbb R^3)$.
\label{prop-far}
\end{proposition}

\begin{proof}
Given $F\in X_q(\mathbb R^3)$, we set $U(t)={\mathcal S}(t)G={\mathcal S}(t+1)F$ with
$G={\mathcal S}(1)F=e^{-A}\widetilde F\in D_q(A)$ and $u(t)=U(t)|_\Omega$.
Then, for the proof of \eqref{est-far0}--\eqref{sup-far}, it suffices to deduce
\begin{equation}
\|\nabla^ju(t)\|_{q,\mathbb R^3\setminus B_R}\leq Ct^{-\min\{j/2,3/2q\}}
\big(\|A\widetilde G\|_{q,\mathbb R^3}+\|\widetilde G\|_{q,\mathbb R^3}\big)
\label{est-far2}
\end{equation}
\begin{equation}
\|u(t)\|_{W^{1,\infty}(\mathbb R^3\setminus B_R)}
\leq Ct^{-3/2q}
\big(\|A\widetilde G\|_{q,\mathbb R^3}+\|\widetilde G\|_{q,\mathbb R^3}\big)
\label{sup-far2}
\end{equation}
for all $t>2$ and $j\in \{0,1,2\}$, which together with \eqref{sm-ef0}--\eqref{sm-ef} imply the desired estimates 
for $t>3$.

To this end, with use of the same cut-off function $\phi$
as in \eqref{cut}
and the Bogovskii operator $\mathbb B=\mathbb B_D$ in the domain $D=B_3\setminus \overline{B_1}$,
we consider
\[
v(t)=(1-\phi)u(t)+\mathbb B[u(t)\cdot\nabla\phi], \qquad p_v(t)=(1-\phi)p(t),
\]
where the pressure $p(t)$ is singled out as in Corollary \ref{cor-led}, that is,
\[
\int_{\Omega_R}p(t)\,dx=0.
\]
Then $v(t)$ obeys the equation
\[
v(t)=e^{t\Delta}G^\prime+\int_0^t e^{(t-s)\Delta}P_{\mathbb R^3}H(s)\,ds,
\]
where $P_{\mathbb R^3}=I+{\mathcal R}\otimes {\mathcal R}$ denotes the Fujita-Kato projection on the whole space $\mathbb R^3$ with 
${\mathcal R}=(-\Delta)^{-1/2}\nabla$ being the Riesz transform,
$G^\prime=(1-\phi)G+\mathbb B[G\cdot\nabla\phi]$ and
\begin{equation}
H(t)=2\nabla\phi\cdot\nabla u(t)+(\Delta\phi)u(t)-\Delta\mathbb B[u(t)\cdot\nabla\phi]
+\mathbb B[\partial_tu(t)\cdot\nabla\phi]-(\nabla\phi)p(t).
\label{far-gar}
\end{equation}
Our task is to estimate $v(t)$ 
since $v(t)=u(t)$ outside $B_R$ for $R\geq 3$.
As for the term $e^{t\Delta}G^\prime$, 
we have
\[
\|\nabla^je^{t\Delta}G^\prime\|_{q,\mathbb R^3}\leq Ct^{-j/2}\|G^\prime\|_{q,\mathbb R^3}
\leq Ct^{-j/2}\|G\|_{q,\mathbb R^3}=Ct^{-j/2}\|\widetilde G\|_{q,\mathbb R^3}
\]
and
\[
\|\nabla^m e^{t\Delta}G^\prime\|_{\infty,\mathbb R^3}
\leq Ct^{-m/2-3/2q}\|\widetilde G\|_{q,\mathbb R^3}
\]
for $t>0$, $j\in \{0,1,2\}$ and $m\in \{0,1\}$. 

Let $r\in (1,q]$, that plays a role as a parameter.
Since $P_{\mathbb R^3}$ is bounded on $W^{1,r}(\mathbb R^3)$, 
we employ \eqref{bog-est1}, \eqref{sm-led} and \eqref{hess-led} to see that the remaining term $H(t)$ given by \eqref{far-gar} satisfies
\begin{equation}
\begin{split}
\|P_{\mathbb R^3}H(t)\|_{W^{1,r}(\mathbb R^3)}
&\leq C\|H(t)\|_{W^{1,r}(\mathbb R^3)}  \\
&\leq C\big(\|u(t)\|_{W^{2,q}(\Omega_R)}+\|\partial_tu(t)\|_{q,\Omega_R}+\|p(t)\|_{W^{1,q}(\Omega_R)}\big)  \\
&\leq C(1+t)^{-3/2q}\big(\|A\widetilde G\|_{q,\mathbb R^3}+\|\widetilde G\|_{q,\mathbb R^3}\big),
\end{split}
\label{est-H}
\end{equation}
from which it follows that
\begin{equation}
\begin{split}
&\quad \left\|\partial_k^{j-1}\partial_\ell\int_0^t e^{(t-s)\Delta} P_{\mathbb R^3}H(s)\,ds\right\|_{q,\mathbb R^3}  \\
&\leq \int_0^{t-1}
\|\nabla^je^{(t-s)\Delta}P_{\mathbb R^3}H(s)\|_{q,\mathbb R^3}\,ds
+\int_{t-1}^t\|\partial_\ell e^{(t-s)\Delta}\partial_k^{j-1}P_{\mathbb R^3}H(s)\|_{q,\mathbb R^3}\,ds  \\
&\leq C
\Big(\int_0^{t-1}(t-s)^{-j/2}(1+t-s)^{-\frac{3}{2}(\frac{1}{r}-\frac{1}{q})}(1+s)^{-3/2q}\,ds  \\
&\qquad\quad +\int_{t-1}^t(t-s)^{-1/2}(1+s)^{-3/2q}\,ds\Big)
\big(\|A\widetilde G\|_{q,\mathbb R^3}+\|\widetilde G\|_{q,\mathbb R^3}\big)
\end{split}
\label{far-split}
\end{equation}
for $t>2$, $j\in \{1,2\}$ and $k,\,\ell\in \{1,2,3\}$,
where we have used 
\[
\partial_k(e^{t\Delta}h)(x)=-\int_{\mathbb R^3}\partial_{y_k}\Gamma(x-y,t)h(y)\,dy
=\int_{\mathbb R^3}\Gamma(x-y,t)(\partial_kh)(y)\,dy
\]
with $h=P_{\mathbb R^3}H$ for the latter integral above.
We split the former integral of \eqref{far-split} into $\int_0^{t/2}+\int_{t/2}^{t-1}$ and
take a suitable $r\in (1,q]$ (for given $q$)
to conclude that the right-hand side of \eqref{far-split} is bounded from above by
$Ct^{-\min\{j/2,3/2q\}}$, which implies
\eqref{est-far2} with $j\in \{1,2\}$.
The proof for the case $j=0$ is similar and thus may be omitted. 

Let us next consider \eqref{sup-far2}.
The same manner with use of \eqref{est-H} leads us to
\begin{equation*}
\begin{split}
&\quad \left\|\partial_k^j\int_0^t e^{(t-s)\Delta} P_{\mathbb R^3}H(s)\,ds\right\|_{\infty,\mathbb R^3}  \\
&\leq 
\int_0^t\|e^{(t-s)\Delta}\partial_k^{j}P_{\mathbb R^3}H(s)\|_{\infty,\mathbb R^3}\,ds  \\
&\leq C
\int_0^{t}(t-s)^{-3/2q}(1+t-s)^{-\frac{3}{2}(\frac{1}{r}-\frac{1}{q})}(1+s)^{-3/2q}\,ds\,
\big(\|A\widetilde G\|_{q,\mathbb R^3}+\|\widetilde G\|_{q,\mathbb R^3}\big)
\end{split}
\end{equation*}
for $t>2$ and $j\in \{0,1\}$ as long as $q>3/2$.
The same splitting of the integral as before with $r\in (1,q]$ chosen appropriately implies \eqref{sup-far2}.
We have thus completed the proof. 
\end{proof}
\noindent
{\it Proof of Theorem \ref{main}}.
With estimates obtained so far at hand, it is not hard to complete the proof of \eqref{decay-1}--\eqref{decay-tempo} and 
\eqref{decay-2}--\eqref{decay-4}.
Among them, \eqref{decay-body} has been already claimed in Corollary \ref{cor-rigid}.

{\it 1}.
Let us fix $R\geq 3$.
We collect \eqref{reg-S0},
\eqref{2nd-led}, \eqref{hess-led2} and \eqref{est-far0}--\eqref{est-far} to see that
\begin{equation}
\begin{split}
\|{\mathcal S}(t)F\|_{q,\mathbb R^3}
&\leq C\|\widetilde F\|_{q,\mathbb R^3}, \\
\|\nabla^j{\mathcal S}(t)F\|_{q,\Omega}
&\leq C
t^{-j/2}(1+t)^{\max\{\frac{j}{2}-\frac{3}{2q},0\}}
\|\widetilde F\|_{q,\mathbb R^3},
\end{split}
\label{decay-deri}
\end{equation}
for all $t>0$, $F\in X_q(\mathbb R^3)$ and $j\in \{1,2\}$.
We will show that
\begin{equation}
\|{\mathcal S}(t)F\|_{r,\mathbb R^3}\leq Ct^{-\frac{3}{2}(\frac{1}{q}-\frac{1}{r})}\|\widetilde F\|_{q,\mathbb R^3}
\label{S-decay2-noch}
\end{equation}
for all $t>0$ and $F\in X_q(\mathbb R^3)$ whenever $q\in (1,\infty)$ and $r\in [q,\infty]$.

We see from the interpolation inequality that
$\mbox{\eqref{decay-deri}}_2$ with $j=1$ implies
\begin{equation*}
\|{\mathcal S}(t)F\|_{r,\Omega}\leq Ct^{-\frac{3}{2}(\frac{1}{q}-\frac{1}{r})}\|\widetilde F\|_{q,\mathbb R^3}
\end{equation*}
provided
\begin{equation}
1<q<3, \qquad 0\leq \frac{1}{q}-\frac{1}{r}\leq \frac{1}{3}
\label{restr-1}
\end{equation}
which along with \eqref{decay-body} leads us to \eqref{S-decay2-noch}
for such $(q,r)$. 
On the other hand, we gather \eqref{reg-S},
\eqref{sup-led} and \eqref{sup-far} to find
\begin{equation*}
\|{\mathcal S}(t)F\|_{\infty,\Omega}\leq Ct^{-3/2q}
\|\widetilde F\|_{q,\mathbb R^3}
\end{equation*}
for all $t>0$ when $q\in (3/2,\infty)$.
By \eqref{decay-body} again we also have
\begin{equation}
\|{\mathcal S}(t)F\|_{\infty,\mathbb R^3}\leq Ct^{-3/2q}
\|\widetilde F\|_{q,\mathbb R^3}
\label{sup-decay}
\end{equation}
for all $t>0$
when $q\in (3/2,\infty)$.
This together with \eqref{S-decay2-noch} for $(q,r)$ satisfying 
\eqref{restr-1} leads to \eqref{sup-decay} for every $q\in (1,\infty)$
by use of the semigroup property \eqref{S-semi} and by taking into account \eqref{S-rela3}.
Successively, this combined with 
$\mbox{\eqref{decay-deri}}_1$ implies
\eqref{S-decay2-noch} for all $q\in (1,\infty)$ and $r\in [q,\infty]$, yielding \eqref{decay-1}.

Let $k\geq 1$ be an integer, then we have
\[
\partial_t^k e^{-tA}F=\partial_t^k {\mathcal S}(t)F=\partial_t^ke^{-tA}\widetilde F.
\]
We utilize \eqref{bdd-anal} to get
\[
\|\partial_t^k e^{-tA}F\|_{r,\mathbb R^3}\leq Ct^{-k}\|e^{-\frac{t}{2}A}\widetilde F\|_{r,\mathbb R^3}
=Ct^{-k}\|{\mathcal S}(t/2)F\|_{r,\mathbb R^3}
\]
for all $t>0$ and $F\in X_q(\mathbb R^3)$.
Thus, \eqref{decay-tempo} follows from \eqref{S-decay2-noch}.

{\it 2}.
With the aid of \eqref{S-semi} and \eqref{S-rela3} again,
we combine
\eqref{decay-deri} with \eqref{S-decay2-noch} to obtain \eqref{decay-2} except for $\|\nabla{\mathcal S}(t)F\|_{\infty,\Omega}$.
As to this latter one, we know from $\mbox{\eqref{reg-S}}_2$,
\eqref{sup-gr-led} and \eqref{sup-far} that
\begin{equation*}
\|\nabla {\mathcal S}(t)F\|_{\infty,\Omega}\leq Ct^{-1/2-3/2q}(1+t)^{1/2}\|\widetilde F\|_{q,\mathbb R^3}
\end{equation*}
for all $t>0$ and $F\in X_q(\mathbb R^3)$ when $q\in (3,\infty)$.
However, this restriction is removed by using  \eqref{S-decay2-noch} together with the splitting ${\mathcal S}(t)={\mathcal S}(\frac{t}{2}){\mathcal S}(\frac{t}{2})$.
This completes the proof of \eqref{decay-2}.

{\it 3}.
From the assumptions on $G$, \eqref{S-rela} and \eqref{proj-ang2} together with \eqref{data-modi}, we observe
\begin{equation*}
e^{-tA}\mathbb P\mbox{div $G$}={\mathcal S}(t)\mathbb P\mbox{div $G$}.
\end{equation*}
Then we find from \eqref{sym-P} and \eqref{S-sym} that
\begin{equation*}
\begin{split}
&\quad |\langle e^{-tA}\mathbb P\mbox{div $G$},\phi\rangle_{\mathbb R^3,\rho}|
=|\langle \mbox{div $G$}, {\mathcal S}(t)\phi\rangle_{\mathbb R^3,\rho}|  \\
&=\left|
(1-\rho)\int_{\partial\Omega}(G\nu)\cdot {\mathcal S}(t)\phi\,dx
-\rho\int_{\mathscr B}G:\nabla {\mathcal S}(t)\phi\,dx
-\int_\Omega G:\nabla {\mathcal S}(t)\phi\,dx
\right|  \\
&\leq\rho\|G\|_{q,{\mathscr B}}\|\nabla {\mathcal S}(t)\phi\|_{q^\prime,{\mathscr B}}
+\|G\|_{q,\Omega}\|\nabla {\mathcal S}(t)\phi\|_{q^\prime,\Omega}
\end{split}
\end{equation*}
for all $\phi\in {\mathcal E}(\mathbb R^3)$
on account of $(G\nu)|_{\partial\Omega}=0$.
In view of \eqref{S-body}, we have
\[
\|\nabla {\mathcal S}(t)\phi\|_{q^\prime,{\mathscr B}}
\leq C\left\{
\begin{array}{ll}
0&\mbox{\rm case (i)},  \\
|\gamma\times\omega_\phi(t)|\qquad&\mbox{\rm case (ii)}, \\
|\omega_\phi(t)|&\mbox{\rm case (iii)},
\end{array}
\right.
\]
where $\omega_\phi(t)=J^{-1}\int_{\mathscr B}x\times ({\mathcal S}(t)\phi)(x)\rho\,dx$.
By duality, the inequality above
together with \eqref{decay-body} and
\eqref{decay-2} concludes \eqref{decay-3}
except the case $r=\infty$, that is covered by the semigroup property and \eqref{S-decay2-noch} with $r=\infty$.

{\it 4}.
Finally, we combine \eqref{decay-2} with \eqref{decay-3} to furnish \eqref{decay-4}.
The proof is complete.
\hfill
$\Box$

\section{Alternative approach to Proposition \ref{prop-led}} 
\label{alter}

Spectral analysis done in section \ref{from-to} is known as a strong tool to get a sharp large time decay of the semigroup.
Indeed, such analysis was quite successful in the literature \cite{I89,KS98} for the exterior problem, however, one needs
some other approaches when the linearized system is non-autonomous, for which the proof of even uniformly boundedness
is mostly involved. 
In \cite{H18,H20} the first author 
has provided an alternative route to develop $L^q$-$L^r$ decay estimates of a genelarized Oseen
evolution operator arising from time-dependent rigid motions in exterior domains
and, successively in \cite{H24}, this strategy has been well adapted to the Oseen-structure evolution operator. 
Those papers cover the autonomous regime as a particular case.
In this section, we make it clear that the similar approach is indeed possible for our problem. 
Since the fluid-structure semigroup under consideration is already known to be uniformly bounded, see \eqref{bdd-anal} with $k=0$, 
which implies that so is the solution operator ${\mathcal S}(t)$ given by \eqref{obje-0}, that is,
\begin{equation}
\|{\mathcal S}(t)F\|_{q,\mathbb R^3}
=\|e^{-tA}\widetilde F\|_{q,\mathbb R^3}
\leq C\|\widetilde F\|_{q,\mathbb R^3},
\label{unif-bdd-S}
\end{equation} 
the proof could be shortened.
On the other hand, there is a remarkable difference between the Navier and no-slip conditions;
in fact, $|\eta|$ and $|\omega|$ can be easily controlled by $\|\mathbb Du\|_{2,\Omega}$ under the latter condition, whereas
this issue is involved for the case of the Navier condition.
We will explain in the next why it is important and how one can overcome the difficulty.

We focus on \cite[Lemma 4.1]{H24}, which claims that, under the no-slip boundary condition,
the uniformly boundedness implies the $L^q$-$L^r$ estimates, in which one of keys of the proof is the energy relation.
Under the Navier boundary condition, the same argument as in deduction of \eqref{energy} implies that 
\begin{equation}
\frac{1}{2}\frac{d}{dt}\|U(t)\|_{X_2(\mathbb R^3)}^2+2\|\mathbb Du(t)\|_{2,\Omega}^2
+\int_{\partial\Omega}\alpha(x)|u-\eta-\omega\times x|^2\,d\sigma=0
\label{energy-evo}
\end{equation}
for $U(t)=e^{-tA}F$ with $F\in X_2(\mathbb R^3)$,
where $\|\cdot\|_{X_2(\mathbb R^3)}$ is the energy norm \eqref{ener} and
$(u,\eta,\omega)=i(U(t))$.
This relation is reduced
to \eqref{energy-evo} 
with the first term replaced by 
$\frac{1}{2}\frac{d}{dt}\|{\mathcal S}(t)F\|_{X_2(\mathbb R^3)}^2$, where
\begin{equation*}
\|{\mathcal S}(t)F\|_{X_2(\mathbb R^3)}^2=
\left\{
\begin{array}{ll}
\|u\|_{2,\Omega}^2+m|\eta|^2&\mbox{case (i)}, \\
\|u\|_{2,\Omega}^2+m|\eta|^2+\big[J\{(\gamma\times\omega)\times\gamma\}\big] \cdot\{(\gamma\times\omega)\times\gamma\}\qquad
 &\mbox{case (ii)}, \\
\|u\|_{2,\Omega}^2+m|\eta|^2+(J\omega)\cdot\omega&\mbox{case (iii)},
\end{array}
\right.
\end{equation*}
because of \eqref{sym-inertia}, \eqref{const-ii} and \eqref{const-i}.
If, in particular, the friction $\alpha\in C^1(\partial\Omega)$ is not identically zero for the case (iii),
there are a nonempty open subset $\Gamma$ of $\partial\Omega$ (with respect to the induced topology)
and a constant $\alpha_0>0$ such that
\[
\alpha(x)\geq \alpha_0\qquad\forall\,x\in\Gamma.
\]
Then we have
\begin{equation}
\frac{1}{2}\frac{d}{dt}\|{\mathcal S}(t)F\|_{X_2(\mathbb R^3)}^2+2\|\mathbb Du(t)\|_{2,\Omega}^2
+\alpha_0\int_{\Gamma}|u-\eta-\omega\times x|^2\,d\sigma\leq 0.
\label{energy-evo2}
\end{equation}

Let $q\in (1,2)$, then the interpolation inequality along with the uniformly boundedness \eqref{unif-bdd-S} implies that 
\begin{equation}
\|{\mathcal S}(t)F\|_{X_2(\mathbb R^3)}\leq C\|\widetilde F\|^{1-\theta}_{q,\mathbb R^3}\|{\mathcal S}(t)F\|_{6,\mathbb R^3}^\theta
\label{inter}
\end{equation}
where $\theta\in (0,1)$ satisfies $1/2=\theta/6+(1-\theta)/q$.
As for $U(t)=e^{-tA_{\rm no-slip}}F$ being the fluid-structure semigroup 
under the no-slip condition (\cite{H24}), the term $\|U(t)\|_{6,\mathbb R^3}$,
in which the rigid motion $\eta+\omega\times x$ is involved outside $\Omega$, can be controlled by
$\|\nabla U(t)\|_{2,\mathbb R^3}=\sqrt{2}\|\mathbb Du(t)\|_{2,\Omega}$; however, this no longer holds 
(even $\nabla U(t)\notin L^2(\mathbb R^3)$)
under the Navier boundary condition and this is why an appropriate estimate of $(\eta,\omega)$ must be an issue here.

To overcome the difficulty, we need two ingredients.
One is Lemma \ref{lem-rigid-0} below, the other is the Korn inequality in exterior domains
\begin{equation}
\|\nabla u\|_{2,\Omega}^2\leq C\|\mathbb Du\|_{2,\Omega}^2
\label{korn}
\end{equation}
for all $u\in \widehat W^{1,2}(\Omega)^3$, which is defined as the completion of $C_0^\infty(\overline{\Omega})^3$
with respect to the norm $\|\nabla(\cdot)\|_{2,\Omega}$, see
Shibata and Soga \cite[Theorem 1.5]{SSo89} and Ito \cite[Theorem 6.5]{It90}.
What is remarkable is that \eqref{korn} holds true without any boundary condition at $\partial\Omega$ nor solenoidal condition.
\begin{lemma}
Let ${\mathscr B}$ 
be a bounded domain as in \eqref{dom}. 
\begin{enumerate}
\item
There is a constant $C=C({\mathscr B})>0$ such that
\begin{equation}
\begin{split}
&\mbox{(I)}\qquad
|\eta|^2\leq C\int_{\partial{\mathscr B}}|\nu\cdot\eta|^2\,d\sigma,  \\
&\mbox{(II)}\qquad
|\eta|^2+|\gamma\times \omega|^2
\leq C\int_{\partial{\mathscr B}}
\Big|\nu\cdot\big[\eta+\{(\gamma\times\omega)\times\gamma\}\times x\big]\Big|^2\,d\sigma,  \\
&\mbox{(III)}\qquad
|\eta|^2+|\omega|^2\leq C\int_{\partial{\mathscr B}}|\nu\cdot (\eta+\omega\times x)|^2\,d\sigma,
\end{split}
\label{rigid-est-0}
\end{equation}
for all $(\eta,\omega)\in \mathbb R^3\times \mathbb R^3$ provided, respectively, that

(I) ${\mathscr B}$ is a ball;

(II) ${\mathscr B}$ is axisymmetric about a single axis $\gamma\mathbb R$ for some $\gamma\in\mathbb R^3$, $|\gamma|=1$;

(III) ${\mathscr B}$ is not axially symmetric. 

Note that 
$\nu\cdot\big[\eta+\{(\gamma\times\omega)\times\gamma\}\times x\big]=\nu\cdot(\eta+\omega\times x)$
in the right-hand side of \eqref{rigid-est-0}--(II) 
on account of \eqref{sym-1}.

\item
Let $\Gamma$ be a nonempty open subset of $\partial {\mathscr B}$ with respect to the induced topology.
Then there is a constant $C=C(\Gamma,{\mathscr B})>0$ such that
\begin{equation}
|\eta|^2+|\omega|^2\leq C\int_\Gamma |\eta+\omega\times x|^2\,d\sigma
\label{rigid-est}
\end{equation}
for all $(\eta,\omega)\in\mathbb R^3\times\mathbb R^3$, regardless of the shape of ${\mathscr B}$. 
\end{enumerate}
\label{lem-rigid-0}
\end{lemma}

\begin{proof}
Let us focus on \eqref{rigid-est-0}--{\it (II)}.
Since $|\gamma\times \omega|=|(\gamma\times \omega)\times\gamma|$, we may show
\[
|\eta|^2+|(\gamma\times \omega)\times\gamma|^2\leq C\int_{\partial{\mathscr B}}
\Big|\nu\cdot\big[\eta+\{(\gamma\times\omega)\times\gamma\}\times x\big]\Big|^2\,d\sigma,
\]
that suits us for the contradiction argument.
Suppose the contrary, then there is a sequence $\{(\eta_j,\omega_j)\}\subset \mathbb R^3\times\mathbb R^3$ such that
\[
|\eta_j|^2+|(\gamma\times \omega_j)\times \gamma|^2=1, \qquad
\lim_{j\to\infty}\int_{\partial{\mathscr B}}
\Big|\nu\cdot\big[\eta_j+\{(\gamma\times\omega_j)\times\gamma\}\times x\big]\Big|^2\,d\sigma=0.
\]
Along a subsequence, that is denoted by the same symbol for simplicity,
there is $(\eta_\infty,\zeta_\infty)\in \mathbb R^3\times \mathbb R^3$ satisfying
\[
\lim_{j\to\infty}\eta_j=\eta_\infty, \qquad \lim_{j\to\infty}(\gamma\times\omega_j)\times \gamma=\zeta_\infty.
\]
Then the limt should obey
\[
|\eta_\infty|^2+|\zeta_\infty|^2=1, \qquad
\gamma\cdot\zeta_\infty=0, \qquad
\int_{\partial{\mathscr B}} |\nu\cdot(\eta_\infty+\zeta_\infty\times x)|^2\,d\sigma=0,
\]
the latter of which implies that
$\nu\cdot (\eta_\infty+\zeta_\infty\times x)|_{\partial{\mathscr B}}=0$. 
By Lemma \ref{lem-rigid-tan} we are led to $(\eta_\infty,\zeta_\infty)=(0,0)$, yielding a contradiction.
The other cases {\it (I)} and {\it (III)} are discussed by use of Lemma \ref{lem-rigid-tan} in the same fashion.

One can also show \eqref{rigid-est} by the similar contradiction argument.
In fact, suppose the contrary, then there is a sequence $\{(\eta_j,\omega_j)\}$ such that
\[
|\eta_j|^2+|\omega_j|^2=1, \qquad
\lim_{j\to\infty}\int_\Gamma |\eta_j+\omega_j\times x|^2\,d\sigma=0.
\]
The limit $(\eta_\infty,\omega_\infty)$ along a subsequence fulfills
\[
|\eta_\infty|^2+|\omega_\infty|^2=1, \qquad
\int_\Gamma |\eta_\infty+\omega_\infty\times x|^2\,d\sigma=0,
\]
and, thereby, $\eta_\infty+\omega_\infty\times x=0$ for every $x\in\Gamma$.
As observed in the proof of Proposition \ref{uni}
for the case when $\alpha$ is not identically zero as well as $\lambda=0$, we infer $(\eta_\infty,\omega_\infty)=(0,0)$.
The proof is complete.
\end{proof}
\begin{remark}
The cases (I)--(III) in Lemma \ref{lem-rigid-0} is partly different from the ones \emph{(i)--(iii)} indicated in Theorem \ref{main}.
For instance, the ball case when $\alpha(x)$ is not identically zero is included in \emph{(iii)}, however, one can control merely $\eta$ by
\eqref{rigid-est-0}--(I).
\label{rem-classification}
\end{remark}

Combining Lemma \ref{lem-rigid-0} with the Korn inequality \eqref{korn} allows us to control $(\eta,\omega)$ by $\mathbb Du$.
\begin{lemma}
Let ${\mathscr B}$ and $\Omega$ be domains as in \eqref{dom}.
Let $U\in X_2(\mathbb R^3)$ and set $(u,\eta,\omega)=i(U)$, see \eqref{mono}.
We assume that $u\in W^{1,2}(\Omega)$.

\begin{enumerate}
\item
Consider three cases (I)--(III) indicated in Lemma \ref{lem-rigid-0}.
Then there is a constant $C=C(\Omega)>0$ such that
\begin{equation}
\left.
\begin{array}{lr}
\mbox{(I)}&|\eta|^2 \\
\mbox{(II)}\qquad\qquad&|\eta|^2+|\gamma\times \omega|^2 
   \\
\mbox{(III)}&|\eta|^2+|\omega|^2
\end{array}
\right\}
\leq C\|\mathbb Du\|_{2,\Omega}^2. \qquad\qquad
\label{rigid-control-0}
\end{equation}

\item
Let $\Gamma$ be a nonempty open subset of $\partial{\mathscr B}$ with respect to the induced topology.
Then there is a constant $C=C(\Gamma,\Omega)>0$ such that
\begin{equation}
|\eta|^2+|\omega|^2\leq C
\left(
\int_{\Gamma}|u-\eta-\omega\times x|^2\,d\sigma+\|\mathbb Du\|_{2,\Omega}^2
\right),
\label{rigid-control}
\end{equation}
regardless of the shape of ${\mathscr B}$. 
\end{enumerate}
\label{lem-rigid}
\end{lemma}

\begin{proof}
We collect the trace estimate, embedding relation and \eqref{korn} to furnish
\begin{equation}
\|u\|_{2,\Gamma}^2\leq \|u\|_{2,\partial\Omega}^2
\leq C\|u\|_{W^{1,2}(\Omega_3)}^2
\leq C\big(\|u\|_{6,\Omega_3}^2+\|\nabla u\|_{2,\Omega_3}^2\big)
\leq C\|\nabla u\|_{2,\Omega}^2
\leq C\|\mathbb Du\|_{2,\Omega}^2.
\label{tr-korn}
\end{equation}
This along with \eqref{rigid-est} and simply the triangle inequality at once yield \eqref{rigid-control},
where $u$ can be independent of $(\eta,\omega)$ unlike \eqref{rigid-control-0}.

Recall that $U\in X_2(\mathbb R^3)$ implies \eqref{normal-va} and, particularly,
\begin{equation*}
\begin{array}{ll}
\mbox{\it (I)}&\nu\cdot(u-\eta)=0, \\
\mbox{\it (II)}\qquad&\nu\cdot\big[u-\eta-\{(\gamma\times \omega)\times\gamma\}\times x\big]=0,
\end{array} \qquad\qquad
\end{equation*}
at $\partial\Omega$ on account of specific shape for those cases, respectively, see \eqref{sym-1} for the case {\it (II)}. 
Due to this boundary condition together with \eqref{rigid-est-0}, the left-hand side of \eqref{rigid-control-0} for every case
is bounded from above by $C\int_{\partial\Omega}|\nu\cdot u|^2\,d\sigma$, 
which concludes \eqref{rigid-control-0} because of \eqref{tr-korn}.
\end{proof}

Let us go back to \eqref{inter} for $F\in {\mathcal E}(\mathbb R^3)$, see \eqref{nice-cl}, and consider 
$\|{\mathcal S}(t)F\|_{6,\mathbb R^3}$.
We make use of
\begin{equation*}
\begin{split}
&\mbox{\eqref{rigid-control-0}--{\it (I)} for the case (i)}; \\
&\mbox{\eqref{rigid-control-0}--{\it (II)} for the case (ii)}; \\
&\mbox{\eqref{rigid-control-0}--{\it (III)} for the case (iii-1): ${\mathscr B}$ is not axially symmetric}; \\ 
&\mbox{\eqref{rigid-control} for the case (iii-2): $\alpha\in C^1(\partial\Omega)$ is not identically zero},
\end{split}
\end{equation*}
where three cases (i)--(iii) are indicated in Theorem \ref{main}.
The case (iii) is covered by the latter two listed above.
In view of \eqref{obje-0} and by \eqref{rigid-control-0}--\eqref{rigid-control} together with \eqref{korn}, we get
\begin{equation}
\begin{split}
\|{\mathcal S}(t)F\|_{6,\mathbb R^3}
&\leq \|u(t)\|_{6,\Omega}+C
\left\{
\begin{array}{ll}
|\eta(t)|&\mbox{case (i)}, \\
|\eta(t)|+|\gamma\times\omega(t)| 
\qquad&\mbox{case (ii)}, \\
|\eta(t)|+|\omega(t)|&\mbox{case (iii)},
\end{array}
\right.   \\
&\leq C 
\left\{
\begin{array}{ll}
\|\mathbb Du(t)\|_{2,\Omega}&\mbox{case (i), (ii), (iii-1)}, \\
\|u(t)-\eta(t)-\omega(t)\times x\|_{2,\Gamma}+\|\mathbb Du(t)\|_{2,\Omega}\qquad&\mbox{case (iii-2)}.
\end{array}
\right.
\end{split}
\label{L6}
\end{equation}

We combine \eqref{inter}, 
in which one may replace $\|\widetilde F\|_{q,\mathbb R^3}$ by $\|F\|_{q,\mathbb R^3}$, with \eqref{L6} to find
\[
\frac{\|{\mathcal S}(t)F\|_{X_2(\mathbb R^3)}^{2/\theta}}{\|F\|_{q,\mathbb R^3}^{2(1/\theta-1)}}
\leq C
\left\{
\begin{array}{ll}
\|\mathbb Du(t)\|_{2,\Omega}^2&\mbox{case (i), (ii), (iii-1)}, \\
\|u(t)-\eta(t)-\omega(t)\times x\|_{2,\Gamma}^2+\|\mathbb Du(t)\|_{2,\Omega}^2\qquad&\mbox{case (iii-2)},
\end{array}
\right.
\]
and, successively,
\begin{equation}
\frac{d}{dt}\|{\mathcal S}(t)F\|_{X_2(\mathbb R^3)}^2+
\frac{C\|{\mathcal S}(t)F\|_{X_2(\mathbb R^3)}^{2/\theta}}{\|F\|_{q,\mathbb R^3}^{2(1/\theta-1)}}\leq 0
\label{diff-ineq}
\end{equation}
for all $t>0$ and 
$F\in{\mathcal E}(\mathbb R^3)\setminus \{0\}$ by using the energy relations \eqref{energy-evo} with $\alpha\equiv 0$ and
\eqref{energy-evo2}.
With \eqref{diff-ineq} for every $q\in (1,2)$ at hand, we are led to 
\begin{equation}
\|{\mathcal S}(t)F\|_{r,\mathbb R^3}\leq Ct^{-\frac{3}{2}(\frac{1}{q}-\frac{1}{r})}\|F\|_{q,\mathbb R^3}
\label{S-decay}
\end{equation}
for all $t>0$, $F\in X_q(\mathbb R^3)$ and every $(q,r)$ with $1< q\leq r<\infty$ in the following way:
Indeed, we solve the differential inequality \eqref{diff-ineq} to
show \eqref{S-decay} with $q\in (1,2)$ and $r=2$, where $F$ can be taken from $X_q(\mathbb R^3)$ since
${\mathcal E}(\mathbb R^3)$ is dense there.
This together with the uniformly boundedness \eqref{unif-bdd-S}
implies \eqref{S-decay} for $1<q\leq r\leq 2$. 
By duality using \eqref{S-sym} and by the semigroup property \eqref{S-semi} we conclude \eqref{S-decay} for $1<q\leq r<\infty$. 
Since we know ${\mathcal S}(t)F={\mathcal S}(t)\widetilde F$, see \eqref{S-rela3}, we rewrite \eqref{S-decay} as
\begin{equation}
\|{\mathcal S}(t)F\|_{r,\mathbb R^3}\leq Ct^{-\frac{3}{2}(\frac{1}{q}-\frac{1}{r})}\|\widetilde F\|_{q,\mathbb R^3}
\label{S-decay2}
\end{equation}
for all $t>0$, $F\in X_q(\mathbb R^3)$ and $1<q\leq r<\infty$.
We thus obtain \eqref{decay-1} for such $(q,r)$.

Finally, let us discuss the local energy decay property for compactly supported data.
Let $F\in L^q_{{\rm R},[d]}(\mathbb R^3)$ and set $V(t)={\mathcal S}(t)\mathbb PF={\mathcal S}(t)\mathbb P\widetilde F$,
see \eqref{obje}.
Once we have \eqref{S-decay2}, it is obvious that, given $\varepsilon>0$ arbitrarily,
\[
\|V(t)\|_{q,B_R}
\leq C\|V(t)\|_{r_0,B_R}
\leq C\|V(t)\|_{r_0,\mathbb R^3}
\leq Ct^{-3/2+\varepsilon}\|\widetilde F\|_{q_0,\mathbb R^3}
\leq Ct^{-3/2+\varepsilon}\|\widetilde F\|_{q,\mathbb R^3}
\]
by taking $r_0\in (q,\infty)$ large enough and $q_0\in (1,q)$ close enough to $1$ with the aid of \eqref{proj-ang2}.
Set $u(t)=V(t)|_\Omega$.
By the same reasoning as above together with \eqref{ellip-0} and by taking into account \eqref{obje} and \eqref{proj-ang2}, we find
\begin{equation*}
\begin{split}
\|u(t)\|_{W^{2,q}(\Omega_R)}
&\leq C\Big(
\|A{\mathcal S}(t)\mathbb PF\|_{r_0,\mathbb R^3}
+\|{\mathcal S}(t)\mathbb PF\|_{r_0,\mathbb R^3}
\Big)  \\
&\leq C\|e^{-(t-1)A}\mathbb P\widetilde F\|_{r_0,\mathbb R^3} 
  \\
&= C\|{\mathcal S}(t-1)\mathbb P\widetilde F\|_{r_0,\mathbb R^3}  \\
&\leq Ct^{-3/2+\varepsilon}\|\widetilde F\|_{q,\mathbb R^3}
\end{split}
\end{equation*}
for all $t\geq 2$, where $\varepsilon>0$ is arbitrarily small. 
In this way, without relying on spectral analysis,
we can deduce the local energy decay \eqref{1st-led} with the rate $t^{-3/2}$ replaced by less sharp one.
As we have seen in subsection \ref{sub-led2} and section \ref{pr-main},
less sharp rate does not cause any problem to conclude further decay properties \eqref{decay-2}--\eqref{decay-4}
as well as \eqref{decay-1} with $r=\infty$.

\bigskip\bigskip
\noindent
{\bf Acknowledgments}.
The first author is partially supported by the Grant-in-aid for Scientific Research 25K07083 from JSPS.
The second author is partially supported by the National Research Foundation of Korea(NRF) grant funded by the Korea government(MSIT)(RS-2026-2546904).

\bigskip
\noindent
{\bf Declarations}.
The authors state that there is no conflict of interest.
Data sharing not applicable to this article as no datasets were generated or analyzed during the current study.

\begin{flushleft}
Toshiaki Hishida \\
Graduate School of Mathematics \\
Nagoya University \\
Nagoya 464-8602 \\
Japan \\
hishida@math.nagoya-u.ac.jp \\
$\,$ \\
Bum Ja Jin \\
Department of Mathematics Education \\
Mokpo National University \\
Muan 534-729 \\
South Korea \\
bumjajin@mnu.ac.kr
\end{flushleft}

\end{document}